\documentclass{article}
\usepackage{amssymb,amsmath,amsfonts,amsbsy}
\usepackage{mathrsfs}
\usepackage{mathtools}
\usepackage{dsfont}
\usepackage{upgreek}
\usepackage{setspace}
\usepackage{enumitem}
\usepackage{verbatim}
\usepackage[ocgcolorlinks]{hyperref}
 \hypersetup{urlcolor=blue, citecolor=blue,linkcolor=blue}
\usepackage{relsize}
\usepackage[margin=1in]{geometry}
\let\svthefootnote\thefootnote
\newcommand\blankfootnote[1]{%
  \let\thefootnote\relax\footnotetext{#1}%
  \let\thefootnote\svthefootnote%
}
\newcommand{\ssection}[1]{\vspace{0.09in}
\begin{center}\textbf{#1}\end{center}
\vspace{0.09in}}
\newcommand{\bi}{\boldsymbol}
\newcommand{\be}{\begin{equation}}
\newcommand{\ee}{\end{equation}}
\newcommand{\nn}{\nonumber}
\newtheorem{Def}{Definition}[section]
\newcommand{\bd}{\begin{Def}}
\newcommand{\ed}{\end{Def}}
\newtheorem{theorem}{Theorem}
\newcommand{\bt}{\begin{theorem}}
\newcommand{\et}{\end{theorem}}
\newcommand{\bs}{\begin{subequations}}
\newcommand{\es}{\end{subequations}}
\newtheorem{proposition}[theorem]{Proposition}
\newcommand{\bp}{\begin{proposition}}
\newcommand{\ep}{\end{proposition}}
\newtheorem{lemma}[theorem]{Lemma}
\newcommand{\bl}{\begin{lemma}}
\newcommand{\el}{\end{lemma}}
\newtheorem{corollary}[theorem]{Corollary}
\newcommand{\bc}{\begin{corollary}}
\newcommand{\ec}{\end{corollary}}
\newtheorem{remark}[theorem]{\it Remark}
\newcommand{\br}{\begin{remark}}
\newcommand{\er}{\end{remark}}
\newcommand{\norm}[3]{\ensuremath{\big\|#1\big\|^{#2}_{#3}}}
\newcommand{\norms}[3]{\ensuremath{\|#1\|^{#2}_{#3}}}
\newcommand{\vac}{\big|0\big\rangle}

\newcommand{\ex}{\big|\psi_\mathrm{ex}\big\rangle}
\newcommand{\ap}{\big|\psi_\mathrm{ap}\big\rangle}
\newcommand{\red}{\big|\psi_{\mathrm{red}}\big\rangle}

\title{\Large Exact Evolution versus Mean Field with Second-order correction for Bosons Interacting via Short-range Two-body Potential}
\author{Elif Kuz\\{\small University of Maryland, College Park, MD 20742, USA }\\{\small E-mail: elifkuz@math.umd.edu}}
\date{}

\begin{document}

\maketitle
\begin{abstract} We consider the evolution of $N$ bosons, where $N$ is large,  with two-body interactions  of the form $N^{3\beta}v(N^{\beta}\bold{\cdot})$, $0\leq\beta\leq 1$. The parameter $\beta$ measures the strength of interactions. We compare the exact evolution with an approximation which considers the evolution of a mean field coupled with an appropriate description of pair excitations, see \cite{GMM1,GMM2}. For $0\leq \beta < 1/2$, we derive an error bound of the form $p(t)/N^\alpha$, where $\alpha>0$ and $p(t)$ is a polynomial, which implies a specific rate of convergence as $N\rightarrow\infty$.
\end{abstract}

\begin{center}\textbf{1. Introduction}
\end{center}\blankfootnote{2010 MSC: 81V70, 82C10, 35Q40, 35Q55}

The goal of the current work is to study certain aspects of the dynamics of Bose-Einstein condensates that are initially trapped. A Bose-Einstein condensate is a state of matter of a dilute gas of bosons at very low temperatures, in which particles macroscopically occupy the lowest energy state described by a single one particle wave function. This phenomenon was first predicted by Einstein in 1925 for non-interacting massive particles based on the ideas of Bose. The experimental realization of the first condensates was achieved in 1995 \cite{COL, MIT} which has been followed by an increase in the experimental and theoretical activity on the study of the condensates. 

In  experiments, to obtain a condensate, weakly interacting atoms trapped by external potentials are cooled below a certain temperature depending on the density of the gas. Traps are then removed to observe the evolution of the condensate. The properties of interest are the macroscopic properties of the system describing the typical behavior of the particles resulting from averaging over a large number of particles. The limiting behavior as the number of particles go to infinity is expected to be a good approximation for the macroscopic properties  observed in the experiments for a system of large but finite number of particles. 
We can describe the corresponding mathematical model as follows.  We consider a system of $N$ weakly interacting Bosons, the dynamics of which is governed by the $N$-body Schr\"{o}dinger equation
\be \label{first}\frac{1}{i}\partial_t\psi_N=H_N\psi_N\quad\text{with}\quad H_N:=H-N^{-1}V \ee
where
\[H:=\sum_{j=1}^N\Delta_{x_j}\quad\text{and}\quad V:=(1/2)\sum_{j\neq k}N^{3\beta}v(N^\beta(x_j-x_k));\quad 0\leq\beta\leq 1;\quad v\geq 0 \text{ symmetric}\]
both acting on the wave functions 
\[\psi_N\in L^2_{\mathrm{s}}(\mathbb{R}^{3N})\quad\text{with}\quad\|\psi_N\|_{L^2_{\mathrm{s}}(\mathbb{R}^{3N})}=1.\]
$L^2_{\mathrm{s}}(\mathbb{R}^{3N})$ stands for the subspace of $L^2(\mathbb{R}^{3N})$ consisting of symmetric functions in $x_1,x_2,\dots x_N$. The potential $v\in L^1(\mathbb{R}^3)\cap L^\infty(\mathbb{R}^3)$ models two body interactions  and the scaling parameter $\beta$ describes the range and the strength of interactions. As noted in \cite{EESY}, the case of $\beta>1/3$ is more interesting since it represents the self-interaction range. We can explain this as follows. We consider the evolution in \eqref{first} with initial data  coming from the ground state of the following Hamiltonian describing the initially trapped gas: 
\begin{align}\nn H_N^{\mathrm{trap}}:=\sum_{j=1}^N\Big(-\Delta_{x_j}+\underbrace{V_{\mathrm{ext}}(x_j)}_{\mathclap{\substack{\qquad\sim|x_j|^2\text{ typically}\\\qquad\text{in available experiments}}}}\Big)+\frac{1}{2N}\sum_{j\neq k}N^{3\beta}v(N^\beta(x_j-x_k)).\end{align}
The ground state of $H_N^\mathrm{trap}$ looks like a factorized state:
\be\label{2nd} \psi_N(0, x_1,\dots,x_N)\simeq\phi_0^{\otimes N}:= \prod_{j=1}^{N}\phi_0(x_j)\ee
as justified by the works of \cite{LSei,LSY} which showed
\be\gamma_N^{(1)}(0,x,x'):=\int_{\mathbb{R}^{3(N-1)}}\psi_N(0, x,\textbf{x}_{N-1})\overline{\psi_N}(0, x',\textbf{x}_{N-1})\mathrm{d}\textbf{x}_{N-1}\rightarrow \phi_0(x)\bar\phi_0(x')\ee
in trace norm as $N\rightarrow\infty$ and for $\phi_0$ minimizing the Gross-Pitaevskii  energy functional subject to $\|\phi_0\|_{L^2(\mathbb{R}^3)}=1$.
The initial data in \eqref{2nd} is localized in space so we can consider all $N$ particles in a box of unit volume. Then the inter-particle distance is $O(1/N^{1/3})$ compared to  the range of the scaled interaction which is $O(1/N^\beta)$. Hence for $\beta>1/3$, each particle feels only the potential generated by itself, which is consistent with the Gross-Pitaevskii theory proposing to model the many body effects by a strong on-site self interaction.

Concerning the time evolution of an initially factorized (or approximately factorized state), it has been  proved in a series of papers \cite{Hepp,Spohn,GiVe,BGM,EY, EESY,ESY1,ESY2,ESY3,LNS}  under varying assumptions on the interaction and the scaling parameter $\beta$ that the evolution is still approximately factorized  at later times i.e.
\be\label{poc}\psi_N(t,x_1,\dots,x_N)\simeq \prod_{j=1}^{N}\phi(t,x_j)\ee
where the limiting one-particle condensate wave function $\phi$ satisfies the Schr\"{o}dinger equation
\begin{equation}\label{limitdyn}\begin{array}{lcl}\displaystyle{\frac{1}{i}}\partial_t\phi&=&\Delta \phi-g_\beta(|\phi|^2)\phi\\\phi(0,\cdot)&=&\phi_0 \end{array}\end{equation}
with $g_\beta$ being
\[g_\beta(|\phi|^2)=\left\{\begin{array}{lcl}v\ast|\phi|^2,&\quad&\mathrm{\beta=0}\\
(\int v)|\phi|^2,&\quad& 0<\beta<1\\
8\pi a |\phi|^2,&\quad&\beta=1\,\,(a\,\,\mathrm{is\,\,scattering \,\,length\,\,corresp. \,\,to}\,\,v).
\end{array}\right.\]
\cite{Hepp} assumed regularity of pair-potentials and  \cite{GiVe} extended the results to singular potentials and both of the mentioned works dealt with the problem in the framework of the second quantization.
In \cite{BGM, EY} for $\beta=0$ and Coulomb potential and in \cite{ESY1, ESY2, ESY3} for values of $\beta$ different than zero, with a strategy based on the work of \cite{Spohn}, \eqref{poc} was established in the sense of marginals i.e.
\be\label{cgamma}\gamma_N^{(1)}(t,x,x'):=\int_{\mathbb{R}^{3(N-1)}}\psi_N(t, x,\textbf{x}_{N-1})\overline{\psi_N}(t, x',\textbf{x}_{N-1})\mathrm{d}\textbf{x}_{N-1}\rightarrow \phi(t,x)\bar\phi(t,x')\ee
in trace norm as $N\rightarrow \infty$ and similarly for the higher order marginals $\gamma_N^{(k)}$ where $k$ is fixed. See \cite{CP,CH,CH2} for an approach based on the space-time norms introduced in \cite{KM} and also \cite{LNS} for an approximation of the exact dynamics $\psi_N(t,\cdot)$ in $L^2(\mathbb{R}^{3N})$ in the Hartree case ($\beta=0$). 

One would also like to obtain results on the rate of convergence in \eqref{poc}. Using the framework of the second quantization in \cite{RS, CLS, BDS}, the following results are obtained.
\begin{align} \label{rocBDRS}\mathrm{Trace}\,\big|\gamma_N^{(1)}(t,\cdot)-\phi(t,\cdot)\otimes\overline\phi(t,\cdot)\big|\lesssim\left\{\begin{array}{ll}e^{Ct}/N  &\scriptsize\text{by \cite{CLS}, }\text{for }\scriptstyle{\beta=0}\text{ and }\\&\scriptsize\text{singular potentials}\\
&\scriptsize\text{including }\scriptstyle{v(x)=|x|^{-1},}\\
\exp(c_1 \exp(c_2t))/N^{1/4} &\scriptstyle\text{by \cite{BDS},}\text{ for }\beta=1\text{ and}\\&\scriptstyle
v\in L^1\cap L^3(\mathbb{R}^3, (1+|x|^6)\mathrm{d}x).\end{array}\right.\end{align}

Convergence in the sense  of marginals provides with partial information about the system since most of the variables are averaged out. Also, although the Hartree equation (corresp. to $\beta=0$) or  the cubic nonlinear Schr\"{o}dinger equation provide a good description of the limiting behavior for the mean field represented by the condensate wave function, they fail to describe pair excitations i.e. the scattering of particles in pairs from the condensate to other states. 
Hence, in \cite{GMM1, GMM2}, inspired by but being different than  that of  \cite{Wu}, a Fock space approximation of the exact dynamics which considers pair excitations as a correction to the mean field has been introduced. Error bounds  deteriorating more slowly in time compared to the bounds in \eqref{rocBDRS} have been obtained. Those results are valid for $\beta<1/3$. Our approach in this note is more along the lines 
of \cite{GMM1,GMM2,GM} and we aim  to extend their results to higher $\beta$ values by using an elliptic estimate, the end-point Strichartz estimates, energy methods and an iteration scheme. 
We will outline the approach of \cite{GMM1, GMM2, GM}, describe our general strategy and provide with the organization of the current note in the next section. 

\begin{center}\textbf{2. Earlier Results and General Strategy}
\end{center}

We will deal with the problem described in the introduction by embedding the $N$-body dynamics in the Fock space
which is defined as:
\[\mathbb{F}=\bigoplus_{n=0}^\infty \mathbb{F}_n; \quad\mathbb{F}_n=L_\mathrm{s}^2(\mathbb{R}^{3n})\,\,\mathrm{for}\,\,n>0\quad\mathrm{and}\quad \mathbb{F}_0=\mathbb{C}\]
containing vectors of the form
\[|\psi\rangle= (\psi_0,\psi_1,\psi_2,\dots).\]
The norm is defined as
\[\||\psi\rangle\|^2= |\psi_0|^2+\sum_{n\geq 1}\|\psi_n\|^2_{L^2(\mathbb{R}^{3n})}.\]  
We will use the following notation for what is known as the vacuum state:
$$\vac=(1,0,0,\dots).$$

We have the following annihilation and creation operators acting on the sectors of Fock vectors as:
\begin{subequations}\begin{align}&\displaystyle
a_x(\psi_{n+1})=\sqrt {n+1}\psi_{n+1}(x,x_1, \dots,x_n),\\
&a^{\ast}_x(\psi_{n-1})={\displaystyle\frac{1}{\sqrt n}\sum_{j=1}^{n}}\delta(x-x_j)\psi_{n-1}(x_1,\dots,x_{j-1},x_{j+1},\dots,x_n). \end{align}\end{subequations}
which satisfy
$\label{8c}[a_x,a_y^{\ast}]=\delta(x-y)\mathrm{\,\,and\,\,}a_x\vac=0.$

In order to embed the $N$-body system  in the Fock space we need to define the Fock Hamiltonian acting on $\mathbb{F}$ as
\begin{subequations}
\begin{align}
&\label{originalham}\mathcal{H}:=\mathcal{H}_1-N^{-1}\mathcal{V}\quad\mathrm{where,}\\
&\mathcal{H}_1=\int\Delta_x\delta(x-y)a_{x}^{\ast} a_y\,\mathrm{d}x\,\mathrm{d}y,\\
&\label{11c} \mathcal{V}:=\frac{1}{2}\int v_N(x-y)a_x^{\ast}a_y^{\ast}a_{x}a_{y}\,\mathrm{d}x\,\mathrm{d}y\\
&\nonumber\mathrm{with}\,\,\,v_N(x):=N^{3\beta}v(N^\beta x)
\end{align}
\end{subequations}
Before introducing the initial value problem in the Fock space we need to define the coherent  states which we will use as our initial data. First let's define the operators
\[a(\bar\phi):=\int\mathrm{d}x\, \bar\phi(x)a_x\quad\mathrm{and}\quad a^\ast(\phi):=\int\mathrm{d}x\,\phi(x)a^\ast_x\quad\text{for}\,\,\phi\in L^2({\mathbb{R}^3})\]
via which we define
\begin{equation}
\mathcal{A}(\phi):=a(\bar\phi)-a^\ast(\phi).
\end{equation}
The problem in the Fock space can be written as:
\begin{subequations}\begin{align}
&\label{Fockeqn}\frac{1}{i}\partial_t|\psi\rangle=\mathcal{H}|\psi\rangle,\\
&\label{FockID}|\psi(0)\rangle=e^{-\sqrt N\mathcal{A}(\phi_0)}\vac=\bigg(\dots, c_n\prod_{j=1}^n\phi_0(x_j),\dots
\bigg)\,\,\mathrm{with}\,\,c_n=(e^{-N}N^n/n!)^{1/2}
\end{align}\end{subequations}
\noindent so that in the $N$-th sector we have the $N$-body equation \eqref{first} with the initial data $c_N\phi_0^{\otimes N}$ where $c_N\simeq(2\pi N)^{-1/4}$. The operator $e^{-\sqrt N\mathcal{A}(\phi_0)}$ in \eqref{FockID} is called the Weyl operator and the initial data defined using it, which has a tensor product in each sector, is called a coherent state. 

The solution to the initial value problem \eqref{Fockeqn}-\eqref{FockID} can  be written formally as
\begin{equation}
\label{ex}\ex:=e^{it\mathcal{H}}e^{-\sqrt{N}\mathcal{A}(\phi_0)}\vac\quad\mathrm{(exact\,\,evolution)}
\end{equation}
and  the mean field evolution is described by $e^{-\sqrt N\mathcal{A}(\phi)}\vac$ where $\phi$ satisfies
\be\label{cf}\frac{1}{i}\partial_t\phi-\Delta\phi+\big(v_N\ast|\phi|^2\big)\phi=0\text{\,\,with\,\,}\phi(0,\cdot)=\phi_0.\ee
 However the mean field evolution does not track the exact dynamics in the Fock space norm. Hence, in \cite{GMM1,GMM2}, via the operator
\begin{equation}
\mathcal{B}(k):=\frac{1}{2}\int \big\{\overline{k}(x,y)a_x a_y-k(x,y)a_x^{\ast} a_y^{\ast}\big \}\,\mathrm{d}x\,\mathrm{d}y,
\end{equation}
a second order correction  was introduced, namely, a state of the form
\begin{equation}
\label{ap}\ap:=e^{iN\chi(t)}e^{-\sqrt{N}\mathcal{A}(\phi)}e^{-\mathcal{B}(k)}\vac\quad\mathrm{(approximate\,\,evolution)}
\end{equation}
where $k(t,x,y)$ is a  function describing the pair excitations and whose evolution is to be determined from the dynamics of \eqref{Fockeqn}. The explicit form of the phase  $\chi(t)$ is irrelevant to the present discussion and we refer to \cite{GMM1,GMM2,GM} for more details. The evolution of pairs is expressed through 
\begin{subequations}\begin{align}
&\label{sh}\mathrm{sh}(k):=k+\frac{1}{3!}k\circ\bar k\circ k+\dots,\\
&\label{delta+p}\mathrm{ch}(k):=\delta(x-y)+p=\delta(x-y)+\frac{1}{2!}\bar k\circ k+\dots\end{align}\end{subequations}
where the products above are of the operator type i.e.
\[(k\circ l)(x,y):=\int k(x,z)l(z,y)\mathrm{d}z\]
for $k$ and $l$ symmetric Hilbert-Schmidt operators on $L^2(\mathbb{R}^3)$.  To describe the evolution of pairs we also need to define 
\begin{subequations}\begin{align}
&\label{112a}g(t,x,y)=\big(-\Delta_x+(v_N\ast|\phi|^2)(t,x)\big)\delta(x-y)+v_N(x-y)\overline{\phi}(t,x)\phi(t,y),\\
&\label{ohss}m(t,x,y)=-v_N(x-y)\phi(t,x)\phi(t,y)\end{align}\end{subequations}
which will help us introduce the operators
\begin{subequations}\label{GMeqns}\begin{align}
&\label{bs}\bi{\mathrm{S}}(s):=\frac{1}{i}\partial_t s+g^T\circ s+s\circ g\text{\,\,\,\,(Schr\"{o}dinger-type)},\\
&\label{dbly}\bi{\mathrm{W}}(p):=\frac{1}{i}\partial_t p+[g^T,p]\text{\,\,\,\,(Wigner-type)}.\end{align}\end{subequations}
We can now state the system of equations (see \cite{GM}) satisfied by the pair excitations function $k$ which we see in \eqref{ap}:
\begin{subequations}\label{efk}
\begin{align}
&\label{efka}\bi{\mathrm{S}}(\mathrm{sh}(2k))=m\circ\mathrm{ch}(2k)+\overline{\mathrm{ch}(2k)}\circ m,\\
&\label{efkb}\bi{\mathrm{W}}(\overline{\mathrm{ch}(2k)})=m\circ\overline{\mathrm{sh}(2k)}-\mathrm{sh}(2k)\circ \bar m,\\
&\nonumber\mathrm{with}\,\, k(0,\cdot)=0.
\end{align}
\end{subequations}
The assumption $k(0,\cdot)\equiv 0$ is to consider only the coherent initial states. However one can consider more general initial data of the form $e^{-\sqrt N\mathcal{A}(\phi_0)}e^{-\mathcal{B}(k_0)}\vac$ with $k_0=k(0,\cdot)\neq 0$ (see Remark 3, \cite{GMM1}).

The main result in \cite{GM}, building on the results in \cite{GMM1, GMM2} was the following:
\bt Let $\phi$ satisfy \eqref{cf} with $\phi_0\in W^{m,1}(\mathbb{R}^3)$ ($m$ derivatives in $L^1$) with $m\geq 2$  and  $k$ satisfy \eqref{efka}-\eqref{efkb}. Then, recalling the definitions \eqref{ex} and \eqref{ap},
\be\label{GMest}\big\|\ex-\ap\big\|_\mathbb{F}\lesssim \frac{(1+t)\log^4(1+t)}{N^{(1-3\beta)/2}}\ee
provided $0<\beta<1/3$.
\et

Our main result in this paper is extending the error estimates to the case of $0<\beta<1/2$ in the following way:
\bt\label{mainresult} Let $\phi$ satisfy \eqref{cf} with $\phi(0,\cdot)\in W^{m,1}(\mathbb{R}^3)$ for $m\geq 6$  and $k$ satisfy \eqref{efka}-\eqref{efkb} with some regularity assumptions on $\big(\partial_t\mathrm{sh}(2k)\big)(0,\cdot)$ to be specified later in the proof. Then, for any $\epsilon>0$ and $j$ positive integer,
\be\label{mrest}\big\|\ex-\ap\big\|_\mathbb{F} \lesssim_{\epsilon,j} t^{\frac{j+3}{2}}\log^6(1+t)\,\mathlarger{\bi\cdot}\left\{\begin{array}{ll} N^{-\frac{1}{2}+\beta(1+\epsilon)} & \mathsmaller{\mathrm{for}}\,\, \frac{1}{3}\leq\mathsmaller\beta < \frac{2j}{(1-2\epsilon +4j)},\\
N^{\frac{-3+7\beta}{2}+(j-1)(-1+2\beta)}&\mathsmaller{\mathrm{for}}\,\,\frac{1+2j}{3+4j}>\mathsmaller\beta \geq\frac{2j}{(1-2\epsilon +4j)}.
\end{array}\right.
\ee
The above estimate implies  a decay as $N\rightarrow\infty$ for $\beta$ as close as desired to 1/2 if we choose $\epsilon$ sufficiently small and $j$ sufficiently large.
\et
\noindent\textbf{Remarks.}
\textit{\begin{itemize}[leftmargin=.55cm]
\item[(i)] \cite{GM15} extended the estimates to the case of $\beta<2/3$ by considering a coupled system introduced in \cite{GM13}, instead of the uncoupled one consisting of \eqref{cf} and  \eqref{efk} that is considered in the current work. Also, similar Fock sapce estimates have recently been obtained in \cite{BCS} for $\beta\in(0,1)$ using a certain class of initial data and an explicit choice of pair excitation function $k$. However the dependence of the error bounds in \cite{BCS} on time is exponential.
\item[(ii)] The above estimate also implies a rate of convergence in the sense of marginals in trace norm  which deteriorates more slowly in time compared to the bounds in \eqref{rocBDRS} since 
\begin{align} \nn &\mathrm{Trace}\,\big|\gamma_N^{(1)}(t,\cdot)-\phi^{(N)}(t,\cdot)\otimes\bar\phi^{(N)}(t,\cdot)\big|\lesssim N^{1/4}\Big(1+\|\mathrm{sh}(k)\|_{L^2(\mathbb{R}^6)}^4\Big)\big\|\ex-\ap\big\|_\mathbb{F}\\
&\mathrm{Trace}\label{ctld}\,\big|\phi^{(N)}(t,\cdot)\otimes\bar\phi^{(N)}(t,\cdot)-\phi(t,\cdot)\otimes\bar\phi(t,\cdot)\big|\lesssim \|\phi^{(N)}-\phi\|_{L^2(\mathbb{R}^3)}\lesssim 1/\sqrt{N} \end{align}
as \hspace{.1cm}proved\hspace{.1cm} in\hspace{.1cm} \cite{EK},\hspace{.1cm} partly\hspace{.1cm} based\hspace{.1cm} on\hspace{.1cm} some\hspace{.1cm} arguments\hspace{.1cm} in \cite{BDS}.\hspace{.1cm} In\hspace{.1cm} the\hspace{.1cm} above\hspace{.1cm} lines,\hspace{.1cm} $\quad\gamma_N^{(1)}(t,x,x')=\\\int_{\mathbb{R}^{3(N-1)}}\psi_N(t, x,\textbf{x}_{N-1})\overline{\psi_N}(t, x',\textbf{x}_{N-1})\mathrm{d}\textbf{x}_{N-1}$ and $\phi^{(N)}$ and $\phi$ solve the equations \eqref{cf} and \eqref{limitdyn} respectively.\end{itemize}}

The proof of Theorem \ref{mainresult} is based on estimating the deviation of the evolution from the vacuum state defined as
\begin{subequations}\begin{align}
&\label{saytilde}|\tilde{\psi}\rangle:=e^{-iN\chi(t)}\red-\vac\text{\,\,\,\,where\,\,}\\
&\label{red}\red:=e^{\mathcal{B}(k)}e^{\sqrt{N}\mathcal{A}(\phi)}\ex\quad\mathrm{(reduced\,\,dynamics)}
\end{align}
\end{subequations}
which satisfies
$$\big\||\tilde{\psi}\rangle\big\|_\mathbb{F}=\big\||\psi_\mathrm{ex}\rangle-|\psi_\mathrm{ap}\rangle\big\|_\mathbb{F}$$
due to $e^{-\sqrt N \mathcal{A}}$ and $e^{-\mathcal B}$ being unitary. We can obtain the evolution for $|\tilde{\psi}\rangle$ as follows. A straightforward computation gives  the evolution of the reduced dynamics:
\be \label{redev}\frac{1}{i}\partial_t\red=\mathcal{H}_\mathrm{red}\red\ee
where
\be\mathcal{H}_\mathrm{red}:=\frac{1}{i}\big(\partial_t e^{\mathcal B}\big)e^{-\mathcal B}+e^{\mathcal B}\left(\frac{1}{i}\big(\partial_te^{\sqrt N \mathcal{A}}\big)e^{-\sqrt N \mathcal{A}}+e^{\sqrt N \mathcal{A}}\mathcal{H}e^{-\sqrt N \mathcal{A}}\right)e^{-\mathcal B}.\ee
As shown in section 2 of \cite{GM}, if \eqref{cf} and \eqref{efk} hold, then
\be\label{twf}\mathcal{H}_\mathrm{red}=N\mu(t)+\int \mathrm{d}x\mathrm{d}y \left\{L(t,x,y)a_x^\ast a_y\right\}-N^{-1/2}\mathcal{E}(t)\ee
which together with \eqref{saytilde}, \eqref{redev} and the fact that $a\vac=0$ implies
\be\boxed{\hspace{1.5cm} \label{oe}\bigg(\frac{1}{i}\partial_t-\mathop{\mathcal{L}}_{\mathclap{\substack {\\{\mathlarger{\swarrow}}\,\,\,\,\, \\ {\scriptsize\overbrace{\int L(t,x,y)a_x^\ast a_y\mathrm{d}x\,\mathrm{d}y-N^{-1/2}\mathcal{E}(t)}}\qquad}}}\bigg)|\mathop{\tilde{\psi\rangle}}
=\underbrace{-N^{-1/2}\mathcal{E}(t)\vac}_{\mathclap{\substack{\quad\searrow \\\qquad\qquad\qquad\quad\quad\quad\quad=:(0,F_1,F_2,F_3,F_4,0,\dots)\\\hspace{4cm}\text{since }\mathcal{E}(t)\text{ is  4th order in }(a,a^\ast)\\\hspace{4cm}\text{as to be explained shortly}}}}\text{ with }|\tilde{\psi}(0)\rangle=0\qquad}\ee
where $\mu(t)$ in \eqref{twf} is related to the phase $\chi(t)$ in \eqref{saytilde} via $\chi(t)=\int \mu(t)$.
 The integral term in \eqref{twf} is the second quantization of the self-adjoint one-particle operator $L(t,x,y)$ which can be considered to be the sum of some kinetic and ``potential" parts as follows:
\be\label{ltxy}\begin{split}L(t,x,y):=&\overbrace{\Delta_x\delta(x-y)-\big(v_N\ast |\phi|^2\big)(t,x)\delta(x-y)-v_N(x-y) \phi(t,x)\bar\phi(t,y)}^{-g(t,y,x)\text{ as defined in \eqref{112a}}}\\&
+\frac{1}{2}\Big((\mathop{\bar c}_{\mathclap{\substack{\\ \mathlarger{\downarrow} \\\overline{ \mathrm{ch}(k)}}}}) ^{-1}\circ \mathop{m}_{\mathclap{\substack{\\\mathlarger{\mathlarger{\mathlarger{\Big\downarrow}}}\\-v_N(x-y)\phi(t,x)\phi(t,y)\\\text{from \eqref{ohss}}}}}\circ \bar u+\mathop{u}_{\mathclap{\substack{\\ \mathlarger{\downarrow} \\ \mathrm{sh}(k)}}}\circ\bar m\circ(\bar c) ^{-1}+\big[\underbrace{\bi{\mathrm W}(\bar c)}_{\mathclap{\substack{\\\mathlarger{\downarrow}\\\frac{1}{i}\partial_t\bar c+[g^T,\bar c]\\ \text{see top line for }g}}},(\bar c) ^{-1}\big]\Big).\end{split}\ee

$N^{-1/2}\mathcal{E}(t)$ in \eqref{twf} is an  error term containing polynomials in $(a, a^\ast)$ up to degree four. It is a self-adjoint operator which can be  written in increasing order in terms of the degrees of polynomials  in the following way where sub-indices  show the degree of the corresponding contribution and superscript ``sa" means ``terms involved are self-adjoint"\footnote{See Section 5, \cite{GM} for the computations leading to this explicit form of $N^{-1/2}\mathcal{E}(t)$.}:
\begin{equation}\label{ese} N^{-1/2}\mathcal{E}(t)=\mathcal{E}_1(t)+\mathop{\mathcal{E}_1^\ast(t)}^{ \mathclap{\substack{\text{denotes}\\\text{adjoint}\\\downarrow\,\,\,\,}}}+\mathcal{E}_2(t)+\mathcal{E}_2^\ast(t)+\mathcal{E}_2^{\mathrm{sa}}(t)+\mathcal{E}_3(t)+\mathcal{E}_3^\ast(t)+\mathcal{E}_4(t)+\mathcal{E}_4^\ast(t)+\mathcal{E}_4^{\mathrm{sa}}(t)\end{equation}
and we define $j$th order terms $\mathcal{E}_j(t)$ and $\mathcal{E}_j^{\mathrm{sa}}(t)$ using the notation $\mathcal{D}_{xy}:=a_x^\ast a_y$, $\mathcal{Q}_{xy}^\ast:=a_x^\ast a_y^\ast$, $\mathcal{Q}_{xy}:=a_x a_y$ and suppressing  the time dependence of the functions $\phi$, $c:=\mathrm{ch}(k)=\delta(x-y)+p$ and $u:=\mathrm{sh}(k)$ in the following list of terms, which needs to be given here explicitly for future reference. We will estimate them in various ways to  be explained later.
\begin{subequations}\label{wad}\begin{align} &\label{22a}\mathcal{E}_1(t):=N^{-1/2}\int \mathrm{d}x_1\mathrm{d}x_2\mathrm{d}y_1\Big\{\big( u\circ c\big)(x_1,x_2)v_N(x_1-x_2)\bar\phi(x_2)\bar u(y_1,x_1)a_{y_1}\\
&\label{22b}\hspace{4.58cm}+c(y_1,x_1)v_N(x_1-x_2)\phi(x_2)\big(c\circ \bar u\big)(x_1,x_2)a_{y_1}\Big\}\\
&\label{22c}\mathcal{E}_2(t):=\frac{1}{2N}\int \mathrm{d}x_1\mathrm{d}x_2\mathrm{d}y_1\mathrm{d}y_2\Big\{\big(\bar u\circ \bar c\big)(x_1,x_2)v_N(x_1-x_2)c(y_1,x_1)u(x_2,y_2)\mathcal{D}_{y_2y_1}\\
&\label{22d}\hspace{4.68cm}+\big(\bar u\circ\bar c\big)(x_1,x_2)v_N(x_1-x_2)u(y_1,x_1)\bar c(x_2,y_2)\mathcal{D}_{y_1y_2}\\
&\label{22e}\hspace{4.68cm}+\big(\bar u\circ\bar c\big)(x_1,x_2)v_N(x_1-x_2)c(y_1,x_1)\bar c(x_2,y_2)\mathcal{Q}_{y_1y_2}\\
&\label{22f}\hspace{4.68cm}+\big(u\circ c\big)(x_1,x_2)v_N(x_1-x_2)\bar u(y_1,x_1)\bar u(x_2,y_2)\mathcal{Q}_{y_1y_2}\Big\}\displaybreak\\
&\label{22g}\mathop{\mathcal{E}^{\mathrm{sa}}_2}^{\mathclap{\smash{\substack{\quad\qquad\text{terms involved}\\\quad\qquad\text{are self-adjoint}\\\,\,\,\,\downarrow\\\color{white} 34\\\color{white}34\\\color{white}34}}}}(t):=\frac{1}{2N}\int \mathrm{d}x_1\mathrm{d}x_2\mathrm{d}y_1\mathrm{d}y_2\Big\{\big(u\circ\bar u)(x_1,x_2)v_N(x_1-x_2)\bar u(y_1,x_1)u(x_2,y_2)\mathcal{D}_{y_2y_1}\\
&\label{22h}\hspace{4.90cm}+2\big(u\circ\bar u)(x_1,x_1)v_N(x_1-x_2)\bar u(y_2,x_2)u(y_1,x_2)\mathcal{D}_{y_1y_2}\Big\}\\
&\label{22i}\mathcal{E}_3(t):=N^{-1/2}\int \mathrm{d}x_1\mathrm{d}x_2\mathrm{d}y_1\mathrm{d}y_2\mathrm{d}y_3\Big\{\bar u(y_1,x_1)v_N(x_1-x_2)\phi(x_2)c(x_2,y_2)c(y_3,x_1)\mathcal{D}_{y_2y_1}a_{y_3}\\
&\label{22j}\hspace{3cm}+\bar c(y_1,x_1)v_N(x_1-x_2)\phi(x_2)\bar u(x_2,y_2)c(y_3,x_1)\mathcal{D}_{y_1y_2}a_{y_3}\\
&\label{22k}\hspace{3cm}+\bar c(y_1,x_1)v_N(x_1-x_2)\bar \phi(x_2)c(y_2,x_1)\bar c(x_2,y_3)a^\ast_{y_1}\mathcal{Q}_{y_2y_3}\\
&\label{22l}\hspace{3cm}+\bar u(y_1,x_1)v_N(x_1-x_2)\phi(x_2)\bar u(x_2,y_2)c(y_3,x_1)\mathcal{Q}_{y_1y_2}a_{y_3}\\
&\label{22m}\hspace{3cm}+\bar u(y_1,x_1)v_N(x_1-x_2)\bar \phi(x_2)u(y_2,x_1)\bar c(x_2,y_3)a_{y_1}\mathcal{D}_{y_2y_3}\\
&\label{22n}\hspace{3cm}+\bar u(y_1,x_1)v_N(x_1-x_2)\bar \phi(x_2)c(y_2,x_1)u(x_2,y_3)a_{y_1}\mathcal{D}_{y_3y_2}\\
&\label{22o}\hspace{3cm}+\bar u(y_1,x_1)v_N(x_1-x_2)\bar \phi(x_2)c(y_2,x_1)\bar c(x_2,y_3)a_{y_1}\mathcal{Q}_{y_2y_3}\\
&\label{22p}\hspace{3cm}+\bar u(y_1,x_1)v_N(x_1-x_2)\phi(x_2)\bar u(x_2,y_2)u(y_3,x_1)\mathcal{Q}_{y_1y_2}a^\ast_{y_3}\Big\}\\
&\nn\mathcal{E}_4(t):=\\
&\label{22q}\frac{1}{2N}\int \mathrm{d}x_1\mathrm{d}x_2\mathrm{d}y_1\mathrm{d}y_2\mathrm{d}y_3\mathrm{d}y_4\Big\{\bar u(y_1,x_1)c(x_2,y_2)v_N(x_1-x_2)c(y_3,x_1)u(x_2,y_4)\mathcal{D}_{y_2y_1}\mathcal{D}_{y_4y_3}\\
&\label{22r}\hspace{3cm}+\bar c(y_1,x_1)\bar u(x_2,y_2)v_N(x_1-x_2)c(y_3,x_1)\bar c(x_2,y_4)\mathcal{D}_{y_1y_2}\mathcal{Q}_{y_3y_4}\\
&\label{22s}\hspace{3cm}+\bar u(y_1,x_1)\bar u(x_2,y_2)v_N(x_1-x_2)c(y_3,x_1)u(x_2,y_4)\mathcal{Q}_{y_1y_2}\mathcal{D}_{y_4y_3}\\
&\label{22t}\hspace{3cm}+\bar u(y_1,x_1)c(x_2,y_2)v_N(x_1-x_2)c(y_3,x_1)\bar c(x_2,y_4)\mathcal{D}_{y_2y_1}\mathcal{Q}_{y_3y_4}\\
&\label{22u}\hspace{3cm}+\bar u(y_1,x_1)\bar u(x_2,y_2)v_N(x_1-x_2)u(y_3,x_1)\bar c(x_2,y_4)\mathcal{Q}_{y_1y_2}\mathcal{D}_{y_3y_4}\\
&\label{22v}\hspace{3cm}+\bar u(y_1,x_1)\bar u(x_2,y_2)v_N(x_1-x_2)c(y_3,x_1)\bar c(x_2,y_4)\mathcal{Q}_{y_1y_2}\mathcal{Q}_{y_3y_4}\Big\}\\
&\nn\mathcal{E}^{\mathrm{sa}}_4(t):=\\
&\label{22w}\frac{1}{2N}\int \mathrm{d}x_1\mathrm{d}x_2\mathrm{d}y_1\mathrm{d}y_2\mathrm{d}y_3\mathrm{d}y_4\Big\{\bar c(y_1,x_1)\bar u(x_2,y_2)v_N(x_1-x_2)c(y_3,x_1)u(x_2,y_4)\mathcal{D}_{y_1y_2}\mathcal{D}_{y_4y_3}\\
&\label{22x}\hspace{3cm}+\bar u(y_1,x_1)c(x_2,y_2)v_N(x_1-x_2)u(y_3,x_1)\bar c(x_2,y_4)\mathcal{D}_{y_2y_1}\mathcal{D}_{y_3y_4}\\
&\label{22y}\hspace{3cm}+\bar c(y_1,x_1)c(x_2,y_2)v_N(x_1-x_2)c(y_3,x_1)\bar c(x_2,y_4)\mathcal{Q}^\ast_{y_1y_2}\mathcal{Q}_{y_3y_4}\\
&\label{22z}\hspace{3cm}+\bar u(y_1,x_1)\bar u(x_2,y_2)v_N(x_1-x_2)u(y_3,x_1)u(x_2,y_4)\mathcal{D}_{y_3y_1}\mathcal{D}_{y_4y_2}\Big\}.
\end{align}\end{subequations}

Based on the explicit form of $N^{-1/2}\mathcal{E}(t)$ in \eqref{ese} and recalling $c:=ch(k)=\delta(x-y)+p$, the sectors of the forcing term $-N^{-1/2}\mathcal{E}(t)\vac$ in \eqref{oe} can be computed (up to symmetrization in the 2nd, 3rd and 4th sectors) as\footnote{The main idea of this computation is to commute $a$ (if there is any), to the right hand side, with $a^\ast$ operators in those terms in \eqref{ese} which do not annihilate the vacuum. This produces some lower order terms (contributions of which we see in \eqref{F1a}-\eqref{F4d}) and terms which annihilate $\vac$ since $a\vac=0$. See Section 5, \cite{GM} for the details.}: \newline

$\bullet$ Sector $\mathbb{F}_1$:
\begin{subequations}\label{F1}\begin{align}
F_1(t,y_1):=-N^{-1/2}\Bigg(&\int \mathrm{d}x_1\mathrm{d}x_2v_N(x_1-x_2)\Big\{\label{F1a}u(y_1,x_2)(\bar u\circ u)(x_1,x_1)\bar\phi(x_2)\\
&\,\,\quad\qquad\qquad\qquad\qquad\quad \label{F1b}+\bar p(y_1,x_2)(u\circ\bar u)(x_1,x_1)\phi(x_2)\\
&\,\,\quad\qquad\qquad\qquad\qquad\quad \label{F1c}+u(y_1,x_1)(\bar u\circ u)(x_1,x_2)\bar\phi(x_2)\\
&\,\,\quad\qquad\qquad\qquad\qquad\quad \label{F1d}+\bar p(y_1,x_1)(\bar p\circ u)(x_1,x_2)\bar \phi(x_2)\\
&\,\,\quad\qquad\qquad\qquad\qquad\quad \label{F1e}+\bar p(y_1,x_1)(u\circ\bar u)(x_1,x_2)\phi(x_2)\\
&\,\,\quad\qquad\qquad\qquad\qquad\quad \label{F1f}+u(y_1,x_1)(\bar u\circ\bar p)(x_1,x_2)\phi(x_2)\\
&\,\,\quad\qquad\qquad\qquad\qquad\quad \label{F1g}+\bar p(y_1,x_1)u(x_1,x_2)\bar \phi(x_2)\\
&\,\,\quad\qquad\qquad\qquad\qquad\quad \label{F1h}+u(y_1,x_1)\bar u(x_1,x_2)\phi(x_2)
\Big\}\displaybreak\\
&\label{F1i}+\int\mathrm{d}x_1\, v_N(y_1-x_1)\Big\{u(y_1,x_1)\bar\phi(x_1)\\
&\label{F1j}\qquad\qquad\qquad\qquad\qquad+ (u\circ\bar u)(y_1,x_1)\phi(x_1)\\
&\label{F1k}\qquad\qquad\qquad\qquad\qquad+(\bar p\circ u)(y_1,x_1)\bar\phi(x_1)\\
&\label{F1l}\qquad\qquad\qquad\qquad\qquad+(u\circ\bar u)(x_1,x_1)\phi(y_1)\Big\}\Bigg)
\end{align}\end{subequations}

$\bullet$ Sector $\mathbb{F}_2$:
\begin{subequations}\label{F2}\begin{align}
F_2(t,y_1,y_2):=-\frac{1}{2N}\Bigg(&\label{F2a}v_N(y_1-y_2)\Big\{u(y_1,y_2)+(\bar p\circ u)(y_1,y_2)\Big\}\\ \quad &+\int \mathrm{d}x_1\mathrm{d}x_2 v_N(x_1-x_2)\Big\{2\label{F2b}\bar p(y_1,x_2)u(x_2,y_2)(\bar u\circ u)(x_1,x_1)\\
&\,\,\quad\qquad\qquad\qquad\qquad\qquad \label{F2c}+2\bar p(y_1,x_2)u(x_1,y_2)(\bar u\circ u)(x_1,x_2)\\
&\,\,\quad\qquad\qquad\qquad\qquad\qquad \label{F2d}+u(y_1,x_1)u(x_2,y_2)(\bar u\circ\bar p)(x_1,x_2)\\
&\,\,\quad\qquad\qquad\qquad\qquad\qquad \label{F2e}+\bar p(y_1,x_1)p(x_2,y_2)(\bar p\circ u)(x_1,x_2)\\
&\,\,\quad\qquad\qquad\qquad\qquad\qquad \label{F2f}+u(y_1,x_1) u(x_2,y_2)\bar u(x_1,x_2)\\
&\,\,\quad\qquad\qquad\qquad\qquad\qquad \label{F2g}+\bar p (y_1,x_1)p(x_2,y_2)u(x_1,x_2)\Big\}\\
&\label{F2h}+\int\mathrm{d}x_1v_N(y_1-x_1) \Big\{2u(y_1,y_2)(\bar u\circ u)(x_1,x_1)\\
&\,\,\,\,\,\quad\qquad\qquad\qquad\qquad \label{F2i}+\bar p(y_2,x_1)u(x_1,y_1)\\
&\,\,\,\,\,\quad\qquad\qquad\qquad\qquad \label{F2j}+2u(x_1,y_2)(\bar u\circ u)(x_1,y_1)\\
&\,\,\,\,\,\quad\qquad\qquad\qquad\qquad \label{F2k}+\bar p(y_2,x_1)(\bar p \circ u)(y_1,x_1)
\Big\}\\
&\label{F2l}+\int \mathrm{d}x_1\,v_N(x_1-y_2)\bar p(y_1,x_1)(\bar c\circ u)(x_1,y_2)\Bigg)
\end{align}\end{subequations}

$\bullet$ Sector $\mathbb{F}_3$:
\begin{subequations}\label{F3}\begin{align}
\label{F3a}F_3(t,y_1,y_2,y_3):=-N^{-1/2}\Big\{&v_N(y_1-y_2)\phi(y_2)u(y_3,y_1)\\
&\label{F3b}+\int \mathrm{d}x\{v_N(y_1-x)\bar\phi(x)u(x,y_3)\}u(y_2,y_1)\\
&\label{F3c}+\int \mathrm{d}x\{\bar p(y_1,x)v_N(x-y_2)u(y_3,x)\}\phi(y_2)\\
&\label{F3d}+\int \mathrm{d}x\{\bar p(y_2,x)v_N(y_1-x)\phi(x)\}u(y_3,y_1)\\
&\label{F3e}+\int\mathrm{d}x_1\mathrm{d}x_2\{\bar p(y_1,x_1)v_N(x_1-x_2)\bar \phi(x_2)u(y_2,x_1)u(x_2,y_3)\}\\
&\label{F3f}+\int\mathrm{d}x_1\mathrm{d}x_2\{\bar p(y_1,x_1)p(x_2,y_2)v_N(x_1-x_2)\phi(x_2)u(y_3,x_1)\}\Big\}
 \end{align}\end{subequations}

$\bullet$ Sector $\mathbb{F}_4$:
\begin{subequations}\label{F4}\begin{align}
\label{F4a}F_4(t,y_1,y_2,y_3,y_4):=&-(1/2N)\Big\{v_N(y_1-y_2)u(y_3,y_1)u(y_2,y_4)\\
&\label{F4b}+\int \mathrm{d}x \{\bar p(y_2,x)v_N(y_1-x)u(x,y_4)\}u(y_3,y_1)\\
&\label{F4c}+\int \mathrm{d}x \{\bar p(y_1,x)v_N(x-y_2)u(y_3,x)\}u(y_2,y_4)\\
\label{F4d}&+\int\mathrm{d}x_1\mathrm{d}x_2 \{\bar p(y_1,x_1)p(x_2,y_2)v_N(x_1-x_2)u(y_3,x_1)u(x_2,y_4)\}\Big\}.
\end{align}\end{subequations}

A standard enrgy estimate applied to \eqref{oe} using self-adjointness of $\mathcal{L}$ implies
\begin{align}
\nn \big\||\psi_\mathrm{ex}(t)\rangle-|\psi_\mathrm{ap}(t)\rangle\big\|_\mathbb{F}=\||\tilde{\psi}(t)\rangle\|_\mathbb{F}\leq N^{-1/2} \int_0^t\big\|\mathcal{E}(t_1)\vac\big\|_\mathbb{F}\,\mathrm{d}t_1.
\end{align}
For an estimate of the right hand side of the above inequality, we need $L^2$-norm estimates of the terms in \eqref{F1a}-\eqref{F4d}. This was done in \cite{GM} using the decay estimate
$\|\phi(t,\cdot)\|_{L^\infty(\mathbb{R}^3)}\lesssim 1/(1+t^{3/2})$  and the  estimate $\|u(t,\cdot)\|_{L^2(\mathbb{R}^6)}\lesssim \log(1+t)$. However  $0<\beta<1/3$ had to be assumed there for the final estimate in \eqref{GMest} to be meaningful.

While, in the current work, we extend the estimates on the error to the case of $\beta<1/2$ as stated in our main result (Theorem \ref{mainresult}), we can also provide here with the following heuristic argument suggesting that the uncoupled system consisting of \eqref{cf} and \eqref{efk} does not provide an approximation for $\beta\geq 1/2$. Indeed, we can write $|\tilde{\psi}\rangle$ e.g. as
\[|\tilde{\psi}\rangle=(0, \psi_{\text{\eqref{F1i}}},0,\dots)+\text{other contributions}\] 
where $\psi_{\text{\eqref{F1i}}}$ satisfies
\be\label{hear}\Big(\frac{1}{i}\partial_t-\Delta_{\mathbb{R}^3}\Big)\psi_{\text{\eqref{F1i}}}(t,y)=N^{-1/2}\int v_{N}(y-x)u(t,y,x)\bar\phi^{(N)}(t,x)\,\mathrm{d}x\quad\text{with}\quad\psi_{\text{\eqref{F1i}}}(0)=0\ee
in which the integral term on the right hand side comes from \eqref{F1i}. We added  the superscript $(N)$ to $\phi$ for recalling that it solves \eqref{cf} and hence it is $N$-dependent. We could have checked other similar contributions coming from \eqref{F1a}-\eqref{F4d} but \eqref{F1i} will serve the purpose. 
At this point using \eqref{efka} we can consider an approximate equation for $u$. Recalling $\mathrm{sh}(2k)=2u\circ c=2u+2u\circ p$, let's just look at
\[\frac{1}{i}\partial_tu-\Delta u+v_N(y_1-y_2)\phi^{(N)}(t,y_1)\phi^{(N)}(t,y_2)=0.\]
If we make the change of variables $x_1:=y_1-y_2$ and $x_2:=y_1+y_2$ then we have
\[\Big(\frac{1}{i}\partial_t-2\big(\Delta_{x_1}+\Delta_{x_2}\big)\Big)u(\mathsmaller{t,\frac{x_1+x_2}{2},\frac{x_2-x_1}{2}})=-v_N(x_1)\phi^{(N)}(\mathsmaller{t,\frac{x_1+x_2}{2}})\phi^{(N)}(\mathsmaller{t,\frac{x_2-x_1}{2}}).\]
Hence one can consider an ``approximate" solution 
\[u(t,y_1,y_2)=-N^\beta w(N^\beta(y_1-y_2))\phi^{(N)}(t,y_1)\phi^{(N)}(t,y_2)\quad\text{where}\quad \Delta w=-\frac{1}{2}v.\]
Inserting the above ansatz in \eqref{hear} gives
\[\Big(\frac{1}{i}\partial_t-\Delta_{\mathbb{R}^3}\Big)\psi_{\eqref{F1i}}=-N^{\beta-1/2} \underbrace{\Big\{\Big(N^{3\beta}v(N^\beta\cdot)w(N^\beta\cdot)\Big)\ast|\phi^{(N)}|^2\Big\}(t,y)\phi^{(N)}(t,y)}_{\mathclap{\substack{\text{converges to }(\int vw)|\phi(t,y)|^2\phi(t,y)\text{ as }N\rightarrow\infty\\\text{since }\phi^{(N)}\rightarrow\phi\text{ in }L^2\text{as in \eqref{ctld}}}}}\]
where $\phi^{(N)}$ and $\phi$ solve equations \eqref{cf} and \eqref{limitdyn} respectively. Hence, to ensure a decay for $\psi_{\eqref{F1i}}$ as $N\rightarrow\infty$, we have to consider $\beta<1/2$.

Finally in this section, let's describe our general strategy. 
From now on $\phi$ will always denote the solution to the equation \eqref{cf} as was the case in all definitions preceding the above heuristic argument. If we look at the terms in \eqref{F1}-\eqref{F4}, in all of them except \eqref{F2a}, \eqref{F3a} and \eqref{F4a}, the singularity associated with the interaction $v_N(x)=N^{3\beta}v(N^\beta x)$, which converges to $\big(\int v)\delta(x) $ as $N\rightarrow \infty$, 
is smoothed out due to the integration against functions with sufficient integrability properties. Hence we separate  $F_l(t,\cdot)$ defined in \eqref{F1}-\eqref{F4} into their regular and singular parts as follows, where super-scripts ``r" ans ``s" stand for ``regular" and ``singular" respectively: 
\begin{subequations}\label{Fsingular}\begin{align}
& \label{F2s}F_2^\mathrm{s}(t,y_1,y_2):=-(1/2N)v_N(y_1-y_2)\Big\{u(t,y_1,y_2)+(\bar p\circ u)(t,y_1,y_2)\Big\},\\
&\label{F3s}
F_3^\mathrm{s}(t,y_1,y_2,y_3):=-N^{-1/2}v_N(y_1-y_2)\phi(t,y_2)u(t,y_3,y_1),\\
&\label{F4s} F_4^\mathrm{s}(t,y_1,y_2,y_3,y_4):=-(1/2N)v_N(y_1-y_2)u(t,y_3,y_1)u(t,y_2,y_4)\text{ and}\\
& F_l^{\mathrm{r}}:=F_l- F_l^{\mathrm{s}}\text{ for }l=2,3,4.
\end{align}\end{subequations}
Using this, we split $|\tilde{\psi}\rangle$ first into its singular and regular parts as 
\begin{subequations}\begin{align}&\nn |\tilde{\psi}\rangle=|\tilde{\psi}^\mathrm{r}\rangle+|\tilde{\psi}^\mathrm{s}\rangle\text{ where }\\
&\label{9a}\Big(\frac{1}{i}\partial_t-\mathcal{L}\Big)|\tilde{\psi}^\mathrm{r}\rangle=(0,F_1,F_2^\mathrm{r},F_3^\mathrm{r},F_4^\mathrm{r},0,\dots),\\
&\label{9b}\Big(\frac{1}{i}\partial_t-\mathcal{L}\Big)|\tilde{\psi}^\mathrm{s}\rangle=(0,0,F_2^\mathrm{s},F_3^\mathrm{s},F_4^\mathrm{s},0,\dots),\\
&\nn |\tilde{\psi}^\mathrm{r}(0)\rangle=|\tilde{\psi}^\mathrm{s}(0)\rangle=0\end{align}\end{subequations}
which follows from \eqref{oe}. Energy estimate applied to \eqref{9a} implies
\be\label{eefr1}\||\tilde{\psi}^\mathrm{r}(t)\rangle\|_\mathbb{F}\lesssim\int_0^t\Big(\|F_1(t_1)\|_{L^2(\mathbb{R}^3)}+\sum_{l=2}^4\|F_l^\mathrm{r}(t_1)\|_{L^2(\mathbb{R}^{3l})}\Big)\,\mathrm{d}t_1.\ee
Hence we need to obtain $L^2$-norm estimates of $F_1$ and $F_l^{\mathrm{r}}, l=2,3,4$, which we do in section 4 after obtaining a priori estimates on the pair excitations in section 3. 

We will start dealing with the singular part of $|\tilde{\psi}\rangle$ in section 5 in which we will split $|\tilde{\psi}^{\mathrm{s}}\rangle$ in \eqref{9b} into its approximate and error parts as follows
\begin{subequations}\begin{align}&\nn |\tilde{\psi}^\mathrm{s}\rangle=|\tilde{\psi}_1^\mathrm{a}\rangle+|\tilde{\psi}_1^\mathrm{e}\rangle\text{ where }\\
&\label{12a}\Big(\frac{1}{i}\partial_t-\int L(t,x,y)a_x^\ast a_y\,\mathrm{d}x\mathrm{d}y\Big)|\tilde{\psi}_1^\mathrm{a}\rangle=(0,0,F_2^\mathrm{s},F_3^\mathrm{s},F_4^\mathrm{s},0,\dots),\\
&\label{12b}\Big(\frac{1}{i}\partial_t-\mathcal{L}\Big)|\tilde{\psi}_1^\mathrm{e}\rangle=-N^{-1/2}\mathcal{E}(t)|\tilde{\psi}_1^\mathrm{a}\rangle,\\&\nn |\tilde{\psi}_1^\mathrm{a}(0)\rangle=|\tilde{\psi}_1^\mathrm{e}(0)\rangle=0.\end{align}\end{subequations}
First we will obtain estimates on $|\tilde{\psi}_1^\mathrm{a}\rangle$ using an elliptic estimate and also Strichartz estimates along with Christ-Kiselev Lemma  after a suitable change of variables. Those will not provide us with sufficient integrability properties for the forcing term in \eqref{12b}. Hence we will also discuss the necessity to iterate the splitting procedure for some finitely many times before applying a final energy estimate to the error part of the solution at the final step of iteration.  We will prove the inductive step of the iteration and discuss its implications in section 6 leading to our main result in Theorem \ref{mainresult}.

\begin{center}\textbf{ 3. A priori estimates for the pair excitations}
\end{center} 

In this section we will prove estimates on mixed $L^p$ and Sobolev norms of the pair excitations which will be needed in estimating the terms in \eqref{F1}-\eqref{F4}. To keep the notation simple in what follows let's define $s_2$ and $p_2$  as $s_2:=\mathrm{sh}(2k)=2\mathrm{sh}(k)\circ\mathrm{ch}(k)$ and $p_2:=\mathrm{ch}(2k)-\delta(x-y)$ respectively. Then \eqref{efk} becomes
\begin{subequations}\begin{align} &\label{s2}\mathbf{S}(s_2)=2m+m\circ p_2+\bar p_2\circ m,\\ 
&\label{p2}\mathbf{W}(\bar p_2)=m\circ\bar s_2-s_2\circ\bar m,\\
&s_2(0,\cdot)=p_2(0,\cdot)=0.\nn\end{align}\end{subequations}
Let's also recall our notation $u:=\mathrm{sh}(k)$, $c:=\mathrm{ch}(k)=\delta(x-y)+p$ from the previous section.
Our main result in this section is the following:
\bt \label{mrs2} Let the initial data $\phi_0$ in \eqref{cf} be in $W^{m,1}(\mathbb{R}^3)$ ($m$ derivatives in $L^1$) for $m\geq 6$ and let $(\partial_ts_2)(0,\cdot)$ be sufficiently regular (to be specified later in the proof). Then the following estimates hold:
\begin{align}&\label{s2hth}\|\partial_t^js_2(t,\cdot)\|_{H^{3/2}}\lesssim_\epsilon N^{\beta(1+\epsilon)}\log(1+t)\quad\text{for }j=0,1\\
&\label{uhth}\|u(t,\cdot)\|_{H^{3/2}}\lesssim_\epsilon N^{\beta(1+\epsilon)}\log(1+t)\\
&\label{umlp}\|u(t,x,y)\|_{L^\infty(\mathrm{d}y;L^2(\mathrm{d}x))}:=\left\|\|u(t,x,y)\|_{L^2(\mathrm{d}x)}\right\|_{L^\infty(\mathrm{d}y)}\lesssim_\epsilon N^{\beta(1+\epsilon)}\log(1+t)\end{align}
for any $\epsilon>0$ and $0<\beta\leq 1$.\et

We will need the following lemmas for the proof Theorem \ref{mrs2}:
\bl \label{lemma1}\text{\textnormal{(Proposition 3.3,  Corollary 3.4, Corollary 3.5 in \cite{GM})}}
Let $\phi$ be a solution of \eqref{cf}.  
\begin{itemize}
\item[(i)] There exists $C_s$ depending only on $\|\phi_0\|_{H^s(\mathbb{R}^3)}$ such that
\be\label{Hsuni}\|\phi(t,\cdot)\|_{H^s(\mathbb{R}^3)}\leq C_s\text{ uniformly  in time.}\ee
\item[(ii)] Assuming $\phi_0\in W^{m,1}$ for $m\geq 2$,
\begin{align}
&\label{finfty}\|\partial_t^j\phi(t,\cdot)\|_{L^\infty(\mathbb{R}^3)}\lesssim (1+t^{3/2})^{-1}\text{ and }\\
&\label{ft}\|\partial_t^j\phi(t,\cdot)\|_{L^3(\mathbb{R}^3)}\lesssim (1+t^{1/2})^{-1}
\text{ for }j=0,1.\end{align}\end{itemize}\el

\noindent\textbf{Remark.} \textit{Note that in case of $j=0$, \eqref{ft} follows by interpolating \eqref{finfty} with mass conservation and in case of $j=1$, by interpolating \eqref{finfty} with 
\be\label{dtfi}\|\partial_t\phi(t,\cdot)\|_{L^2(\mathbb{R}^3)}\lesssim \|\phi(t,\cdot)\|_{H^2(\mathbb{R}^3)}+\underbrace{\|(v_N\ast|\phi(t,\cdot)|^2)\phi(t,\cdot)\|_{L^2(\mathbb{R}^3)}}_{\leq\|v\|_1\|\phi\|_4^2\|\phi\|_\infty}\leq \text{const.}\ee
We will also frequently use 
\be\label{fay4}\|\partial_t^j\phi(t,\cdot)\|_{L^4(\mathbb{R}^3)}\lesssim (1+t^{3/4})^{-1}\quad\text{for }j=0,1\ee
which follows again by interpolation.}
\bc\label{htd}\eqref{finfty}-\eqref{ft} hold for $j\geq2$ if $\phi_0\in W^{m,1}$ for $m$ sufficiently large.\ec

\noindent\textit{Proof Sketch}. Since  we will  only need the estimates on the second and third order time derivatives, let's  provide here with an outline of the proof in case of the second order time derivative, which can be modified to obtain estimates on higher time derivatives. We claim that if $\phi$ solves \eqref{cf} with $\phi_0\in W^{m,1}$ for $m\geq 4$ ($m\geq6$ in case of third order time derivative) then we have 
\be \label{dt2fest}\|\partial_t^2\phi(t)\|_{L^\infty(\mathrm{R}^3)}\lesssim \frac{1}{1+t^{3/2}}.\ee
To prove this estimate let's differentiate \eqref{cf} with respect to time twice and solve the resulting equation for $\partial_t^2\phi$ by Duhamel's formula. Then we have 
\begin{align}\big \|\partial_t^2\phi(t)\big\|_\infty &\leq \big\| e^{it\Delta}\big(\partial_t^2\phi\big)_0\big\|_\infty +\int_0^t \Big\|e^{i(t-s)\Delta}\partial_s^2\big[(v_N\ast|\phi|^2)\phi(s)\big]\Big\|_\infty\,\mathrm{d}s.
\end{align}  
Assuming $t>1$, we split the above integral and, to estimate the integrand,  we use the standard $L^\infty L^1$ decay estimate for the linear equation when we integrate over $(0,t-1)$. For the part of the same integral on $(t-1,t)$,  we first use the Sobolev embedding $W^{3+\epsilon,1}(\mathbb{R}^3)\hookrightarrow L^\infty(\mathbb{R}^3)$ and then the $L^{3}L^{3/2}$ decay estimate for the linear equation, up to modifying the exponents by a small amount. Hence we obtain
\begin{align}
\nn\big \|\partial_t^2\phi(t)\big\|_\infty\lesssim& \frac{\|\big(\partial_t^2\phi\big)_0\|_1}{1+t^{3/2}}+\int_0^{t-1}\frac{1}{1+|t-s|^{3/2}}\big\|\partial_s^2\big[(v_N\ast|\phi|^2)\phi(s)\big]\big\|_1\mathrm{d}s\\
&\label{dt2fi}+\int_{t-1}^t\frac{1}{1+|t-s|^{1/2+\epsilon}}\big\|\nabla\partial_s^2\big[(v_N\ast|\phi|^2)\phi(s)\big]\big\|_{3/2-\epsilon'}\mathrm{d}s.
\end{align}
Now we have $\|\partial^\alpha\phi\|_2\lesssim \|\phi\|_{H^{|\alpha|}}\leq C_{|\alpha|}$ and also $\|\partial^{\alpha}\phi\|_\infty\lesssim\|\partial^\alpha
\phi\|_{H^2}\leq C_{2+|\alpha|}$ which follow from \eqref{Hsuni}. Interpolating, we obtain 
\begin{align}\label{reg}\left.\begin{array}{ll}\bullet \,\|\partial^\alpha \phi\|_p\leq C_{\alpha, p} \text{ for } p\geq 2 \text{ and for spatial derivatives } \partial^\alpha \text{ of all orders}. 
\\ \bullet \text{ We can extend this to the case of derivatives including the time variable  if}\\ \,\,\text{ we take derivatives in \eqref{cf} as needed and use the estimates obtained so far.}\end{array}\right\}\end{align}
These regularity and integrability properties together with \eqref{finfty} imply
\begin{align}
\nn&\big\|\partial_s^2\big[(v_N\ast|\phi|^2)\phi(s)\big]\big\|_1\lesssim \|\phi(s)\|_\infty\lesssim\frac{1}{1+s^{3/2}},\\
&\nn\big\|\nabla\partial_s^2\big[(v_N\ast|\phi|^2)\phi(s)\big]\big\|_{3/2-\epsilon'}\lesssim \|\phi(s)\|_\infty+\|\partial_s\phi(s)\|_\infty \lesssim \frac{1}{1+s^{3/2}}.
\end{align}
Inserting these in \eqref{dt2fi} implies our claim in \eqref{dt2fest}. We also have  
\be\label{dt2f3}\|\partial_t^2\phi(t)\|_3\lesssim \frac{1}{1+t^{1/2}}\ee
by interpolation between \eqref{dt2fest} and $L^2$-norm which is uniformly bounded. \hfill$\Box$

\vspace{0.4cm}

Before stating the next lemma, let's write the kinetic and the potential parts of $g$ (see \eqref{112a}) separately as 
\be\label{gpot} g=-\Delta_x\delta(x-y)+g_\mathrm{pot}.\ee
then we can define $V$  as follows \be V(u):=g^\mathrm{T}_\mathrm{pot}\circ u+u\circ g_\mathrm{pot}.\ee
Explicitly,
\begin{align}\label{Vex} V(u)(t,x,y)=&\Big((v_N\ast|\phi|^2)(t,x)+(v_N\ast|\phi|^2)(t,y)\Big)u(x,y)\\
&\nn+\int v_N(x-z)\phi(t,x)\bar\phi(t,z)u(z,y)\mathrm{d}z+\int u(x,z)v_N(z-y)\bar\phi(t,z)\phi(t,y)\mathrm{d}z.\end{align}
This allows us to write the potential part of  $\mathbf{S}$ (see \eqref{bs}) separately:
\be\mathbf{S}(\cdot)=\Big(\frac{1}{i}\partial_t-\Delta\Big)(\cdot)+V(\cdot).\ee
We will split $s_2$ satisfying \eqref{s2} as
\be \label{s2fsplit}s_2=s_a+s_e\ee 
where $s_a$ satisfies the equation $\bi{\mathrm{S}}(s_a)=2m=-2v_N(x-y)\phi(t,x)\phi(t,y)$ and it represents the singular part of $s_2$ since
\be\label{m32est}\|m(t,\cdot)\|_{L^2(\mathbb{R}^6)}= \Big(v_N^2\ast|\phi(t,\cdot)|^2\mathlarger{\bi{,}} \,|\phi(t,\cdot)|^2\Big)^{\frac{1}{2}}\lesssim \|v_N\|_{L^2(\mathbb{R}^3)}\|\phi(t,\cdot)\|_{L^4(\mathbb{R}^3)}^2\mathop{\lesssim}^{\substack{\text{by}\\ \eqref{fay4}}} N^{3\beta/2}(1+t^{3/2})^{-1}\ee blows up as $N\rightarrow\infty$. We further split $s_a$ into its approximate and error parts as 
\be\label{sasplit}s_a=s_a^0+s_a^1\ee
and we have the following set of equations  being equivalent to \eqref{s2}:\begin{subequations}
\begin{align}
&\label{mostsing}\Big(\frac{1}{i}\partial_t-\Delta\Big)s^0_a=2m\\
&\label{sa1}\mathbf{S}(s^1_a)=-V(s^0_a)\\
& \label{se}\mathbf{S}(s_e)=m\circ p_2+\bar p_2\circ m\\
& \nn s_a^0(0)=s_a^1(0)=s_e(0)=0.
\end{align}
\end{subequations} 
We are ready to state the next lemma:
\bl\label{lemma2}
Assuming $\phi_0\in W^{m,1}$ for $m\geq 2$ in \eqref{cf}, the following estimates hold: 
\be\label{sa0sa1}\|s_a^0(t,\cdot)\|_{L^2(\mathbb{R}^6)}\lesssim \log(1+t),\qquad \|s_a^1(t,\cdot)\|_{L^2(\mathbb{R}^6)}\lesssim 1\ee
which imply
\be \label{sep2}\|s_e(t,\cdot)\|_{L^2(\mathbb{R}^6)}\lesssim 1,  \qquad\|p_2(t,\cdot)\|_{L^2(\mathbb{R}^6)}\lesssim 1.\ee
Since $s_2=s_a^0+s_a^1+s_e$, we also have
\begin{equation}\|s_2(t,\cdot)\|_{L^2(\mathbb{R}^6)}\lesssim\log(1+t).\end{equation}
Finally since $s_2=\mathrm{sh}(2k)=2\mathrm{sh}(k)\circ \mathrm{ch}(k)$ and $\|\mathrm{ch}(k)^{-1}\|_\mathrm{operator}$ is uniformly bounded, recalling the notation $u=\mathrm{sh}(k)$ and $p=\mathrm{ch}(k)-\delta(x-y)$, we have
\be\label{pul2}\|p(t,\cdot)\|_{L^2(\mathbb{R}^6)}\leq\|u(t,\cdot)\|_{L^2(\mathbb{R}^6)}\lesssim \log(1+t)\ee
where the first inequality follows from taking traces in the relation $p\circ p+2p=\bar u\circ u$ and using $p(x,x)\geq 0$.   Constants involved in the above estimates depend only on $\|\phi_0\|_{W^{m,1}}$.\el
\textbf{Remark.} \textit{For the proof of the first inequality in \eqref{sa0sa1}, one solves equation \eqref{mostsing} by Duhamel's formula and, after an integration by parts, uses the elliptic estimates below (Lemma 4.3 in \cite{GM}) along with \eqref{ft}:
\begin{subequations}\label{elliptic}\begin{align}&\int \frac{|\hat{m}(t,\xi,\eta)|^2}{(|\xi |^2+|\eta |^2)^2}\mathrm{d}\xi\mathrm{d}\eta\lesssim \|\phi(t,\cdot)\|_3^4,\\&\int \frac{|\partial_t\hat m(t,\xi,\eta)|^2}{(|\xi |^2+|\eta |^2)^2}\mathrm{d}\xi\mathrm{d}\eta\lesssim \|\phi(t,\cdot)\|_3^2\|\partial_t\phi(t,\cdot)\|_3^2\text{ and}\\&\nn\text{similar estimates hold for higher time derivatives.}\end{align}\end{subequations}
 The proof of the second inequality in \eqref{sa0sa1} is achieved by applying an energy estimate to the equation \eqref{sa1} and using the first inequality in \eqref{sa0sa1}. A final application of energy estimates to the equations \eqref{se} and \eqref{p2} together with the estimates in \eqref{sa0sa1} implies the estimates in \eqref{sep2}. We refer for more details to the proofs of Lemma 4.4 and Lemma 4.5 in \cite{GM}.} 

\vspace{0.4cm}
Now with the help of Lemma \ref{lemma1}, Corollary \ref{htd} and Lemma \ref{lemma2}, we can prove Theorem \ref{mrs2}.\newline

\noindent\textit{Proof of Theorem \ref{mrs2}.} \textbf{Proof of \eqref{s2hth}.} Recalling that $s_2=s_a^0+s_a^1+s_e$ from \eqref{s2fsplit} and \eqref{sasplit}, we will prove \eqref{s2hth} in two steps.
 
\textbf{Step 1} \underline{Estimates on  $\|s_a^0\|_{H^{3/2}}$ and $\|\partial_ts_a^0\|_{H^{3/2}}$}: We will first estimate $H^2$ and $H^{1/2-\epsilon}$-norms and then interpolate. 

Differentiating  \eqref{mostsing} as needed, solving the corresponding equations by Duhamel's formula and using integration by parts give
\begin{align}\label{dhmldt}&\nn\big(\partial_t^j\widehat{s_a^0}\big)(t,\xi,\eta)=\overbrace{\big(\partial_t^j\widehat{s_a^0}\big)(0,\xi,\eta}^{\mathclap{\text{equals }0\text{ if }j=0}})\\
&+\frac{2e^{-it(|\xi|^2+|\eta|^2)}}{|\xi|^2+|\eta|^2}\Big(\partial_t^j\hat m(t,\xi,\eta)e^{it(|\xi|^2+|\eta|^2)}-\big(\partial_t^j\hat m\big)(0,\xi,\eta)-\int_0^te^{is(|\xi|^2+|\eta|^2)}\partial_s^{j+1}\hat m(s,\xi,\eta)\mathrm{d}s\Big)\end{align} 
which implies
\begin{align}
\label{deltadtjsa0}\|\Delta \partial_t^js_a^0(t,\cdot)\|_2\lesssim \underbrace{\|\Delta\big(\partial_t^js_a^0\big)(0,\cdot)\|_2}_{\mathclap{\substack{\text{equals }0\text{ if }j=0\text{ and}\\\text{assumed to be finite for }j=1}}}+\|\big(\partial_t^jm\big)(0,\cdot)\|_2+\|\partial_t^jm(t,\cdot)\|_2+\int_0^t\|\partial_s^{j+1} m(s,\cdot)\|_2\mathrm{d}s\end{align}
Applying estimate \eqref{m32est} and the following estimates
\begin{align}
&\nonumber\|\partial_s m(s,\cdot)\|_2\lesssim \Big(v_N^2\ast|\partial_s\phi(s,\cdot)|^2\mathlarger{\bi{,}} \,|\phi(s,\cdot)|^2\Big)^{1/2}\leq \|\phi(s,\cdot)\|_\infty\|v_N\|_2\|\partial_s\phi(s,\cdot)\|_2\mathop{\,\lesssim\,}_{\mathclap{\substack{\text{by}\\\eqref{finfty}
\text{ and }\eqref{dtfi}}}} N^{3\beta/2}(1+s^{3/2})^{-1}\\
&\nn\|\partial_s^2 m(s,\cdot)\|_2\lesssim \Big(v_N^2\ast|\partial_s^2\phi(s,\cdot)|^2\mathlarger{\bi{,}} \,|\phi(s,\cdot)|^2\Big)^{1/2}+\Big(v_N^2\ast|\partial_s\phi(s,\cdot)|^2\mathlarger{\bi{,}} \,|\partial_s\phi(s,\cdot)|^2\Big)^{1/2}\\&\nn\hspace{1.9cm}\leq \|\phi(s,\cdot)\|_\infty\|v_N\|_2\|\partial_s^2\phi(s,\cdot)\|_2+\|\partial_s\phi(s,\cdot)\|_\infty\|v_N\|_2\|\partial_s\phi(s,\cdot)\|_2\mathop{\lesssim}_{\mathclap{\substack{\text{by }\\\eqref{finfty}\text{ and } \eqref{reg}}}}N^{3\beta/2}(1+s^{3/2})^{-1}\end{align}
to \eqref{deltadtjsa0}, we obtain
\be\boxed{\label{sa0h3h}\|\partial_t^js_a^0(t,\cdot)\|_{H^2}\lesssim N^{3\beta/2}\quad\text{for }j=0,1.}\ee

We will next estimate $\|\partial_t^js_a^0(t,\cdot)\|_{H^{1/2-\epsilon'}}$, $j=0,1$ for $\epsilon'>0$ small and to be determined later. Again by using  \eqref{dhmldt}
\begin{align}
&\nn \|D^{1/2-\epsilon'}\partial_t^js_a^0(t,\cdot)\|_2\simeq\|(|\xi |+|\eta |)^{1/2-\epsilon'}\partial_t^j\widehat{s_a^0}(t,\xi,\eta)\|_2\\
&\lesssim\nn\|\overbrace{D^{1/2-\epsilon'}\big(\partial_t^js_a^0\big)(0,\cdot)}^{\text{equals }0\text{ if }j=0}\|_2+\Big\|\frac{\big(\partial_t^j\hat m\big)(0,\xi,\eta)}{(|\xi|+|\eta|)^{3/2+\epsilon'}}\Big\|_2+\Big\|\frac{\partial_t^j\hat m(t,\xi,\eta)}{(|\xi|+|\eta|)^{3/2+\epsilon'}}\Big\|_2\\
&\hspace{0.4cm}\label{D1h2}+ \int_0^t\Big\|\frac{\partial_s^{j+1}\hat m(s,\xi,\eta)}{(|\xi|+|\eta|)^{3/2+\epsilon'}}\Big\|_2\mathrm{d}s.\end{align}
Now we need estimates of 
$$\Big\|\frac{\partial_t^j\hat m(t,\xi,\eta)}{(|\xi|+|\eta|)^{3/2+\epsilon'}}\Big\|_2 \quad\text{for }j=0,1,2.$$
We will prove the estimates on the above terms similarly to the proof of \eqref{elliptic}. Let's do it first for the case $j=0$. Writing
$$-m(t,x,y)=v_N(x-y)\phi(t,x)\phi(t,y)=\int\delta(x-y-z)v_N(z)\phi(t,x)\phi(t,y)\mathrm{d}z$$
and considering the Fourier transform of $\delta(x-y-z)\phi(t,x)\phi(t,y)$ in the variables $x,\,y$: 
$$e^{iz\cdot\eta}\widehat{\phi\phi_z}(t,\xi+\eta)\text{ where }\phi_z(x)=\phi(x-z)$$
we can write
$$|\hat m(t,\xi,\eta)|^2=\Big|\int v_N(z)e^{iz\cdot\eta}\widehat{\phi\phi_z}(t,\xi+\eta)\mathrm{d}z\Big|^2\leq \|v\|_1\int|v_N(z)||\widehat{\phi\phi_z}(t,\xi+\eta)|^2\mathrm{d}z.$$
Hence after a change of variables
\begin{align}\nonumber\Big\|\frac{\hat m(t,\xi,\eta)}{(|\xi|+|\eta|)^{3/2+\epsilon'}}\Big\|_2^2\lesssim \int |v_N(z)|\frac{|\widehat{\phi\phi_z}(t,\xi)|^2}{(|\xi|+|\eta|)^{3+2\epsilon'}}\,\mathrm{d}\xi\,\mathrm{d}\eta\,\mathrm{d}z\lesssim\frac{1}{\epsilon'}\int|v_N(z)|\frac{|\widehat{\phi\phi_z}(t,\xi)|^2}{|\xi|^{2\epsilon'}}\mathrm{d}\xi\mathrm{d}z.\end{align}
Combining this last estimate with
\begin{align}\nonumber\int\frac{|\widehat{\phi\phi_z}(t,\xi)|^2}{|\xi|^{2\epsilon'}}\mathrm{d}\xi\lesssim \|D^{-\epsilon'}(\phi\phi_z)\|_2^2\mathop{\,\lesssim_{\epsilon'}\,}_{\mathclap{\substack{\text{by}\\\text{Hardy-}\\\text{Littlewood-}\\\text{Sobolev}}}}\|\phi\phi_z\|_{2-\epsilon''}^2\leq \|\phi\|_{4-2\epsilon''}^4\text{ where }\epsilon''=\frac{4\epsilon'}{3+2\epsilon'}\end{align} 
gives
\be\label{fmtep1}\Big\|\frac{\hat m(t,\xi,\eta)}{(|\xi|+|\eta|)^{3/2+\epsilon'}}\Big\|_2 \lesssim_{\epsilon'}\|\phi\|_{4-2\epsilon''}^2. \ee
We can prove similarly in general  the following estimate:
\begin{align}
&\label{fmtep3}\Big\|\frac{\partial_t^j\hat m(t,\xi,\eta)}{(|\xi|+|\eta|)^{3/2+\epsilon'}}\Big\|_2\lesssim_{\epsilon'} \sum_{l=0}^j\|\partial_t^l\phi\|_{4-2\epsilon''}\|\partial_t^{j-l}\phi\|_{4-2\epsilon''}.\end{align}
Inserting estimates \eqref{fmtep1}-\eqref{fmtep3} into  \eqref{D1h2} gives
\begin{align}&\nn\|s_a^0(t,\cdot)\|_{H^{1/2-\epsilon'}}\lesssim_{\epsilon'}\|\phi(0,\cdot)\|_{4-2\epsilon''}^2+\|\phi(t,\cdot)\|_{4-2\epsilon''}^2+\int_0^t\|\phi(s,\cdot)\|_{4-2\epsilon''}\|\partial_s\phi(s,\cdot)\|_{4-2\epsilon''}\mathrm{d}s\\
&\nn\|\partial_ts_a^0(t,\cdot)\|_{H^{1/2-\epsilon'}}\lesssim_{\epsilon'}\|\big(\partial_ts_a^0\big)(0,\cdot)\|_{H^{1/2-\epsilon'}}+\|\phi(0,\cdot)\|_{4-2\epsilon''}\|\big(\partial_t\phi\big)(0,\cdot)\|_{4-2\epsilon''}\\[2mm]&\nn\hspace{3.3cm}+\|\phi(t,\cdot)\|_{4-2\epsilon''}\|\partial_t\phi(t,\cdot)\|_{4-2\epsilon''}\\
&\nn\hspace{3.3cm}+\int_0^t\Big(\|\phi(s,\cdot)\|_{4-2\epsilon''}\|\partial_s^2\phi(s,\cdot)\|_{4-2\epsilon''}+\|\partial_s\phi(s,\cdot)\|_{4-2\epsilon''}^2\Big)\mathrm{d}s\end{align}
Using $\|\partial_t^j\phi(t,\cdot)\|_{4-2\epsilon''}\lesssim (1+t^{3/4-\epsilon'/2})^{-1}$ for $j=0,1,2$, which follow by interpolating $L^2$-norm with $L^\infty$ estimates (see Corollary \ref{htd} and \eqref{reg}), we obtain
\be\boxed{\|\partial_t^js_a^0(t,\cdot)\|_{H^{1/2-\epsilon'}}\lesssim_{\epsilon'} 1\quad\text{for }j=0,1.}\ee
Interpolating this with \eqref{sa0h3h} 
gives
$$\|\partial_t^js_a^0\|_{H^{3/2}}\leq \|\partial_t^js_a^0\|_{H^2}^{\frac{2+2\epsilon'}{3+2\epsilon'}}\|\partial_t^js_a^0\|_{H^{1/2-\epsilon'}}^\frac{1}{3+2\epsilon'}\lesssim_{\epsilon'} (N^\frac{3\beta}{2})^{\frac{2+2\epsilon'}{3+2\epsilon'}}\quad\text{for }j=0,1.$$
Hence finally we obtain
\be\label{sa0h32}\boxed{\|\partial_t^js_a^0(t,\cdot)\|_{H^{3/2}}\lesssim_\epsilon N^{\beta(1+\epsilon)}\quad\text{for }j=0,1}\quad \text{ where }\epsilon=\frac{\epsilon'}{3+2\epsilon'}\ee
So for $\epsilon>0$ arbitrarily small, we can choose $\epsilon'=3\epsilon/(1-2\epsilon)$ in the above estimates leading to \eqref{sa0h32}. \newline

\textbf{Step 2} \underline{Estimates on $\|\partial_t^js_a^1\|_{H^{3/2}}$ and $\|\partial_t^js_e\|_{H^{3/2}}$ for $j=0,1$}: We will first estimate $H^2$-norms then we will use the Sobolev embedding $H^2\hookrightarrow H^{3/2}$. We will obtain $H^2$-estimates of $\partial_t^js_a^1$ and $\partial_t^js_e$ by estimating $\partial_t^{j+1}s_a^1$ and $\partial_t^{j+1}s_e$ in $L^2$ first and then using the equations satisfied by $\partial_t^js_a^1$ and $\partial_t^js_e$.
If we take derivative on both sides in \eqref{sa1} and recall that $s_a=s_a^0+s_a^1$ from  \eqref{sasplit}, we can write
\begin{align}&\label{Swdts}\mathbf{S}(\partial_t s_a^1)=-V(\partial_t s_a^0)-\Big(\big(\partial_t g_\mathrm{pot}^\mathrm{T}\big)\circ s_a+s_a\circ \big(\partial_t g_\mathrm{pot}\big)\Big),\\
&\nn\mathbf{S}(\partial_t^2 s_a^1)=-V(\partial_t^2 s_a^0)-2\Big(\big(\partial_t g_\mathrm{pot}^\mathrm{T}\big)\circ \partial_ts_a+(\partial_ts_a)\circ \big(\partial_t g_\mathrm{pot}\big)\Big)\\
&\label{Swd2ts}\hspace{1.57cm}-\Big(\big(\partial_t^2 g_\mathrm{pot}^\mathrm{T}\big)\circ s_a+s_a\circ \big(\partial_t^2 g_\mathrm{pot}\big)\Big)\end{align}
where 
\begin{align}\nn\partial_tg_\mathrm{pot}(t,x,y)=&\Big(v_N\ast\big(2\mathrm{Re}(\bar\phi\partial_t\phi)\big)\Big)(t,x)\delta(x-y)\\
&\nn+v_N(x-y)\partial_t\bar\phi(t,x)\phi(t,y)+v_N(x-y)\bar\phi(t,x)\partial_t\phi(t,y)\end{align}
and we can compute $\partial_t^2 g_\mathrm{pot}$ likewise.
We will apply an energy estimates to \eqref{Swdts}-\eqref{Swd2ts}. Let's define
\be\label{vtdef}(\partial_t^jV)(u):=\big(\partial_t^j g_\mathrm{pot}^\mathrm{T}\big)\circ u+u\circ \big(\partial_t^j g_\mathrm{pot}\big)\quad\text{for }j=1,2.\ee
Using this definition and \eqref{Swdts}-\eqref{Swd2ts} and also recalling \eqref{GMeqns}, we have \vspace{0.1cm}
\begin{align}&\nn\mathbf{W}\big((\partial_t s_a^1)\circ \partial_t \overline{s_a^1}\big)=\mathbf{S}(\partial_t s_a^1)\circ\partial_t \overline{s_a^1}-(\partial_t s_a^1)\circ\overline{\mathbf{S}(\partial_t s_a^1)}\\
&\nn = -V(\partial_t s_a^0)\circ\partial_t\overline{s_a^1}+(\partial_t s_a^1)\circ\overline{V(\partial_t s_a^0)}-\big((\partial_t V)(s_a)\big)\circ\partial_t \overline{s_a^1}
+(\partial_t s_a^1)\circ\overline{(\partial_t V)(s_a)},\\[0.3cm]
&\nn\mathbf{W}\big((\partial_t^2 s_a^1)\circ \partial_t^2 \overline{s_a^1}\big)=\mathbf{S}(\partial_t^2 s_a^1)\circ\partial_t^2 \overline{s_a^1}-(\partial_t^2 s_a^1)\circ\overline{\mathbf{S}(\partial_t^2 s_a^1)}\\
&\nn=-V(\partial_t^2 s_a^0)\circ\partial_t^2\overline{s_a^1}+(\partial_t^2 s_a^1)\circ\overline{V(\partial_t^2 s_a^0)}\\&\nn\hspace{0.4cm}-2\Big[\big((\partial_tV)(\partial_ts_a)\big)\circ\partial_t^2\overline{s_a^1}-(\partial_t^2 s_a^1)\circ\overline{(\partial_tV)(\partial_ts_a)}\Big]\\
&\nn\hspace{0.4cm}-(\partial_t^2V)(s_a)\circ\partial_t^2 \overline{s_a^1}+(\partial_t^2 s_a^1)\circ\overline{(\partial_t^2V)(s_a)}\end{align}
To obtain $L^2$-norm estimates, we take traces on both sides of the above equations and make the following estimates:
\begin{align}&\label{dtsa1}\partial_t\norms{\partial_ts_a^1}{2}{2}\leq 2\Big(\norms{V(\partial_t s_a^0)}{}{2}+\norms{(\partial_t V)(s_a)}{}{2}\Big)\norms{\partial_t s_a^1}{}{2}\\
&\label{dt2sa1}\partial_t\norms{\partial_t^2s_a^1}{2}{2}\leq 2\Big(\|V(\partial_t^2 s_a^0)\|_2+2\|(\partial_tV)(\partial_ts_a)\|_2+\|(\partial_t^2V)(s_a)\|_2\Big)\|\partial_t^2s_a^1\|_2\end{align}
Both $V$ $\big($see \eqref{Vex}$\big)$ and $\partial_t^j V$  for $j=1,2$ $\big($see \eqref{vtdef}$\big)$ are bounded from $L^2$ to $L^2$ since the inequalities\vspace{0.1cm} 
\begin{equation}\label{efVad}
\left.\begin{split}
 & \|\big(v_N\ast\partial_t^j|\phi|^2\big)(t,x)u(x,y)\|_{L^2_{x,y}}\\
&\lesssim \|v_N\|_{L^1(\mathbb{R}^3)}\Big(\sum_{k=0}^j\|\partial_t^k\phi(t,\cdot)\|_{L^\infty(\mathbb{R}^3)}\|\partial_t^{j-k}\phi(t,\cdot)\|_{L^\infty(\mathbb{R}^3)}\Big)\|u\|_{L^2(\mathbb{R}^6)},\\
&\|\int v_N(x-z)\partial_t^j\big(\bar\phi(t,x)\phi(t,z)\big)u(z,y)\mathrm{d}z\|_{L^2_{x,y}}\\
&\lesssim\Big( \sum_{k=0}^j\|\partial_t^k\phi(t,\cdot)\|_{L^\infty(\mathbb{R}^3)}\|\partial_t^{j-k}\phi(t,\cdot)\|_{L^\infty(\mathbb{R}^3)}\Big)\smash{\overbrace{\|\big(v_N\ast\|u(\cdot,y)\|_{L^2_y}\big)(x)\big\|_{L^2_x}}^{\lesssim\|v_N\|_{L^1(\mathbb{R}^3)}\|u\|_{L^2(\mathbb{R}^6)}}}
\end{split}\right\}
\end{equation}
for $j=0,1,2$ and $\|\partial_t^j\phi(t,\cdot)\|_{L^\infty(\mathbb{R}^3)}\lesssim 1/(1+t^{3/2})$ (see Corollary \ref{htd}) imply
\be\label{vops}\left.\begin{split}&\norms{V}{}{\mathrm{op}}\lesssim \norms{\phi(t,\cdot)}{2}{\infty}\lesssim (1+t^3)^{-1},\\&\norms{\partial_tV}{}{\mathrm{op}}\lesssim \norms{\phi(t,\cdot)}{}{\infty}\norms{\partial_t\phi(t,\cdot)}{}{\infty}\lesssim (1+t^3)^{-1},\\
&\norms{\partial_t^2V}{}{\mathrm{op}}\lesssim \norms{\phi(t,\cdot)}{}{\infty}\norms{\partial_t^2\phi(t,\cdot)}{}{\infty}+\norms{\partial_t\phi(t,\cdot)}{2}{\infty}\lesssim (1+t^3)^{-1}.\end{split}\right\}\ee
Hence \eqref{dtsa1}-\eqref{dt2sa1} take the form
\begin{align}&\label{sn}\partial_t\|\partial_ts_a^1(t,\cdot)\|_2\lesssim\frac{1}{1+t^3}\Big(\|\partial_ts_a^0(t,\cdot)\|_2+\underbrace{\|s_a(t,\cdot)\|_2}_{\mathclap{\substack{\quad\qquad\lesssim \log(1+t)\text{ by \eqref{sa0sa1}}\\\quad\qquad\text{since }s_a=s_a^0+s_a^1}}}\Big),\\
&\label{dtsn}\partial_t\|\partial_t^2s_a^1(t,\cdot)\|_2\lesssim\frac{1}{1+t^3}\Big(\|\partial_t^2s_a^0(t,\cdot)\|_2+\|\partial_ts_a(t,\cdot)\|_2+\underbrace{\|s_a(t,\cdot)\|_2}_{\mathclap{\substack{\lesssim \log(1+t)}}}\Big).\end{align}
Now we need estimates of $\|\partial_t^j s_a^0\|_2$, $j=1,2$.
Taking $L^2$-norms in \eqref{dhmldt} and using \eqref{elliptic}, we can obtain the following estimate:
\[\begin{split}\label{dersa0}&\nn\norms{\partial_t^js^0_a(t,\cdot)}{}{2}\\
&\leq\|(\partial_t^js_a^0)(0,\cdot)\|_2+ \norm{\frac{(\partial_t^j\hat m)(0,\xi,\eta)}{|\xi|^2+|\eta|^2}}{}{2}+\norm{\frac{\partial_t^j\hat m(t,\xi,\eta)}{|\xi|^2+|\eta|^2}}{}{2}+\norm{\int_0^t e^{is(|\xi|^2+|\eta|^2)}\frac{\partial_s^{j+1} \hat m(s,\xi,\eta)}{|\xi|^2+|\eta|^2} \mathrm{d}s}{}{2}\\
&\nonumber\lesssim\|(\partial_t^js_a^0)(0,\cdot)\|_2+\sum_{l=0}^j\Big(\norms{(\partial_t^l\phi)(0,\cdot)}{}{3}\norms{(\partial_t^{j-l}\phi)(0,\cdot)}{}{3}+\norms{\partial_t^l\phi(t,\cdot)}{}{3}\norms{\partial_t^{j-l}\phi(t,\cdot)}{}{3}\Big)\\&\quad+\int_0^t\big\{\sum_{l=0}^{j+1}\norms{\partial_t^l\phi(t,\cdot)}{}{3}\norms{\partial_t^{j+1-l}\phi(t,\cdot)}{}{3}\big\}\mathrm{d}s.\end{split}\]
This last estimate considered with  \eqref{ft} and Corollary \ref{htd} imply
\be\label{dtsa0l2}\|\partial_t^js_a^0(t,\cdot)\|_2\lesssim\log(1+t)\quad\text{for }j=1,2.\ee
Inserting this in \eqref{sn} gives
\be\label{dtsa1di}\begin{split}\partial_t\norms{\partial_ts_a^1}{}{2}\lesssim \frac{\log(1+t)}{1+t^3} \end{split}\ee
which implies uniform-in-time boundedness of $\norms{\partial_ts_a^1}{}{2}$. This together with \eqref{dtsa0l2} implies
\be\label{l2dtsa}\|\partial_ts_a\|_2\lesssim \log(1+t)\ee
since $s_a=s_a^0+s_a^1$. Inserting this last estimate  and estimate \eqref{dtsa0l2} in \eqref{dtsn} implies
\be\label{dt2sa1di}\begin{split}\partial_t\norms{\partial_t^2s_a^1}{}{2}\lesssim \frac{\log(1+t)}{1+t^3} \end{split}\ee
yielding uniform-in-time boundedness of $\norms{\partial_t^2s_a^1}{}{2}$.
With the help of the uniform bounds on $\|\partial_t^js_a^1\|_2$, $j=1,2$, we can control $\Delta s_a^1$ and $\Delta \partial_ts_a^1$ using equations \eqref{sa1} and \eqref{Swdts} satisfied by $s_a^1$ and $\partial_ts_a^1$ respectively:
\begin{align}&\label{deltasa1}\norms{\Delta s_a^1}{}{2}\,\,\leq \underbrace{\norms{\partial_t s_a^1}{}{2}}_{\mathrm {unif. \,\,bounded}}+\quad\underbrace{\norms{V(s_a)}{}{2}}_{\mathclap{\hspace{1.5cm}\lesssim\frac{\log(1+t)}{(1+t^3)}\text{ by \eqref{vops}, \eqref{sa0sa1}}}}\,\,\lesssim \,\,\log(1+t),\\
&\label{deltadtsa1}\norms{\Delta \partial_ts_a^1}{}{2}\,\,\leq \underbrace{\norms{\partial_t^2 s_a^1}{}{2}}_{\mathrm {unif. \,\,bounded}}+\quad\underbrace{\norms{V(\partial_ts_a)}{}{2}+\norms{(\partial_tV)(s_a)}{}{2}}_{\lesssim\frac{\log(1+t)}{(1+t^3)}\text{ by \eqref{vops}, \eqref{sa0sa1}, \eqref{l2dtsa}}}\,\,\lesssim \,\,\log(1+t).\end{align}
Since we have $H^2\hookrightarrow H^{3/2}$, we obtain
\be\label{sa1h32}\boxed{\|\partial_t^js_a^1(t,\cdot)\|_{H^{3/2}}\lesssim\log(1+t)\quad\text{for }j=0,1.}\ee

Finally for estimating $\|\partial_t^js_e\|_{H^{3/2}}$ for $j=0,1$, again we will estimate $\partial_t^{j+1}s_e$ in $L^2$ and use the equations satisfied by $\partial_t^js_e$ to estimate $\Delta\partial_t^js_e$ and then the embedding $H^2\hookrightarrow H^{3/2}$. If we take derivatives of equations \eqref{se} and \eqref{p2}, we obtain the following equations to which we will apply energy estimates:
\be\label{SseWp2}\left.\begin{split} &\mathbf{S}(\partial_ts_e)=-(\partial_tV)(s_e)+(\partial_tm)\circ p_2+m\circ\partial_tp_2+(\partial_t\bar p_2)\circ m+\bar p_2\circ(\partial_t m)\\
&\mathbf{W}(\partial_t\bar p_2)=-\big[\partial_tg_\mathrm{pot}^\mathrm{T},\bar p_2\big]+\partial_t \mathop{M}_{\mathclap{\substack{\downarrow\\m\circ\overline{ s_a}-s_a\circ\overline{m}}}}+(\partial_tm)\circ\overline{ s_e}+m\circ \partial_t \overline{ s_e}-(\partial_ts_e)\circ\overline m-s_e\circ\partial_t\overline m\end{split}\right\}\ee
\be\label{dt2SseWp2}\left.\begin{split} &\mathbf{S}(\partial_t^2s_e)=-2(\partial_tV)(\partial_ts_e)-(\partial_t^2V)(s_e)+(\partial_t^2m)\circ p_2+\bar p_2\circ \partial_t^2m\\
&\hspace{1.5cm}\,+2\Big[(\partial_tm)\circ \partial_tp_2+(\partial_t\bar p_2)\circ \partial_tm\Big]+m\circ\partial_t^2p_2+(\partial_t^2\bar p_2)\circ m\\
&\mathbf{W}(\partial_t^2\bar p_2)=-\big[\partial_t^2g_\mathrm{pot}^\mathrm{T},\bar p_2\big]-2\big[\partial_tg_\mathrm{pot}^\mathrm{T},\partial_t\bar p_2\big]+\partial_t^2 M+(\partial_t^2m)\circ\overline{ s_e}-s_e\circ\partial_t^2\overline m\\
&\hspace{1.77cm}+2\Big[(\partial_tm)\circ\partial_t\overline{s_e}-(\partial_ts_e)\circ\partial_t\bar m\Big]+m\circ \partial_t^2 \overline{ s_e}-(\partial_t^2s_e)\circ\overline m\end{split}\qquad\,\right\}\ee
where $M:=m\circ\overline{ s_a}-s_a\circ\overline{m}$. Now we add the equations
\begin{subequations}\begin{align}\label{1}&\mathbf{W}\big((\partial_t^js_e)\circ\partial_t^j\bar s_e\big)=\mathbf{S}(\partial_t^js_e)\circ\partial_t^j\overline{s_e}-(\partial_t^js_e)\circ\overline{\mathbf{S}(\partial_t^js_e)}\nn\\
&\mathbf{W}((\partial_t^j\bar p_2)\circ\partial_t^j\bar p_2)=\mathbf{W}(\partial_t^j\bar p_2)\circ\partial_t^j\bar p_2+(\partial_t^j \bar p_2)\circ\mathbf{W}(\partial_t^j\bar p_2)\nn\end{align}\end{subequations}
side by side and then take traces to make  the following estimate:
\be\label{sevenseven}\partial_t\big(\overbrace{\norms{\partial_t^js_e}{2}{2}+\norms{\partial_t^jp_2}{2}{2}}^{=:E_j^2(t)}\big)\lesssim \|\mathbf{S}(\partial_t^js_e)\|_2\|\partial_t^js_e\|_2+\|\mathbf{W}(\partial_t^j\bar p_2)\|_2\|\partial_t^jp_2\|_2\quad\text{for }j=1,2.
\ee
We already know from \eqref{vops} that  $\norms{\partial_t^jV}{}{\mathrm{op}}\lesssim (1+t^3)^{-1}$ for $j=0,1,2$. Similarly
\be\label{dtjppt}\big\|[\partial_t^jg_\mathrm{pot}^\mathrm{T},(\cdot)]\|_\mathrm{op}\lesssim (1+t^3)^{-1}.\ee
Recalling $m(t,x,y)=-v_N(x-y)\phi(t,x)\phi(t,y)$, the definition of $M$ from \eqref{SseWp2}
and using estimates similar to the second one in \eqref{efVad} we obtain  
\be \label{dtjmM}\left.\begin{split}&\norms{(\partial_t^jm)\circ u}{}{2}\leq (1+t^3)^{-1}\norms{u}{}{2},\\& \|\partial_t^jM\|_2\lesssim (1+t^3)^{-1}\sum_{k=0}^j\|\partial_t^ks_a\|_2\end{split}\right\}\quad \text{for }j=1,2.\ee Considering all these estimates together with \eqref{SseWp2} implies
\begin{align}&\nn\|\mathbf{S}(\partial_ts_e)\|_2\lesssim \frac{1}{1+t^3}\Big(\overbrace{\|s_e\|_2+\|p_2\|_2}^{\mathclap{\substack{O(1)\text{ by \eqref{sep2}}}}}+\|\partial_tp_2\|_2\Big),\\
&\nn \|\mathbf{W}(\partial_t\bar p_2)\|_2\lesssim \frac{1}{1+t^3}\Big(\|p_2\|_2+\underbrace{\|s_a\|_2+\|\partial_ts_a\|_2}_{\mathclap{\substack{\qquad\lesssim \log(1+t) \text{ by \eqref{sa0sa1}, \eqref{l2dtsa}}}}}+\|s_e\|_2+\|\partial_ts_e\|_2\Big).  \end{align}
Inserting the above estimates in \eqref{sevenseven} for $j=1$, we obtain
\be\nn\partial_t\big(\overbrace{\norms{\partial_ts_e}{2}{2}+\norms{\partial_tp_2}{2}{2}}^{=:E_1^2(t)}\big)\lesssim \frac{1}{1+t^3}\Big(\overbrace{\norms{\partial_ts_e}{}{2}}^{\leq E_1(t)}+\overbrace{\norms{\partial_ts_e}{}{2}\norms{\partial_tp_2}{}{2}}^{\lesssim E_1^2(t)}+\big(1+\log(1+t)\big)\overbrace{\norms{\partial_tp_2}{}{2}}^{\leq E_1(t)}\Big)\\
\ee
from which it follows that
\[\partial_tE_1(t)\lesssim \frac{1}{1+t^3}E_1(t)+\frac{1+\log(1+t)}{1+t^3}.\]
This in turn implies that $E_1(t)$ is uniformly bounded in time. Using this, we can deduce 
\be\label{udtsep2}\|\partial_ts_e(t,x,y)\|_{L^2_{x,y}}\lesssim1\quad\text{and}\quad\|\partial_tp_2(t,x,y)\|_{L^2_{x,y}}\lesssim1.\ee 

Now considering estimates \eqref{vops}, \eqref{dtjppt}, \eqref{dtjmM} together with \eqref{dt2SseWp2} implies
\begin{align}&\nn\|\mathbf{S}(\partial_t^2s_e)\|_2\lesssim \frac{1}{1+t^3}\Big(\overbrace{\|s_e\|_2+\|p_2\|_2}^{\mathclap{\substack{O(1)\text{ by \eqref{sep2}}}}}+\overbrace{\|\partial_ts_e\|_2+\|\partial_tp_2\|_2}^{\mathclap{\substack{O(1)\text{ by \eqref{udtsep2}}}}}+\|\partial_t^2p_2\|_2\Big),\\
&\nn \|\mathbf{W}(\partial_t^2\bar p_2)\|_2\lesssim \frac{1}{1+t^3}\Big(\underbrace{\sum_{j=0}^1\big(\|\partial_t^js_e\|_2+\|\partial_t^jp_2\|_2\big)}_{O(1)\text{ as above}}\,\,\,\,+\,\,\,\,\underbrace{\sum_{j=0}^2\|\partial_t^js_a\|_2}_{\mathclap{\substack{\lesssim \log(1+t)\\ \text{by \eqref{sa0sa1}, \eqref{l2dtsa}, \eqref{dtsa0l2}, \eqref{dt2sa1di}}\\\text{and recalling }s_a=s_a^0+s_a^1}}}\,\,\,\,+\,\,\,\,\|\partial_t^2s_e\|_2\Big).  \end{align}
Inserting the above estimates in \eqref{sevenseven} for $j=2$, we obtain
\be\nn\partial_t\big(\overbrace{\norms{\partial_t^2s_e}{2}{2}+\norms{\partial_t^2p_2}{2}{2}}^{=:E_2^2(t)}\big)\lesssim \frac{1}{1+t^3}\Big(\overbrace{\norms{\partial_t^2s_e}{}{2}}^{\leq E_2(t)}+\overbrace{\norms{\partial_t^2s_e}{}{2}\norms{\partial_t^2p_2}{}{2}}^{\lesssim E_2^2(t)}+\big(1+\log(1+t)\big)\overbrace{\norms{\partial_t^2p_2}{}{2}}^{\leq E_2(t)}\Big)\\
\ee
which yields
\[\partial_tE_2(t)\lesssim \frac{1}{1+t^3}E_2(t)+\frac{1+\log(1+t)}{1+t^3}.\]
This implies that $E_2(t)$ is uniformly bounded in time, which helps us conclude
\be\label{udt2sep2}\|\partial_t^2s_e(t,x,y)\|_{L^2_{x,y}}\lesssim1\quad\text{and}\quad\|\partial_t^2p_2(t,x,y)\|_{L^2_{x,y}}\lesssim1.\ee
Now we can  estimate $\norms{\Delta \partial_t^js_e}{}{2}$, $j=0,1$ using \eqref{se} and the first equation in \eqref{SseWp2} as follows:
\begin{align}\label{deltase}&\norms{\Delta s_e}{}{2}\leq \norms{\partial_t s_e}{}{2}+\norms{V(s_e)}{}{2}+\norms{m\circ p_2}{}{2}+\norms{\bar p_2\circ m}{}{2}\lesssim 1 + \frac{1}{1+t^3}\\
\nn &\|\Delta \partial_ts_e\|_2\leq\|\partial_t^2s_e\|_2+\|(\partial_tm)\circ p_2\|_2+\|\bar p_2\circ\partial_tm\|_2+\|m\circ\partial_tp_2\|_2+\|(\partial_t\bar p_2)\circ m\|_2\\
&\label{deltadtse}\hspace{1.4cm}\lesssim 1 + \frac{1}{1+t^3}\end{align}
where we used \eqref{vops}, \eqref{dtjmM}, \eqref{sep2},  \eqref{udtsep2} and \eqref{udt2sep2}. The estimates above imply
\be\label{seh3h}\boxed{\|\partial_t^js_e(t,\cdot)\|_{H^{3/2}}\lesssim 1\quad\text{for }j=0,1.}\ee
due to the Sobolev embedding $H^{2}\hookrightarrow H^{3/2}$. Recalling $s_2=s_a^0+s_a^1+s_e$ and combining \eqref{sa0h32}, \eqref{sa1h32} and \eqref{seh3h} imply 
\[\|\partial_t^js_2\|_{H^{3/2}}\lesssim_\epsilon N^{\beta(1+\epsilon)}\log(1+t)\quad\text{for }j=0,1\]
which proves \eqref{s2hth}.

\textbf{Proof of \eqref{uhth}}. This is based on the identity
$s_2=2u\circ c=2\bar c\circ u$. We have
\be \label{uDs}D_x^{\sigma}u(t,x,y)=\frac{1}{2}\big(D_x^{\sigma}s_2\big)\circ c^{-1}\text{ and }D_y^{\sigma}u(t,x,y)=\frac{1}{2}\bar c^{-1}\circ D_y^{\sigma}s_2\ee
where $\sigma\in \mathbb{R}$ denotes the order of the derivative. \eqref{uDs} implies
\be\label{Dsigma}\begin{split}\|D^{\sigma}u(t,\cdot)\|_2^2&=\int\big(|\xi|^2+|\eta|^2\big)^{\sigma}|\hat u(t,\xi,\eta)|^2\mathrm{d}\xi\mathrm{d}\eta\\
&\lesssim \int |\xi|^{2\sigma}|\hat u(t,\xi,\eta)|^2\mathrm{d}\xi\mathrm{d}\eta+\int|\eta|^{2\sigma}|\hat u(t,\xi,\eta)|^2\mathrm{d}\xi\mathrm{d}\eta\\
&=\|D_x^{\sigma}u(t,\cdot)\|_2^2+\|D_y^{\sigma}u(t,\cdot)\|_2^2\lesssim \|s_2\|_{H^{\sigma}}^2 \end{split}\ee
where the last inequality follows from \eqref{uDs} since $\|c^{-1}\|_{\mathrm{op}}$ is uniformly bounded. 
Taking $\sigma=3/2$ in \eqref{uDs}-\eqref{Dsigma} proves \eqref{uhth}.
\newline

\textbf{Proof of \eqref{umlp}}. For $\tilde{\epsilon}>0$ small, we can make the following estimate:
\be\label{u2inf}\Big\|\|u(t,x,y)\|_{L^2(\mathrm{d}x)}\Big\|_{L^\infty(\mathrm{d}y)}\leq\Big\|\|u(t,x,y)\|_{L^\infty(\mathrm{d}y)}\Big\|_{L^2(\mathrm{d}x)}\lesssim_{\tilde{\epsilon}}\|D_y^{\frac{3}{2}+\tilde{\epsilon}}u(t,\cdot)\|_2\lesssim\|u\|_{H^{\frac{3}{2}+\tilde{\epsilon}}}\ee
where, for the second last inequality, we have used $H^s(\mathbb{R}^n)\hookrightarrow L^{\infty}(\mathbb{R}^n)$ for $s> n/2$ with $n=3$ (see e.g. Remark 1.4.1 (v) in \cite{Caz}).
Considering $\sigma=2$ in \eqref{uDs}-\eqref{Dsigma}, one can prove $\|u(t,\cdot)\|_{H^2}\lesssim \|s_2(t,\cdot)\|_{H^2}\lesssim N^{3\beta/2}\log(1+t)$ where the last inequality follows from \eqref{sa0h3h}, \eqref{deltasa1} and \eqref{deltase}. Interpolating between this $H^{2}$-norm estimate and the previously obtained $H^{3/2}$-norm estimate (see \eqref{uhth}) gives
\be\label{gsn}\|u(t,\cdot)\|_{H^{3/2+\tilde{\epsilon}}}\lesssim \underbrace{\bi(N^{\beta(1+\epsilon)})^{1-2\tilde{\epsilon}}\big(N^{3\beta/2}\big)^{2\tilde{\epsilon}}}
_{N^{\beta[1+\epsilon+\tilde{\epsilon}(1-2\epsilon)]}}\log(1+t).\ee
This last estimate  considered with \eqref{u2inf} proves \eqref{umlp}.\hfill$\Box$\newline

\noindent\textbf{Remarks.}
\textit{\begin{itemize}[leftmargin=.55cm]
\item[(i)] In the following section we will frequently use an estimate of $\big\|\|u(t,\cdot)\|_{2}\big\|_{4}:=\big\|\|u(t,x,y)\|_{L^2_x}\big\|_{L^4_y}$ to control most of the contributions in \eqref{F1}-\eqref{F4}. This follows by interpolation between $\big\|\|u\|_2\big\|_\infty$ and $\|u\|_2=\big\|\|u\|_2\big\|_2$ i.e. we have 
\be\label{u24}\big\|\|u(t,\cdot)\|_2\big\|_4\leq \big\|\|u(t,\cdot)\|_2\big\|_\infty^{1/2}\|u(t,\cdot)\|_2^{1/2}\lesssim N^{(\beta/2)(1+\epsilon)}\log(1+t)\text{, }\epsilon>0\ee
where for the last inequality we used  \eqref{pul2} and \eqref{umlp}. 
\item[(ii)] Recalling the relation $p\circ p+2p=\bar u\circ u$ and the fact that $(p\circ p)(t,x,x)\geq 0$ and also $p(t,x,x)\geq 0$, we have
$$\|p(t,x,y)\|_{L^2(\mathrm{d}x)}^2=(p\circ p)(t,y,y)\leq (\bar u\circ u)(t,y,y)=\|u(t,x,y)\|_{L^2(\mathrm{d}x)}^2$$
which implies (for any $\epsilon>0$)
\begin{align}&\label{p2i}\big\|\|p(t,\cdot)\|_2\big\|_\infty\leq \big\|\|u(t,\cdot)\|_2\big\|_\infty\lesssim N^{\beta(1+\epsilon)}\log(1+t)\\
&\label{p24}\big\|\|p(t,\cdot)\|_2\big\|_4\leq \big\|\|u(t,\cdot)\|_2\big\|_4\lesssim N^{(\beta/2)(1+\epsilon)}\log(1+t)
\end{align}
using \eqref{umlp} and \eqref{u24}.
\end{itemize}}
\begin{center} \textbf{4. The regular part of $|\tilde{\psi}\rangle$}
\end{center}

Our main result in  this section is the following:
\bt \label{rpfs} We have the following estimate for $|\tilde{\psi}^\mathrm{r}\rangle$ solving equation \eqref{9a}:
\be\||\tilde{\psi}^\mathrm{r}(t)\rangle\|_\mathbb{F}\lesssim_\epsilon N^{-1/2+\beta(1+\epsilon)}t\log^4(1+t)\ee
for any $\epsilon>0$.\et

We will need the following lemma for the proof of Theorem \ref{rpfs}: 
\bl\label{lemmafn} Given the definitions in \eqref{F1}-\eqref{F4} and \eqref{Fsingular}, the following estimates hold:
\be\label{Frest}\|\mathop{F_l^\mathrm{r}}^{\mathclap{\substack{\,\,\,\text{ without } r\\\,\,\,\text{ if }l=1\\\,\,\,\,\,\,\downarrow}}}(t)\|_{L^2(\mathbb{R}^{3l})}\lesssim_\epsilon \left\{\begin{split}&N^{-1/2+\beta(1+\epsilon)}\log^3(1+t)/(1+t^{3/2}),\quad l=1,3\\
&N^{-1+2\beta(1+\epsilon)}\log^4(1+t),\quad l=2,4.\end{split}\right.\ee
for any $\epsilon>0$.\el

\noindent\textit{Proof.} Let's prove \eqref{Frest} for $l=1,2$ first. We need to estimate the $L^2$-norms of the contributions in  \eqref{F1a}-\eqref{F1l} and the ones in \eqref{F2b}-\eqref{F2l}. Estimate for the term in \eqref{F2b} can be made as follows:
\begin{align}&\label{one} N^{-1}\|\int \mathrm{d}x_1\mathrm{d}x_2 v_N(x_1-x_2)p(x_2,y_1)u(x_2,y_2)(\bar u\circ u)(x_1,x_1)\|_{L^2(\mathrm{d}y_1\mathrm{d}y_2)}\\
&\nn\leq N^{-1}\int \mathrm{d}x_1\mathrm{d}x_2v_N(x_1-x_2)\|p(x_2,y_1)\|_{L^2(\mathrm{d}y_1)}\|u(x_2,y_2)\|_{L^2(\mathrm{d}y_2)}(\bar u\circ u)(x_1,x_1)\\
&\nn\leq N^{-1}\big\|\|p\|_2\big\|_{\infty}\big\|\|u\|_2\big\|_{\infty}\|v_N\|_1\underbrace{\|(\bar u\circ u)(x_2,x_2)\|_{L^1(\mathrm{d}x_2)}}_{\|u\|_2^2}\mathop{\,\lesssim_\epsilon\,}_{\mathclap{\substack{\\\text{by}\\\text{\eqref{pul2}},\\ \text{\eqref{umlp} and \eqref{p2i}  }}}}N^{-1+2\beta(1+\epsilon)}\log^4(1+t).\end{align}

Estimates of the terms in \eqref{F2c}-\eqref{F2e} are similar and differ slightly from \eqref{one}. We estimate only for \eqref{F2c}: 
\begin{align}
&\nn\leq N^{-1}\|\int \mathrm{d}x_1\mathrm{d}x_2 v_N(x_1-x_2)p(x_2,y_1)u(x_1,y_2)(\bar u\circ u)(x_1,x_2)\|_{L^2(\mathrm{d}y_1\mathrm{d}y_2)}\\
&\nn\leq N^{-1}\int \mathrm{d}x_1\mathrm{d}x_2 v_N(x_1-x_2)\|p(x_2,y_1)\|_{L^2(\mathrm{d}y_1)}\|u(x_1,y_2)\|_{L^2(\mathrm{d}y_2)}|(\bar u\circ u)(x_1,x_2)|\\
&\nn \leq N^{-1}\big\|\|p\|_2\big\|_{\infty}\big\|\|u\|_2\big\|_{\infty}\int\mathrm{d}x_2\, v_N(x_2)\underbrace{\Big(\int \mathrm{d}x_1\,|(\bar u\circ u)(x_1,x_1-x_2)|\Big)}_{\leq\|u\|_2^2\text{ uniformly in } x_2}\\
&\label{two}\leq N^{-1}\big\|\|p\|_2\big\|_{\infty}\big\|\|u\|_2\big\|_{\infty}\|v_N\|_1\|u\|_2^2\mathop{\,\lesssim_\epsilon\,}_{\mathclap{\substack{\text{by}\\\text{\eqref{pul2}},\\ \text{\eqref{umlp} and \eqref{p2i}  }}}}N^{-1+2\beta(1+\epsilon)}\log^4(1+t). \end{align}

Estimates of \eqref{F1a}-\eqref{F1b} are similar to \eqref{one} and the estimates of \eqref{F1c}-\eqref{F1f} are similar to \eqref{two}; the only difference being that, in \eqref{one}-\eqref{two}, we were able to pull two factors out of the integral in $L^\infty$-norm, each of which is either a $p$-term or a $u$-term whereas in estimates of \eqref{F1a}-\eqref{F1c} there is only one $u$ (or $p$)-term available for us to pull out in the same manner and we also need to pull $\phi$ out of the integral in $L^\infty$-norm. 
This explains the the difference between the powers of $N$ and the time dependence of the bounds in the estimates in \eqref{Frest}, in cases of $l=1$ and $l=2$.

Estimates of \eqref{F2f}-\eqref{F2g} are similar so let's just look at the estimate of \eqref{F2f}:
\begin{align}&\nn (1/2N)\|\int \mathrm{d}x_1\mathrm{d}x_2 v_N(x_1-x_2)u(y_1,x_1)u(x_2,y_2)\bar u(x_1,x_2)\|_{L^2(\mathrm{d}y_1\mathrm{d}y_2)}\\
&\nn\leq(1/2N)\int \mathrm{d}x_1\mathrm{d}x_2v_N(x_1-x_2)\|u(y_1,x_1)\|_{L^2(\mathrm{d}y_1)}\|u(x_2,y_2)\|_{L^2(\mathrm{d}y_2)}|u(x_1,x_2)|\\
&\nn \leq(1/2N)\big\|\|u\|_2\big\|_{\infty} \int\mathrm{d}x_2v_N(x_2)\underbrace{\Big(\int\mathrm{d}x_1\,|u(x_1,x_1-x_2)| \|u(y_1,x_1)\|_{L^2(\mathrm{d}y_1)}\Big)}_{\leq \|u\|_{H^{3/2+\tilde{\epsilon}}}\|u\|_2\text{ unif. in $x_2$}}\\
&\label{three}\leq(1/2N)\big\|\|u\|_2\big\|_{\infty}\|v_N\|_1\|u\|_{H^{3/2+\tilde{\epsilon}}}\|u\|_2\mathop{\,\lesssim_\epsilon\,}_{\mathclap{\substack{\text{by}\\\text{\eqref{pul2},}\\ \text{\eqref{gsn} and \eqref{umlp}}}}} N^{-1+2\beta(1+\epsilon)}\log^3(1+t).\end{align}

Estimate of the more singular term \eqref{F2i} differs slightly from the above estimate:
\begin{align}&\nn(1/2N)\|\int \mathrm{d}x_1 v_N(y_1-x_1)p(x_1,y_2)u(y_1,x_1)\|_{L^2(\mathrm{d}y_1\mathrm{d}y_2)}\\
&\nn \leq (1/2N)\Big\|\int \mathrm{d}x_1v_N(y_1-x_1)\|p(x_1,y_2)\|_{L^2(\mathrm{d}y_2)}u(y_1,x_1)\Big\|_{L^2(\mathrm{d}y_1)}\\
&\nn\leq (1/2N)\big\|\|p\|_2\big\|_{\infty}\int\mathrm{d}x_1\, v_N(x_1)\underbrace{\|u(y_1,y_1-x_1)\|_{L^2(\mathrm{d}y_1)}}_{\leq\|u\|_{H^{3/2+\tilde{\epsilon}}}\text{unif. in $x_1$}} \\
&\label{four}\leq (1/2N)\big\|\|p\|_2\big\|_{\infty}\|v_N\|_1\|u\|_{H^{3/2+\tilde{\epsilon}}}\mathop{\,\lesssim_\epsilon\,}_{\mathclap{\substack{\text{by}\\\text{\eqref{gsn} and \eqref{p2i}}}}} N^{-1+2\beta(1+\epsilon)}\log^2(1+t).\end{align}

Now let's consider the estimate of \eqref{F2h}:
\begin{align}&\nn N^{-1}\|\int \mathrm{d}x_1 v_N(x_1-y_1)u(y_1,y_2)(\bar u\circ u)(x_1,x_1)\|_{L^2(\mathrm{d}y_1\mathrm{d}y_2)}\\
&\nn \leq N^{-1}\Big \|\int \mathrm{d}x_1v_N(x_1-y_1)\|u(y_1,y_2)\|_{L^2(\mathrm{d}y_2)}(\bar u\circ u)(x_1,x_1)\Big\|_{L^2(\mathrm{d}y_1)}\\
&\nn\leq N^{-1}\Big\|\Big(v_N\ast\big((\bar u\circ u)(\cdot,\cdot)\big)\Big)(y_1)\|u(y_1,y_2)\|_{L^2(\mathrm{d}y_2)}\Big\|_{L^2(\mathrm{d}y_1)}\\
&\label{five}\leq N^{-1}\big\|\|u\|_2\big\|_{\infty}\|v_N\|_1\underbrace{\|(\bar u\circ u)(y_1,y_1)\|_{L^2(\mathrm{d}y_1)}}_{\big\|\|u\|_2\big\|_4^2}\mathop{\,\lesssim_\epsilon\,}_{\mathclap{\substack{\\\text{by}\\\text{\eqref{umlp}}\\ \text{and \eqref{u24}  }}}}N^{-1+2\beta(1+\epsilon)}\log^3(1+t).\end{align}

Estimates of \eqref{F2j}-\eqref{F2k} are similar and differ slightly from \eqref{five}. We will estimate for \eqref{F2j} in the following way:
\begin{align}&\nn N^{-1}\|\int \mathrm{d}x_1 v_N(x_1-y_1)u(x_1,y_2)(\bar u\circ u)(x_1,y_1)\|_{L^2(\mathrm{d}y_1\mathrm{d}y_2)}\\
&\nn\leq N^{-1}\Big\|\int\mathrm{d}x_1\, v_N(y_1-x_1)(\bar u\circ u)(x_1,y_1)\|u(x_1,y_2)\|_{L^2(\mathrm{d}y_2)}\Big\|_{L^2(\mathrm{d}y_1)}\\
&\nn\leq N^{-1}\big\|\|u\|_2\big\|_{\infty}\int\mathrm{d}x_1\,v_N(x_1)\underbrace{\|(\bar u\circ u)(y_1-x_1,y_1)\|_{L^2(\mathrm{d}y_1)}}_{\leq \big\|\|u\|_2\big\|_4^2\text{ uniformly in $x_1$}}\\
&\label{six} \leq N^{-1}\big\|\|u\|_2\big\|_{\infty}\|v_N\|_1\big\|\|u\|_2\big\|_4^2\mathop{\,\lesssim_\epsilon\,}_{\mathclap{\substack{\text{by}\\\text{\eqref{umlp} and \eqref{u24}  }}}}N^{-1+2\beta(1+\epsilon)}\log^3(1+t).
\end{align}
 
\eqref{F2l} is similar to the sum of the terms in \eqref{F2i} and \eqref{F2k} whose estimates have already been discussed.

Estimates of \eqref{F1g}-\eqref{F1h} are similar to \eqref{three}. The estimate of \eqref{F1i} resembles \eqref{four}. Estimate of \eqref{F1l} is similar to \eqref{five} and estimates of \eqref{F1j}-\eqref{F1k} resemble \eqref{six}. However, similar to the remarks coming right after \eqref{two}, in \eqref{F1g}-\eqref{F1l}, there is no $u$ (or $p$)-term
available for us to pull out of the integral in the way we did in \eqref{three}-\eqref{six}. Instead, we can pull $\phi$ out in $L^\infty$-norm, which explains the difference in the powers of $N$ and the time dependence of the bounds in \eqref{Frest}, in cases $l=1$, $l=2$. 

In order to prove \eqref{Frest} for $l=3,4$, we need to consider $L^2$-norms of the terms in \eqref{F3b}-\eqref{F3f} and the terms in \eqref{F4b}-\eqref{F4d}. Estimates of \eqref{F4b} and \eqref{F4c} are similar so let's make it for \eqref{F4b}:
\begin{align}&\nn(1/2N)\|\int \mathrm{d}x \bar p(y_2,x)v_N(y_1-x)u(x,y_4)u(y_3,y_1)\|_{L^2(\mathrm{d}y_1\mathrm{d}y_2\mathrm{d}y_3\mathrm{d}y_4)}\\
&\nn\leq(1/2N)\Big\|\int\mathrm{d}x\,v_N(y_1-x)\|p(x,y_2)\|_{L^2(\mathrm{d}y_2)}\|u(x,y_4)\|_{L^2(\mathrm{d}y_4)}\|u(y_3,y_1)\|_{L^2(\mathrm{d}y_3)}\Big\|_{L^2(\mathrm{d}y_1)}\\
&\nn\leq (1/2N)\big\|\|u\|_2\big\|_{\infty}^2 \|\big(v_N\ast \|p(\cdot,y_2)\|_{L^2(\mathrm{d}y_2)})\big)(y_1)\|_{L^2(\mathrm{d}y_1)}\\
&\nn\leq (1/2N)\big\|\|u\|_2\big\|_{\infty}^2\|v_N\|_1\|p\|_2\mathop{\,\lesssim_\epsilon\,}_{\mathclap{\substack{\text{by}\\\text{\eqref{pul2} and \eqref{umlp}}}}} N^{-1+2\beta(1+\epsilon)}\log^3(1+t).\end{align}
Estimates of \eqref{F3b}-\eqref{F3d} are similar but   we need to pull out $\big\|\|u(y_2,y_1)\|_{L^2_{y_2}}\big\|_{L^\infty_{y_1}}$ in \eqref{F3b}, $\big\|\|u(y_3,x)\|_{L^2_{y_3}}\big\|_{L^\infty_{x}} $ in \eqref{F3c}, $\big\|\|u(y_3,y_1)\|_{L^2_{y_3}}\big\|_{L^\infty_{y_1}}$ in \eqref{F3d} and also $\|\phi\|_\infty$ in all three of them, instead of the $\big\|\|u\|_2\big\|_\infty^2$ factor in the above estimate. That again causes the difference in the powers of $N$ and the time dependence of the bounds in \eqref{Frest} in cases $l=3,\,l=4$.

And our last estimate is for \eqref{F4d}:
\begin{align}&\nn  (1/2N)\|\int\mathrm{d}x_1\mathrm{d}x_2 \bar p(y_1,x_1)p(x_2,y_2)v_N(x_1-x_2)u(y_3,x_1)u(x_2,y_4)\|_{L^2(\mathrm{d}y_1\mathrm{d}y_2\mathrm{d}y_3\mathrm{d}y_4)}\\
&\nn\leq (1/2N) \int\mathrm{d}x_1\mathrm{d}x_2\, v_N(x_1-x_2) \|p(x_1,y_1)\|_{L^2(\mathrm{d}y_1)}\|p(x_2,y_2)\|_{L^2(\mathrm{d}y_2)} \|u(x_1,y_3)\|_{L^2(\mathrm{d}y_3)}\|u(x_2,y_4)\|_{L^2(\mathrm{d}y_4)}\\
&\nn\leq (1/2N)\big\|\|p\|_2\big\|_{\infty}\big\|\|u\|_2\big\|_{\infty}\int \mathrm{d}x_1\Big(v_N\ast\big(\|p(\cdot,y_2)\|_{L^2(\mathrm{d}y_2)}\|u(\cdot,y_4)\|_{L^2(\mathrm{d}y_4)}\big)\Big)(x_1)\\
&\nn\leq (1/2N)\big\|\|p\|_2\big\|_{\infty}\big\|\|u\|_2\big\|_{\infty}\|v_N\|_1\|p\|_2\|u\|_2\mathop{\,\lesssim_\epsilon\,}_{\mathclap{\substack{\text{by}\\\text{\eqref{pul2},}\\ \text{\eqref{umlp} and \eqref{p2i}}}}} N^{-1+2\beta(1+\epsilon)}\log^4(1+t).\end{align}
Estimates of \eqref{F3e}-\eqref{F3f} are similar but we need to pull out $\big\|\|u(x_2,y_3)\|_{L^2_{y_3}}\big\|_{L^\infty_{x_2}} $ in \eqref{F3e}, $\big\|\|p(x_2,y_2)\|_{L^2_{y_2}}\big\|_{L^\infty_{x_2}} $ in \eqref{F3f} and also $\|\phi\|_\infty$ in both of them. 
\hfill$\Box$\newline

\noindent\textit{Proof of Theorem \ref{rpfs}}. Recalling the equation \eqref{9a} satisfied by $|\tilde{\psi}^\mathrm{r}\rangle$ and the energy estimate \eqref{eefr1} obtained from it,  one can insert estimates in Lemma \ref{lemmafn} into the energy estimate \eqref{eefr1}  and this implies our claim in Theorem \ref{rpfs}.\hfill$\Box$

 \begin{center}\textbf{5. The singular part of $|\tilde{\psi}\rangle$ }
\end{center}

The singular part of $|\tilde{\psi}\rangle$, denoted by $|\tilde{\psi}^\mathrm{s}\rangle$, satisfies  equation \eqref{9b}. Let's recall from \eqref{12a}-\eqref{12b} how we split  $|\tilde{\psi}^{\mathrm{s}}\rangle$:
\begin{subequations}\begin{align}&\nn |\tilde{\psi}^\mathrm{s}\rangle=|\tilde{\psi}_1^\mathrm{a}\rangle+|\tilde{\psi}_1^\mathrm{e}\rangle\text{ where }\\
&\label{12an}\Big(\frac{1}{i}\partial_t-\int L(t,x,y)a_x^\ast a_y\,\mathrm{d}x\mathrm{d}y\Big)|\tilde{\psi}_1^\mathrm{a}\rangle=(0,0,F_2^\mathrm{s},F_3^\mathrm{s},F_4^\mathrm{s},0,\dots),\\
&\label{12bn}\Big(\frac{1}{i}\partial_t-\mathcal{L}\Big)|\tilde{\psi}_1^\mathrm{e}\rangle=-N^{-1/2}\mathcal{E}(t)|\tilde{\psi}_1^\mathrm{a}\rangle,\\&\nn |\tilde{\psi}_1^\mathrm{a}(0)\rangle=|\tilde{\psi}_1^\mathrm{e}(0)\rangle=0.\end{align}\end{subequations}
First we want to obtain estimates on the components of $|\tilde{\psi}_1^\mathrm{a}\rangle$ and then  use them to estimate the error part $|\tilde{\psi}_1^\mathrm{e}\rangle$. We would like to do the latter by applying an energy estimate to \eqref{12bn}. But the estimates we obtain for $|\tilde{\psi}_1^\mathrm{a}\rangle$ will still not ensure  sufficient $L^2$-integrability for the components of the forcing term in \eqref{12bn} as $N\rightarrow\infty$ and for $\beta$ close to $1/2$. Hence we will need to split  $|\tilde{\psi}_1^\mathrm{e}\rangle$ further into its regular and singular parts and we will repeat similar splitting procedure for some finitely many times before a final application of an energy estimate.

Recalling the explicit formula for $L(t,x,y)$ from \eqref{ltxy}, let's define $\tilde{V}(t,x,y)$ via the equation
\be \label{vt}L(t,x,y)= \Delta_x\delta(x-y)-\tilde{V}(t,x,y).
\ee
Let's also define the operator 
\be\label{sj}\bi{\mathrm{S}_j}=\frac{1}{i}\partial_t-\Delta_{\mathbb{R}^{3j}}+\underbrace{\sum_{k=1}^{j}\overbrace{\Big(\tilde{V}(t)\Big)_k}^{\mathclap{\substack{\text{action of }\\\tilde{V}(t)\text{ on a function}\\\text{in the}\\k\text{th variable} }}}}_{\mathlarger{\bi{\mathrm{V}_j}}}\ee
Hence we have the following set of equations being  equivalent to \eqref{12an}:
\begin{align}\label{psi1a}&\bi{\mathrm{S}_j}\psi_1^{(j)}=F_j^{\mathrm{s}}\text{ with }\psi_1^{(j)}(0)=0\text{ for }j=2,3,4\text{ and}\\&\nn |\tilde{\psi}_1^\mathrm{a}\rangle=(0,0,\psi_1^{(2)},\psi_1^{(3)},\psi_1^{(4)},0,\dots).\end{align}
Our main result in this section is the following:
\bt\label{estpsi1a}  We have the following estimates for $\psi_1^{(j)}$ satisfying \eqref{psi1a}: 
\begin{subequations}\begin{align}
&\label{p12}\|\psi_1^{(2)}\|_{L^2(\mathbb{R}^6)}\lesssim_\epsilon N^{-1+\beta+\beta\epsilon} t \log(1+t)\text{ for any }\epsilon>0,\\
&\label{p13}\|\psi_1^{(3)}\|_{L^2(\mathbb{R}^9)}\lesssim N^{\frac{-1+\beta}{2}},\\
&\label{p14}\|\psi_1^{(4)}\|_{L^2(\mathbb{R}^{12})}\lesssim_\epsilon N^{-1+\frac{3\beta}{2}+\beta\epsilon}\sqrt t \log^2(1+t)\text{ for any }\epsilon>0\end{align}\end{subequations}
which imply the following estimate for $|\tilde{\psi}_1^\mathrm{a}\rangle$ satisfying \eqref{12an}:
\be\label{ipp1a}\||\tilde{\psi}_1^\mathrm{a}\rangle\|_\mathbb{F}\lesssim N^{\frac{-1+\beta}{2}}t\log^2(1+t)\text{ for }\beta<1/2 .\ee
\et

We will need the following lemmas  to prove Theorem \ref{estpsi1a}:
\bl\label{CK} \text{\textnormal{(Christ-Kiselev Lemma, see e.g. Lemma 2.4 in \cite{Tao})}} Let $X,\, Y$ be Banach spaces, let I be a time interval, and let $K\in C^0(I\times I;B(X,Y))$ be a kernel taking values in the space of bounded operators from $X$ to $Y$. Suppose that $1\leq p<q\leq \infty$ is such that
\[\|\int_I K(t,s)f(s)\,\mathrm{d}s\|_{L^q(I,Y)}\lesssim \|f\|_{L^p(I,X)}\]
for all $f\in L^p(I,X)$. Then one also has 
\[\|\int_{s\in I:s<t}K(t,s)f(s)\,\mathrm{d}s\|_{L^q(I,Y)}\lesssim_{p,q}\|f\|_{L^p(I,X)}.\]
\el
{\color{white} new}

\bl\label{onoV}For the operator norm of $\bi{\mathrm{ V}_j}$ defined in \eqref{sj}, we have the following estimate:
\be\label{onoVest}\|\bi{\mathrm{ V}_j}\|_{L^2(\mathbb{R}^{3j})\rightarrow L^2(\mathbb{R}^{3j})}\leq j\|\tilde{V}(t)\|_{L^2(\mathbb{R}^{3})\rightarrow L^2(\mathbb{R}^{3})}\lesssim \frac{j}{1+t^3}\|u(t)\|_{L^2(\mathbb{R}^6)}^4\lesssim \frac{j\log^4(1+t)}{1+t^3}.\ee
\el

\noindent\textit{Proof}. The first inequality follows from the definition of $\bi{\mathrm{ V}_j}$ in \eqref{sj}. For the second inequality let's write $\tilde{V}(t)$ explicitly recalling \eqref{ltxy} and \eqref{vt}:
\begin{subequations}\begin{align}\nn(\tilde{V}(t)f)(x)&=\int\tilde{V}(t,x,y)f(y)\mathrm{d}y\\
&\label{svntna}=\big(v_N\ast |\phi|^2\big)(t,x)f(x)+\int v_N(x-y) \phi(t,x)\bar\phi(t,y)f(y)\mathrm{d}y\\
&\label{svntnb}\quad\,\, -\frac{1}{2}\Big(\big(\bar c)^{-1}\circ m\circ \bar u+u\circ\bar m\circ \big(\bar c)^{-1} +[\bi{\mathrm{W}}(\bar c),\big(\bar c)^{-1}]\Big)\circ f.\end{align}\end{subequations}
We can estimate $L^2$-norms of the   terms in \eqref{svntna} as:
\begin{align}
&\label{17a1}\|\big(v_N\ast |\phi |^2)f \|_2\leq \|v_N\|_1\|\phi\|_\infty^2\|f\|_2\lesssim \frac{1}{1+t^3}\|f\|_2,\\
&\label{17a2}\|\int v_N(x-y) \phi(t,x)\bar\phi(t,y)f(y)\mathrm{d}y\|_2\leq \|\phi\|_\infty^2\|v_N\ast f\|_2\lesssim \frac{1}{1+t^3}\|f\|_2.
\end{align}
where we used $\|\phi(t)\|_{L^\infty(\mathbb{R}^3)}\lesssim 1/(1+t^{3/2})$ from \eqref{finfty}. Similarly to \eqref{17a2}, one can prove for $m(t,x,y)=-v_N(x-y)\phi(t,x)\phi(t,y)$ that
\be\label{oy}\|m\circ l\|_{L^2(\mathbb{R}^6)}\lesssim\frac{1}{1+t^3}\|l\|_{L^2(\mathbb{R}^6)}\ee
for any $l\in L^2(\mathbb{R}^6)$.

Recalling the relation $\bar c^2=\bar c\circ \bar c=\delta(x-y)+u\circ \bar u$ and considering a contour $\Gamma$ enclosing the spectrum of the non-negative Hilbert-Schmidt operator $q:=u\circ \bar u$ one can write
\be\label{os}\bi{\mathrm{W}}(\bar c)=\bi{\mathrm{W}}(\sqrt{1+q})=\frac{1}{2\pi i}\int_\Gamma (q-z)^{-1}\underbrace{\bi{\mathrm{W}}(q)}_{\mathclap{\substack{\downarrow\\ m\circ\bar u\circ \bar c-u\circ c\circ\bar m\\\text{by (92) in \cite{GM}} }}}(q-z)^{-1}\sqrt{1+z}\,\mathrm{d}z\ee
Since $\big(\bar c\big)^{-1}$ and $(q-z)^{-1}$ have uniformly bounded operator norms and $|z|\lesssim \|u\|_2^2$, \eqref{oy} and \eqref{os} help us dominate $L^2$-norm of \eqref{svntnb} with
\[\frac{1}{1+t^3}\|u(t)\|_{L^2(\mathbb{R}^6)}^4\|f\|_{L^2(\mathbb{R}^3)}.\]
This last bound considered together with the estimates in \eqref{17a1}-\eqref{17a2}  proves the second inequality in \eqref{onoVest}. The last inequality in \eqref{onoVest} follows from the estimate $\|u(t)\|_{L^2(\mathbb{R}^6)}\lesssim \log(1+t)$ as we recall from \eqref{pul2}. \hfill$\Box$\newline

\noindent\textit{Proof of Theorem \ref{estpsi1a}}. \eqref{psi1a} is equivalent to the following set of equations: 
\begin{subequations}\begin{align}&\label{ps1js}\psi_1^{(j)}=\psi_{1,a}^{(j)}+\psi_{1,e}^{(j)}
\text{ where}\\&\label{ps1aj}\Big(\frac{1}{i}\partial_t-\Delta_{\mathbb{R}^{3j}}\Big)\psi_{1,a}^{(j)}=F_j^s ,\\&\label{ps1ej}\bi{\mathrm{S}_j}\psi_{1,e}^{(j)}
=-\bi{\mathrm{V}_j}\psi_{1,a}^{(j)},\\
&\nn\psi_{1,a}^{(j)}(0)=\psi_{1,e}^{(j)}(0)=0\text{ for }j=2,3,4.\end{align}\end{subequations}
We will try to obtain estimates on $\|\psi_{1,a}^{(j)}\|_{L^2(\mathbb{R}^{3j})}$  using an elliptic estimate in case of $j=2$ and for the cases $j=3,4$ we will make use of Strichartz estimates along with $TT^\ast$-method (to be explained shortly) and Christ-Kiselev Lemma (see Lemma \ref{CK}). Then we will use the following energy estimate to  control  $\psi_{1,e}^{(j)}$:
\begin{align}
\nn \partial_t\|\psi_{1,e}^{(j)}\|_{L^2(\mathbb{R}^{3j})}^2&\leq -2\mathrm{Im}\,\bi{\Big(}\big(\Delta_{\mathbb{R}^{3j}}-\bi{\mathrm{V}_j}\big)
\psi_{1,e}^{(j)}-\bi{\mathrm{V}_j}\psi_{1,a}^{(j)}\,\,\larger{\bi{,}}\,\,\psi_{1,e}^{(j)}\bi{\Big)}\\
&\nn\lesssim\|\bi{\mathrm{V}_j}\psi_{1,a}^{(j)}\|_{L^2(\mathbb{R}^{3j})}\|\psi_{1,e}^{(j)}\|_{L^2(\mathbb{R}^{3j})}\,\,{\scriptstyle\text{ since }\bi{\mathrm{V}_j}\text{ is self-adjoint}}\\
&\nn\lesssim \frac{j\log^4(1+t)}{1+t^3}\|\psi_{1,a}^{(j)}\|_{L^2(\mathbb{R}^{3j})}\|\psi_{1,e}^{(j)}\|_{L^2(\mathbb{R}^{3j})}\,{\scriptstyle\text{ by Lemma \ref{onoV}}}
\end{align}
which implies
\be\label{eeps1ej}\|\psi_{1,e}^{(j)}(t)\|_{L^2(\mathbb{R}^{3j})}\lesssim\int_0^t\frac{j\log^4(1+t_1)}{1+t_1^3}\|\psi_{1,a}^{(j)}(t_1)\|_{L^2(\mathbb{R}^{3j})}\,\mathrm{d}t_1.\ee \newline

\textbf{Case 1: $j=2$}. For $j=2$, recalling \eqref{F2s}, \eqref{ps1aj} becomes:
\begin{align}
&\label{tta}\Big(\frac{1}{i}\partial_t-\Delta_{\mathbb{R}^6}\Big)\psi_{1,a}^{(2)}=-\frac{1}{2N}v_N(y_1-y_2)\big\{\overbrace{u(t,y_1,y_2)+(\bar p\circ u)(t,y_1,y_2)}^{\bar c\circ u=u\circ c=\frac{1}{2}s_2}\big\}.
\end{align}
Solving \eqref{tta} by Duhamel's formula and using  integration by parts we get 
\begin{align}&\nn\|\psi_{1,a}^{(2)}(t,\cdot)\|_{L^2(\mathbb{R}^6)}\lesssim\left\|\int_0^te^{it_1(|\xi|^2+|\eta|^2)}\widehat{F^\mathrm{s}_2}(t_1,\xi,\eta)\mathrm{d}t_1\right\|_{L^2(\mathbb{R}^6)}\\
&\label{ps12duh}\lesssim \left\|\frac{\widehat{F^\mathrm{s}_2}(0,\xi,\eta)}{|\xi|^2+|\eta|^2}\right\|_{\mathrlap{L^2(\mathbb{R}^6)}}\qquad+\left\|\frac{\widehat{F^\mathrm{s}_2}(t,\xi,\eta)}{|\xi|^2+|\eta|^2}\right\|_{\mathrlap{L^2(\mathbb{R}^6)}}\qquad+\left\|\int_0^te^{it_1(|\xi|^2+|\eta|^2)}\frac{\partial_{t_1}\widehat{F^\mathrm{s}_2}(t_1,\xi,\eta)}{|\xi|^2+|\eta|^2}\mathrm{d}t_1\right\|_{\mathrlap{L^2(\mathbb{R}^6)}}.\end{align}
Now we need estimates of 
$$\left\|\frac{\partial_t^j\widehat{F^\mathrm{s}_2}(t,\xi,\eta)}{|\xi|^2+|\eta|^2}\right\|_{L^2(\mathbb{R}^6)}\quad\text{for }j=0,1.$$
Writing 
\be\nn \partial_t^jF_2^\mathrm{s}(t,x,y)=-\frac{1}{4N}v_N(x-y)\partial_t^js_2(t,x,y)=-\frac{1}{4N}\int\delta(x-y-z)v_N(z)\partial_t^js_2(t,x,y)\mathrm{d}z\ee
and considering the Fourier transform of $\delta(x-y-z)\partial_t^js_2(t,x,y)$ in the variables $x$, $y$:
\be\nn e^{iz\cdot\eta}\,\partial_t^j\widehat{s_2^z}(t,\xi+\eta)\quad\text{where}\quad s_2^z(t,x)=s_2(t,x,x-z)\ee
we can write 
\be\nn |\partial_t^j\widehat{F_2^\mathrm{s}}(t,\xi,\eta)|^2=\frac{1}{16N^2}\Big|\int v_N(z)e^{iz\cdot\eta}\,\partial_t^j\widehat{s_2^z}(t,\xi+\eta)\mathrm{d}z\Big|^2\lesssim\frac{ \|v\|_1}{N^2}\int |v_N(z)||\partial_t^j\widehat{s_2^\mathrm{z}}(t,\xi+\eta)|^2\mathrm{d}z.\ee
Hence after a change of variables
\begin{align}&\nn\left\|\frac{\partial_t^j\widehat{F^\mathrm{s}_2}(t,\xi,\eta)}{|\xi|^2+|\eta|^2}\right\|_{L^2(\mathbb{R}^6)}^2\lesssim\frac{1}{N^2}\int|v_N(z)|\frac{|\partial_t^j\widehat{s_2^z}(t,\xi)|^2}{\big(|\xi|^2+|\eta|^2\big)^2}\mathrm{d}\xi\mathrm{d}\eta\mathrm{dz}\\
&\nn\hspace{3.2cm}\,\lesssim \frac{1}{N^2}\int |v_N(z)|\Big(\underbrace{\int\frac{|\partial_t^j\widehat{s_2^z}(t,\xi)|^2}{|\xi|}\mathrm{d}\xi}_{\lesssim\|D^{-1/2}\partial_t^js_2^z\|_{L^2(\mathbb{R}^3)}^2}\Big)\mathrm{d}z\\
&\label{psi12n}\hspace{3.2cm}\,\,
\,\,\mathop{\lesssim}_{\mathclap{\substack{\text{by}\\\text{Trace theorem}}}} \,\,\frac{1}{N^2}\|\partial_t^js_2\|_{H^{\scriptsize{\frac{3}{2}}+}(\mathbb{R}^6)}^2.\end{align}
Now since $s_2=s_a^0+s_a^1+s_e$,  $\|\partial_t^js_2\|_{H^2}\lesssim N^{3\beta/2}\log(1+t)$ by \eqref{sa0h3h}, \eqref{deltasa1}-\eqref{deltadtsa1} and \eqref{deltase}-\eqref{deltadtse}. Interpolating this $H^2$-norm estimate with $\|\partial_t^js_2\|_{H^{3/2}}\lesssim N^{\beta(1+\epsilon)}\log(1+t)$ (see \eqref{s2hth}) and applying the resulting estimate in \eqref{psi12n} imply
\be\left\|\frac{\partial_t^j\widehat{F^\mathrm{s}_2}(t,\xi,\eta)}{|\xi|^2+|\eta|^2}\right\|_{L^2(\mathbb{R}^6)}\lesssim_\epsilon N^{-1+\beta+\beta\epsilon}\log(1+t)\quad\text{for any } \epsilon>0\text{ and for }j=0,1.\ee
This inserted in \eqref{ps12duh} implies
\be\label{psi1a2est}\|\psi_{1,a}^{(2)}(t)\|_{L^\infty((0,t);L^{2}({\mathbb{R}^6}))}\lesssim_\epsilon N^{-1+\beta+\beta\epsilon}t\log(1+t)\quad\text{for any }\epsilon>0.\ee
\eqref{psi1a2est} considered with \eqref{eeps1ej} for $j=2$ gives
\begin{align}
\nn \|\psi_{1,e}^{(2)}(t)\|_{L^2(\mathbb{R}^6)}\lesssim_\epsilon N^{-1+\beta+\beta\epsilon}\int_0^t 
\frac{t_1\log^5(1+t_1)}{1+t_1^3}\,\mathrm{d}t_1\lesssim N^{-1+\beta+\beta\epsilon}t\log(1+t).
\end{align}
Since $\psi_1^{(2)}=\psi_{1,a}^{(2)}+\psi_{1,e}^{(2)}$ as we  recall from \eqref{ps1js}, we can combine our last estimate with \eqref{psi1a2est} to obtain
\be\label{c1j2}\boxed{\|\psi_{1}^{(2)}(t)\|_{L^{2}({\mathbb{R}^6})}\lesssim_\epsilon N^{-1+\beta+\beta\epsilon}t\log(1+t)}\quad\text{for any }\epsilon>0.\ee\newline

\textbf{Case 2:} $j=3$. For $j=3$, recalling \eqref{F3s}, \eqref{ps1aj} becomes:
\begin{align} 
\label{thtoa1}\Big(\frac{1}{i}\partial_t-\Delta_{\mathbb{R}^9}\Big)\psi_{1,a}^{(3)}=-N^{-1/2}v_N(y_1-y_2)\phi(t,y_2)u(t,y_1,y_3).\end{align}
To put the forcing term in a more suitable form  for the mixed space-time norm estimates so that the power of $N$ will depend on $\beta$ in the desired way, we need the change of variables $x_1=y_1-y_2$ and $x_2=y_1+y_2$ which is inspired by the technique introduced in Lemma 4.6, \cite{CH} (see also the remark following  Lemma 5.3 in \cite{CH2}). So \eqref{thtoa1} takes the form
\begin{align}\label{smm}&\nn\Big(\frac{1}{i}\partial_t-2\big(\Delta_{x_1}+\Delta_{x_1}\big)-\Delta_{y_3}\Big)\psi_{1,a}^{(3)}(\mathsmaller{t,\frac{x_1+x_2}{2},\frac{x_2-x_1}{2}},y_3)\\
&=-N^{-1/2}v_N(x_1)\phi(\mathsmaller{t,\frac{x_2-x_1}{2}})u(\mathsmaller{t,\frac{x_1+x_2}{2},y_3}).\end{align}
Now if we consider the solution operator $T:=e^{it\{2(\Delta_{x_1}+\Delta_{x_2})+\Delta_{y_3}\}}$  for the corresponding free Schr\"{o}odinger equation, we have the following estimate: 
\begin{align}&\nn\|Tf_0\|_{L^2_tL^6_{x_1}L^2_{x_2y_3}}=\big\|\underbrace{\|Tf_0\|_{L^2_{x_2y_3}}}_{\mathclap{\qquad\,\,=\|e^{2it\Delta_{x_1}}f_0\|_{L^2_{x_2y_3}}}}\big\|_{L^2_tL^6_{x_1}}\leq \big\|\underbrace{\|e^{2it\Delta_{x_1}}f_0\|_{L^2_tL^6_{x_1}}}_{\mathclap{\substack{\lesssim \|f_0(x_1,x_2,y_3)\|_{L^2_{x_1}}\\\text{by Strichartz estimates}\\\text{in dimension 3}}}}\big\|_{L^2_{x_2y_3}}\lesssim\|f_0\|_{L^2(\mathbb{R}^9)}\end{align}
which proves
\be\label{t1} T: L^2(\mathbb{R}^9)\rightarrow L^2_tL^6_{x_1}L^2_{x_2y_3}.\ee
Similarly we also have 
\be \label{t2}T:L^2(\mathbb{R}^9)\rightarrow L^\infty_tL^2_{x_1x_2y_3}.\ee
If we consider 
\be (T^\ast f)(x_1,x_2,y_3)=\int_\mathbb{R} e^{-is\{2(\Delta_{x_1}+\Delta_{x_2})+\Delta_{y_3}\}}f(s, x_1,x_2,y_3)\,\mathrm{d}s,\ee
then \eqref{t1} is equivalent to 
\be\label{t4} T^\ast :L^2_tL^{6/5}_{x_1}L^2_{x_2y_3}\rightarrow L^2(\mathbb{R}^9).\ee
\eqref{t2} and \eqref{t4} imply
\be \label{t5}TT^\ast:L^2_tL^{6/5}_{x_1}L^2_{x_2,y_3}\rightarrow L^\infty_tL^2_{x_1,x_2,y_3}.\ee
Now using \eqref{t5} and Christ Kiselev Lemma (Lemma \ref{CK})  with $K(t,s)=e^{i(t-s)\{2(\Delta_{x_1}+\Delta_{x_2})+\Delta_{y_3}\}}$, $f$ being the right hand side of \eqref{smm}, $X=L^{6/5}(\mathbb{R}^3;L^2(\mathbb{R}^9))$, $Y=L^2(\mathbb{R}^9)$ and $p=2$, $q=\infty$, we obtain the first inequality in the following estimate:
\begin{align}\nn &\big\|\psi_{1,a}^{(3)}\|_{L^\infty((0,t);L^2(\mathbb{R}^9))}\lesssim N^{-1/2}\|v_N(x_1)\phi(\mathsmaller{t,\frac{x_2-x_1}{2}})u(\mathsmaller{t,\frac{x_1+x_2}{2},y_3})\big\|_{L^2_tL^{6/5}_{x_1}L^2_{x_2,y_3}}\\
&\nn\leq N^{-1/2}\Bigg(\int_0^t\|\phi(t_1)\|_{L^\infty(\mathbb{R}^3)}^2\bigg(\int v_N^{6/5}(x_1)\Big(\int|u(\mathsmaller{t_1,\frac{x_1+x_2}{2},y_3})|^2\,\mathrm{d}x_2\,\mathrm{d}y_3\Big)^{\mathlarger{\frac{1}{2}\cdot\frac{6}{5}}}\,\mathrm{d}x_1\bigg)^{\mathlarger{\frac{5}{6}}\cdot 2}\,\mathrm{d}t_1\Bigg)^\mathlarger{\frac{1}{2}}\\
&\nn\lesssim N^{-1/2}\|v_N\|_{L^{6/5}(\mathbb{R}^3)}\Big(\int_0^t \|\phi(t_1)\|_{L^\infty(\mathbb{R}^3)}^2\|u(t_1)\|_{L^2(\mathbb{R}^6)}^2\,\mathrm{d}t_1\Big)^{1/2}\\
&\label{psi13ae}\mathop{\lesssim}_{\mathclap{\substack{\,\quad\text{by \eqref{finfty}}\\\,\,\quad\text{and \eqref{pul2}}}}} N^{(-1+\beta)/2}\Big(\int_0^t\frac{\log^2(1+t_1)}{1+t_1^3}\mathrm{d}t_1\Big)^{1/2}\lesssim N^{(-1+\beta)/2}.
\end{align}
This inserted in \eqref{eeps1ej} for $j=3$ implies:
\begin{align}
\nn \|\psi_{1,e}^{(3)}(t)\|_{L^2(\mathbb{R}^9)}\lesssim N^{(-1+\beta)/2}\int_0^t\frac{\log^4(1+t_1)}{1+t_1^3}\,\mathrm{d}t_1\lesssim N^{(-1+\beta)/2}.
\end{align}
Combining the last estimate with \eqref{psi13ae} gives
\be\label{c2j3}\boxed{\|\psi_{1}^{(3)}\|_{L^2(\mathbb{R}^9)}\lesssim N^{(-1+\beta)/2}}
\ee
since $\psi_1^{(3)}=\psi_{1,a}^{(3)}+\psi_{1,e}^{(3)}$ by \eqref{ps1js} .\newline

\textbf{Case 3:} $j=4$. Finally for $j=4$, recalling \eqref{F4s}, \eqref{ps1aj} becomes:
\begin{align}\label{4efa}\Big(\frac{1}{i}\partial_t-\Delta_{\mathbb{R}^{12}}\Big)\psi_{1,a}^{(4)}=-\frac{1}{2N}v_N(y_1-y_2)u(t,y_3,y_1)u(t,y_2,y_4).\end{align}
Doing the same change of variables as before, i.e. $x_1=y_1-y_2$ and $x_2=y_1+y_2$ in \eqref{4efa} and letting $T$ denote the corresponding free propagator, this time we have $TT^\ast :L^2_tL^{6/5}_{x_1}L^2_{x_2,y_3,y_4}\rightarrow L^\infty_tL^2_{x_1,x_2,y_3,y_4}$. We again use Lemma \ref{CK} to obtain the first inequality in the following:
\begin{align}\nn&\|\psi_{1,a}^{(4)}\|_{L^\infty((0,t);L^2(\mathbb{R}^{12}))}\lesssim N^{-1}\Big\|v_N(x_1)u(\mathsmaller{t,y_3,\frac{x_1+x_2}{2}})u(\mathsmaller{t,\frac{x_2-x_1}{2},y_4})\Big\|_{L^2_tL^{6/5}_{x_1}L^2_{x_2,y_3,y_4}}\\
&\nn=N^{-1}\Bigg(\int_0^t\bigg(\int v_N^{6/5}(x_1)\Big(\underbrace{\int\,\big\|u(\mathsmaller{t_1,y_3,\frac{x_1+x_2}{2}})\big\|_{L^2_{y_3}}^2\big\|u(\mathsmaller{t_1,\frac{x_2-x_1}{2},y_4})\big\|_{L^2_{y_4}}^2\mathrm{d}x_2}_{\mathclap{\quad\qquad\qquad\qquad\qquad\lesssim \big\|\|u(t_1,x,y)\|_{L^2_x}\big\|_{L^4_y}^4\lesssim_\epsilon N^{2\beta(1+\epsilon)}\log^4(1+t_1)\text{ by \eqref{u24}}}} \Big)^{\mathlarger{\frac{1}{2}\cdot\frac{6}{5}}}\,\mathrm{d}x_1\bigg)^{\mathlarger{\frac{5}{6}}\cdot 2}\,\mathrm{d}t_1\Bigg)^{\mathlarger{\frac{1}{2}}}
\\& \label{psi1a4}\lesssim_{\epsilon} N^{-1+\beta(1+\epsilon)}\|v_N\|_{L^{6/5}}\bigg(\int_0^t\log^4(1+t_1)\,\mathrm{d}t_1\bigg)^{1/2}\lesssim N^{-1+\frac{3\beta}{2}+\beta\epsilon}t^{1/2}\log^2(1+t)\end{align}
for any $\epsilon<0$.  Inserting this in \eqref{eeps1ej} for $j=4$ implies 
\begin{align}
\nn \|\psi_{1,e}^{(4)}(t)\|_{L^2(\mathbb{R}^{12})}\lesssim_\epsilon N^{-1+\frac{3\beta}{2}+\beta\epsilon}\int_0^t \frac{t_1^{1/2}\log^6(1+t_1)}{1+t_1^3}\,\mathrm{d}t_1\lesssim N^{-1+\frac{3\beta}{2}+\beta\epsilon} \log^2(1+t)
\end{align}
which, when combined with \eqref{psi1a4}, gives
\be\label{c3j4}\boxed{\|\psi_{1}^{(4)}(t)\|_{L^2(\mathbb{R}^{12})}\lesssim_\epsilon N^{-1+\frac{3\beta}{2}+\beta\epsilon}t^{1/2} \log^2(1+t)\text{ for any }\epsilon>0}\ee
since $\psi_{1}^{(4)}=\psi_{1,a}^{(4)}+\psi_{1,e}^{(4)}$ by \eqref{ps1js}. \hfill$\Box$\newline

Theorem \ref{estpsi1a}, estimate \eqref{ipp1a} provides us with an estimate of $\big\||\tilde{\psi}_1^{\mathrm{a}}\rangle\|_\mathbb{F}$ in case of $\beta<1/2$, which decays as $N\rightarrow\infty$. We still need to estimate the error part $|\tilde{\psi}_1^{\mathrm{e}}\rangle$. 
Recalling \eqref{12an}-\eqref{12bn}, at this point, one might think of applying the standard $L^2$-energy estimate to \eqref{12bn} to obtain
\be\label{eef}\big\||\tilde{\psi}_1^{\mathrm{e}}(t)\rangle\big\|_\mathbb{F}\lesssim N^{-1/2}\int_0^t\big\|\mathcal{E}(t_1)|\tilde{\psi}_1^{\mathrm{a}}(t_1)\rangle\big\|_\mathbb{F}\,\mathrm{d}t_1\ee
in which we want to estimate the right hand side by using the estimates in Theorem \ref{estpsi1a}.
However, as we will explain shortly,  we will not be able to pick up the desired powers of $N$ from the estimate of $\big\|\mathcal{E}(t)|\tilde{\psi}_1^{\mathrm{a}}(t)\rangle\big\|_\mathbb{F}$
to ensure   a decay as $N\rightarrow\infty$ for $\beta<1/2$. This problem is due to the contribution to $N^{-1/2}\mathcal{E}(t)$ coming from the term \eqref{22y}, considered only with $\delta$-parts of $c(x,y)=\mathrm{ch}(k)(x,y)=\delta(x-y)+p(x,y)$ factors in it, namely
\be\frac{1}{2N}\int\mathrm{d}y_1\mathrm{d}y_2 v_N(y_1-y_2)\mathcal{Q}^\ast_{y_1y_2}\mathcal{Q}_{y_1y_2}.\ee
Notice that this corresponds to the potential part of the original Hamiltonian (see \eqref{originalham}-\eqref{11c}).  So let's define the Fock space operators
\begin{align} 
&\label{htilddef}\tilde{\mathbb{H}}:=\frac{1}{2N}\int\mathrm{d}y_1\mathrm{d}y_2 v_N(y_1-y_2)\mathcal{Q}^\ast_{y_1y_2}\mathcal{Q}_{y_1y_2},\\
& \label{hdef}\mathbb{H}:=N^{-1/2}\mathcal{E}(t)-\tilde{\mathbb{H}}.\end{align}
Then we can rewrite \eqref{eef} as
\begin{align}\label{mfi}\big\||\tilde{\psi}_1^{\mathrm{e}}(t)\rangle\big\|_\mathbb{F}\lesssim \int_0^t\bigg\{\big\|\mathbb{H}\underbrace{|\tilde{\psi}_1^{\mathrm{a}}(t_1)\rangle}_{\mathclap{\substack{\downarrow\\{\scriptsize{\overbrace{(0,0,\psi_1^{(2)}(t_1),\psi_1^{(3)}(t_1),\psi_1^{(4)}(t_1),0\dots)}}}\\\text{from \eqref{psi1a}}}}}\big\|_\mathbb{F}+\big\|\tilde{\mathbb{H}}\,|\tilde{\psi}_1^{\mathrm{a}}(t_1)\rangle\big\|_\mathbb{F}\bigg\} \,\mathrm{d}t_1.\end{align}

We need the following operator norm estimates on $\mathbb{H}$ and $\tilde{\mathbb{H}}$: 
\bl\label{hopest} Based on the definitions \eqref{htilddef}-\eqref{hdef}, we have the following estimates for the actions of $\mathbb{H}$ and $\tilde{\mathbb{H}}$ on the $j$th sector of Fock space:
\bs\begin{align}&\nn\label{hest}\|\mathbb{H}\psi^{(j)}\|_\mathbb{F}\lesssim_{\epsilon ,j}   \Big(N^{-1/2+\beta(1+\epsilon)} \log^{4}(1+t)
 + N^{-1+5\beta/2+\beta\epsilon}\log^{2}(1+t)\\
&\hspace{2.4cm}+N^{(-1+3\beta)/2}\frac{\log(1+t)}{1+t^{3/2}}\Big)\|\psi^{(j)}\|_{L^2(\mathbb{R}^{3j})}\\
&\label{htest}\|\tilde{\mathbb{H}}\psi^{(j)}\|_\mathbb{F}\lesssim N^{-1+3\beta}\|\psi^{(j)}\|_{L^2(\mathbb{R}^{3j})}\end{align}\es
for any $\psi^{(j)}\in L_\mathrm{s}^2(\mathbb{R}^{3j})$ and $\epsilon >0$.
\el

We prove Lemma \ref{hopest} in the Appendix. 

Now turning back to the energy estimate \eqref{mfi}, the inequalities given by \eqref{p12}-\eqref{p14} and \eqref{hest} imply that the first term inside the integral on the left hand side of \eqref{mfi} i.e. $\big\|\mathbb{H}|\tilde{\psi}_1^a\rangle\big\|_\mathbb{F}$ is of order 
$N^{(-1+\beta)/2}N^{(-1+3\beta)/2}$ for $\beta<1/2$ implying a decay as $N\rightarrow\infty$.  However, the second term $\big\|\tilde{\mathbb{H}}|\tilde{\psi}_1^a\rangle\big\|_\mathbb{F}$ is of order $N^{-1+3\beta}N^{(-1+\beta)/2}$  using \eqref{p13} and \eqref{htest}. In that case, we have a decay as $N\rightarrow\infty$ as long as we choose  $\beta<3/7$ which is not good enough but we can improve it as we will describe in the next section.  

\begin{center}\textbf{6. Iterating the splitting method}\end{center}

Let's recall how we split $|\tilde{\psi}\rangle$ which is defined by \eqref{saytilde} and satisfies equation \eqref{oe}. We first split $|\tilde{\psi}\rangle$  into its regular and singular parts as $|\tilde{\psi}^{\mathrm{r}}\rangle+|\tilde{\psi}^{\mathrm{s}}\rangle$ where $|\tilde{\psi}^{\mathrm{r}}\rangle$, $|\tilde{\psi}^{\mathrm{s}}\rangle$ satisfy equations \eqref{9a}-\eqref{9b} respectively. We obtained an estimate on $\||\tilde{\psi}^{\mathrm{r}}\rangle\|_\mathbb{F}$ in Theorem \ref{rpfs}. We then split $|\tilde{\psi}^{\mathrm{s}}\rangle$ into its approximate and error parts as $|\tilde{\psi}_1^{\mathrm{a}}\rangle+|\tilde{\psi}_1^{\mathrm{e}}\rangle$ where $|\tilde{\psi}_1^{\mathrm{a}}\rangle$, $|\tilde{\psi}_1^{\mathrm{e}}\rangle$ satisfy \eqref{12an}-\eqref{12bn} respectively. We obtained an estimate on $\||\tilde{\psi}_1^{\mathrm{a}}\rangle\|_\mathbb{F}$ in Theorem \ref{estpsi1a}. Theorems \ref{rpfs} and \ref{estpsi1a} not only provide with bounds that are slowly deteriorating in time but also imply a decay as $N\rightarrow\infty$ for $\beta<1/2$.  We then considered analyzing $|\tilde{\psi}_1^{\mathrm{e}}\rangle$ to see if we can extend these observations to the case of the full error $\||\psi_\mathrm{ex}\rangle-|\psi_\mathrm{ap}\rangle\|_\mathbb{F}=\||\tilde{\psi}\rangle\|_\mathbb{F}$ since $|\tilde{\psi}\rangle=|\tilde{\psi}^{\mathrm{r}}\rangle+|\tilde{\psi}_1^{\mathrm{a}}\rangle+|\tilde{\psi}_1^{\mathrm{e}}\rangle$. As we discussed at the end of the previous section, an approach  based solely on the energy estimate \eqref{eef}, which is rewritten in \eqref{mfi}, only provides with a bound which is meaningful as long as $\beta<3/7$. The problem is due to the term
$\tilde{\mathbb{H}}|\tilde{\psi}_1^\mathrm{a}\rangle$  on the right hand side of the equation for $|\tilde{\psi}_1^\mathrm{e}\rangle$:
\begin{align*}
\Big(\frac{1}{i}\partial_t-\mathcal{L}\Big)|\tilde{\psi}_1^\mathrm{e}\rangle=\overbrace{-\mathbb{H}\underbrace{|\tilde{\psi}_1^\mathrm{a}\rangle}_{\mathclap{\substack{\uparrow\\\scriptsize{\overbrace{{(0,0,\psi_1^{(2)},\psi_1^{(3)},\psi_1^{(4)},0,\dots)}}}\\\psi_1^{(j)}\text{ satisfy \eqref{psi1a} which is equivalent to \eqref{12an}}}}}-\tilde{\mathbb{H}}|\tilde{\psi}_1^\mathrm{a}\rangle}^{-N^{-1/2}\mathcal{E}(t)|\tilde{\psi}_1^\mathrm{a}\rangle\text{ by \eqref{hdef}}}
\end{align*}
 For an improvement, we now consider splitting $|\tilde{\psi}_1^{\mathrm{e}}\rangle$ into its regular and singular parts as $|\tilde{\psi}_1^{\mathrm{r}}\rangle+|\tilde{\psi}_1^{\mathrm{s}}\rangle$ where
\begin{align*}
&\Big(\frac{1}{i}\partial_t-\mathcal{L}\Big)|\tilde{\psi}_1^\mathrm{r}\rangle=-\mathbb{H}|\tilde{\psi}_1^\mathrm{a}\rangle\text{ with }|\tilde{\psi}_1^\mathrm{r}(0)\rangle=0,\\
&\Big(\frac{1}{i}\partial_t-\mathcal{L}\Big)|\tilde{\psi}_1^\mathrm{s}\rangle=-\tilde{\mathbb{H}}|\tilde{\psi}_1^\mathrm{a}\rangle\text{ with }|\tilde{\psi}_1^\mathrm{s}(0)\rangle=0
\end{align*}
and then we again split $|\tilde{\psi}_1^{\mathrm{s}}\rangle$ into its approximate and error parts as $|\tilde{\psi}_2^{\mathrm{a}}\rangle+|\tilde{\psi}_2^{\mathrm{e}}\rangle$ where
\begin{align*}
&\Big(\frac{1}{i}\partial_t-\int L(t,x,y)a_x^\ast a_y \mathrm{d}x\mathrm{d}y\Big)|\tilde{\psi}_2^\mathrm{a}\rangle=-\tilde{\mathbb{H}}|\tilde{\psi}_1^\mathrm{a}\rangle\text{ with }|\tilde{\psi}_2^\mathrm{a}(0)\rangle=0,\\
&\Big(\frac{1}{i}\partial_t-\mathcal{L}\Big)|\tilde{\psi}_2^{\mathrm{e}}\rangle=-N^{-1/2}\mathcal{E}(t)|\tilde{\psi}_2^\mathrm{a}\rangle=-\mathbb{H}|\tilde{\psi}_2^\mathrm{a}\rangle-\tilde{\mathbb{H}}|\tilde{\psi}_2^\mathrm{a}\rangle\text{ with }|\tilde{\psi}_2^{\mathrm{e}}(0)\rangle=0\\
&\text{where }|\tilde{\psi}_2^\mathrm{a}\rangle=(0,0,\psi_2^{(2)},\psi_2^{(3)},\psi_2^{(4)},0,\dots)\text{ and}\\&
\bi{\mathrm{S}_{l}}\psi_2^{(l)}\mathop{=}_{\mathclap{\substack{\mathlarger{\downarrow}\\[.05cm]\text{up to}\\\text{symmetrizations}}}}-\frac{1}{2N}v_N(y_1-y_2) \psi_1^{(l)}(t,y_1,y_2,\dots,y_l), \quad l=2,3,4\quad\text{ \scriptsize (recalling \eqref{sj}, \eqref{htilddef})}. 
\end{align*}
We will iterate splitting in this manner for $j-1$ times and at $j$th step we will only split into approximate and error parts as $|\tilde{\psi}_{j-1}^{\mathrm{s}}\rangle=|\tilde{\psi}_{j}^{\mathrm{a}}\rangle+|\tilde{\psi}_{j}^{\mathrm{e}}\rangle$ where $j$ is to be determined later. We can summarize our iteration scheme by the following set of equations:
\bs\begin{align}
&\label{52a} |\tilde{\psi}\rangle=|\tilde{\psi}^\mathrm{r}\rangle+|\tilde{\psi}_1^\mathrm{a}\rangle+\overbrace{|\tilde{\psi}_1^\mathrm{r}\rangle+
\overbrace{|\tilde{\psi}_2^\mathrm{a}\rangle+\cdots+\underbrace{|\tilde{\psi}_{j-1}^\mathrm{r}\rangle+\underbrace{|\tilde{\psi}_j^\mathrm{a}\rangle+|\tilde{\psi}_j^\mathrm{e}\rangle}_{|\tilde{\psi}_{j-1}^\mathrm{s}\rangle}}_{|\tilde{\psi}_{j-1}^\mathrm{e}\rangle}}^{|\tilde{\psi}_1^{\mathrm{s}}\rangle}}^{|\tilde{\psi}_1^\mathrm{e}\rangle}\text{ where}\\
&\label{52b} \Big(\frac{1}{i}\partial_t-\mathcal{L}\Big)|\tilde{\psi}^\mathrm{r}\rangle=(0,F_1,F_2^\mathrm{r},F_3^\mathrm{r},F_4^\mathrm{r},0,\dots)\text{ with }|\tilde{\psi}^\mathrm{r}(0)\rangle=0,\\
&\label{52c}\Big(\frac{1}{i}\partial_t-\int L(t,x,y)a_x^\ast a_y \mathrm{d}x\mathrm{d}y\Big)|\tilde{\psi}_1^\mathrm{a}\rangle=(0,0,F_2^\mathrm{s},F_3^\mathrm{s},F_4^\mathrm{s},0,\dots )\text{ with }|\tilde{\psi}_1^\mathrm{a}(0)\rangle=0 , \\
& \label{52d}\Big(\frac{1}{i}\partial_t-\mathcal{L}\Big)|\tilde{\psi}_1^\mathrm{r}\rangle=-\mathbb{H}|\tilde{\psi}_1^\mathrm{a}\rangle\text{ with }|\tilde{\psi}_1^\mathrm{r}(0)\rangle=0,\\
& \label{52e}\Big(\frac{1}{i}\partial_t-\int L(t,x,y)a_x^\ast a_y \mathrm{d}x\mathrm{d}y\Big)|\tilde{\psi}_2^\mathrm{a}\rangle=-\tilde{\mathbb{H}}|\tilde{\psi}_1^\mathrm{a}\rangle\text{ with }|\tilde{\psi}_2^\mathrm{a}(0)\rangle=0,
\\&\hspace{5.45cm}\nn\vdots\\
& \label{52f}\Big(\frac{1}{i}\partial_t-\mathcal{L}\Big)|\tilde{\psi}_{j-1}^\mathrm{r}\rangle=-\mathbb{H}|\tilde{\psi}_{j-1}^\mathrm{a}\rangle\text{ with }|\tilde{\psi}_{j-1}^\mathrm{r}(0)\rangle=0,\\
& \label{52g}\Big(\frac{1}{i}\partial_t-\int L(t,x,y)a_x^\ast a_y \mathrm{d}x\mathrm{d}y\Big)|\tilde{\psi}_j^\mathrm{a}\rangle=-\tilde{\mathbb{H}}
|\tilde{\psi}_{j-1}^\mathrm{a}\rangle\text{ with }|\tilde{\psi}_j^\mathrm{a}(0)\rangle=0,\\
&\label{52h}\Big(\frac{1}{i}\partial_t-\mathcal{L}\Big)|\tilde{\psi}_j^\mathrm{e}\rangle=-N^{-1/2}\mathcal{E}(t)|\tilde{\psi}_j^\mathrm{a}\rangle\text{ with }|\tilde{\psi}_j^\mathrm{e}(0)\rangle=0\text{ where}\\
&\label{52i}|\tilde{\psi}_j^\mathrm{a}\rangle:=(0,0,\psi_j^{(2)},\psi_j^{(3)},\psi_j^{(4)},0\,,\dots)\text{ and}\\
&\nn\bi{\mathrm{S}_{l}}\psi_j^{(l)}\mathop{\simeq_l}_{\mathclap{\substack{\mathlarger{\uparrow}\\[0.05cm]\text{means}\\\text{``equal up to}\\\text{symmetrizations"}}}}-\frac{1}{2N}v_N(y_1-y_2)\psi_{j-1}^{(l)}(t,y_1,y_2,\dots,y_l), \quad l=2,3,4.\quad\text{ \scriptsize (recalling \eqref{sj},  \eqref{htilddef}).} \end{align}\es

 We have the following result  on the inductive step of the iteration:
\bt\label{inductive} Under the above setting and based on the estimates in Theorem \ref{estpsi1a} and Lemma \ref{hopest}, we have the following estimates:
\bs\begin{align}
&\label{53a}\||\tilde{\psi}_j^\mathrm{r}(t)\rangle\|_\mathbb{F}\lesssim N^{j(-1+2\beta)}t^{(j+3)/2}\log^{6}(1+t),\\
&\label{53b} \|\psi_j^{(2)}(t)\|_{L^2(\mathbb{R}^{6})}\lesssim_{\epsilon } N^{(j-1)(-1+2\beta)}N^{-1+\beta+\beta\epsilon} t^{(j+1)/2}\log(1+t),\\
&\label{53c} \|\psi_j^{(3)}(t)\|_{L^2(\mathbb{R}^9)}\lesssim N^{(j-1)(-1+2\beta)}N^{(-1+\beta)/2}t^{(j-1)/2},\\
&\label{53d}\|\psi_j^{(4)}(t)\|_{L^2(\mathbb{R}^{12})}\lesssim_{\epsilon}N^{(j-1)(-1+2\beta)} N^{-1+\beta(3/2+\epsilon)}t^{j/2}\log^{2}(1+t)\end{align}\es
for all $j\geq 1$ and for every $\epsilon>0$.
\et
\noindent\textit{Proof.} Let's prove first prove \eqref{53b}-\eqref{53d}. The case $j=1$ for \eqref{53b}-\eqref{53b} was handled in Theorem  \ref{estpsi1a}. Hence, for the inductive step, assuming 
\eqref{53b}-\eqref{53d}, we will provide with a proof of
the case $j+1$. 

Now let's consider the equation \eqref{52g} by replacing $j$ with $j+1$. It will be equivalent to the following set of equations where we need to recall \eqref{sj}, \eqref{htilddef}:
\bs\begin{align}&\nn\bi{\mathrm{S}_{l}}\psi_{j+1}^{(l)}\simeq_l-\frac{1}{2N}v_N(y_1-y_2)\psi_j^{(l)}(t,y_1,y_2,\dots,y_l)\text{ with }\psi_{j+1}^{(l)}(0)=0\text{ for }l=2,3,4\text{ and }\\
&|\tilde{\psi}_{j+1}^\mathrm{a}\rangle=(0,0, \psi_{j+1}^{(2)}, \psi_{j+1}^{(3)}, \psi_{j+1}^{(4)},0, \dots).\end{align}\es
We can split $\psi_{j+1}^{(l)}$ similar to what we did in Theorem \ref{estpsi1a} as follows:
\begin{subequations}\begin{align}&\psi_{j+1}^{(l)}=\psi_{j+1,a}^{(l)}+\psi_{j+1,e}^{(l)}
\text{ where}\\&\label{56b}\Big(\frac{1}{i}\partial_t-\Delta_{\mathbb{R}^{3l}}\Big)\psi_{j+1,a}^{(l)}\simeq_l-\frac{1}{2N}v_N(y_1-y_2)\psi_j^{(l)}(t,y_1,y_2,\dots,y_l) ,\\&\bi{\mathrm{S}_{l}}\psi_{j+1,e}^{(l)}
=-\bi{\mathrm{V}_{l}}\psi_{j+1,a}^{(l)},\\
&\nn\psi_{j+1,a}^{(l)}(0)=\psi_{j+1,e}^{(l)}(0)=0\text{ for }l=2,3,4\end{align}\end{subequations}
and again after estimating $\|\psi_{j+1,a}^{(l)}\|_{L^2(\mathbb{R}^{3l})}$, we can use the energy estimate
\begin{align}\|\psi_{j+1,e}^{(l)}(t)\|_{L^2(\mathbb{R}^{3l})}&\nn\leq\int_0^t\|\bi{\mathrm{V}_{l}}(t_1)\psi_{j+1,a}^{(l)}(t_1)\|_{L^2(\mathbb{R}^{3l})}\,\mathrm{d}t_1\\
&\label{ee57}\mathop{\lesssim_l}_{\mathclap{\substack{\\\text{ by}\,\,\\\text{Lemma \ref{onoV}}}}}\int_0^t\frac{\log^4(1+t_1)}{1+t_1^3}\|\psi_{j+1,a}^{(l)}(t_1)\|_{L^2(\mathbb{R}^{3l})}\,\mathrm{d}t_1.\end{align} \newline
Hence let's prove the estimate on $\|\psi_{j+1,a}^{(l)}\|_{L^2(\mathbb{R}^{3l})}$ first. Similar to Case 2 in the proof of Theorem \ref{estpsi1a}, after a change of variables in equation \eqref{56b} and using Strichartz estimates, $TT^\ast$-method and Christ-Kiselev Lemma we can make the following estimate:
\begin{align}&\nn\|\psi_{j+1,a}^{(l)}(t)\|_{L^\infty((0,t);L^2(\mathbb{R}^{3l}))}\\
&\nn\lesssim_lN^{-1}\Bigg(\int_0^t\Big\|v_N(x_1)
\psi_{j}^{(l)}(\mathsmaller{t_1,\frac{x_1+x_2}{2},\frac{x_2-x_1}{2},y_3,\dots,y_l})
\Big\|_{L^{6/5}_{x_1}L^2_{x_2y_3\dots y_l}}^2\mathrm{d}t_1\Bigg)^{1/2}\\
\nn &=N^{-1}\Bigg(\int_0^t\bigg( \int \underbrace{v_N^{6/5}(x_1)}_{\substack{\\\downarrow\\ \\5/2}}\underbrace{\Big(\,\int |\psi_j^{(l)}(\mathsmaller{t_1,\frac{x_1+x_2}{2},\frac{x_2-x_1}{2}},y_3,\dots,y_l)|^2\,\mathrm{d}x_2\mathrm{d}y_3\dots\mathrm{d}y_l\Big)^\mathlarger{\frac{1}{2}\cdot\frac{6}{5}}}_{\substack{\\\downarrow\\\\5/3}}\,\mathrm{d}x_1\bigg)^{\mathlarger{\frac{5}{6}}\cdot 2}\,\mathrm{d}t_1\Bigg)^{1/2}\\
&\nn\mathop{\lesssim}_{\mathclap{\substack{\downarrow\\\text{H\"{o}lder}\\
\text{in }x_1\\\text{with}\\(5/2,5/3)}}
}\,N^{-1}\underbrace{\|v_N\|_{L^3(\mathbb{R}^3)}}_{\lesssim N^{2\beta}}\bigg(\int_0^t\|\psi_j^{(l)}(t_1)\|_{L^2(\mathbb{R}^{3l})}^2\,\mathrm{d}t_1\bigg)^{1/2}\quad\text{for }l=2,3,4.
\end{align}
Now inserting the bounds in \eqref{53b}-\eqref{53d} in the last line of the above estimate implies
\begin{align}&\nn \|\psi_{j+1,a}^{(l)}(t)\|_{L^\infty((0,t);L^2(\mathbb{R}^{3l}))}\mathop{\lesssim}_{\mathclap{\substack{\uparrow\\\text{constant}\\\text{involved}\\\text{depends on }\epsilon\text{ for }l=2,4}}} N^{j(-1+2\beta)}\,\mathlarger{\bi{\cdot}}\left\{\begin{array}{ll}
N^{-1+\beta(1+\epsilon)}t^{(j+2)/2}\log(1+t)\text{ for }l=2\\
N^{(-1+\beta)/2}t^{j/2}\text{ for }l=3\\
N^{-1+\beta(3/2+\epsilon)}t^{(j+1)/2}\log^2(1+t)\text{ for }l=4
\end{array}\right.\end{align}
Finally inserting this in \eqref{ee57} yields the same bounds for $\|\psi_{j+1,e}^{(l)}(t)\|_{L^\infty((0,t);L^2(\mathbb{R}^{3l}))}$ because $\log^4(1+t)/(1+t^3)$ inside the integral in line \eqref{ee57} is integrable. Since $\psi_{j+1}^{(l)}=\psi_{j+1,a}^{(l)}+\psi_{j+1,e}^{(l)}$, we completed the inductive step of proving \eqref{53b}-\eqref{53d}.

Now let's move on to proving \eqref{53a}. Replacing $j-1$ with $j$ in \eqref{52f}, applying the $L^2$-energy estimate to the resulting equation and using \eqref{53b}-\eqref{53d} we can make the following estimate for any $j\geq 1$:
\begin{align}&\big\||\tilde{\psi}_{j}^{\mathrm{r}}(t)\rangle\big\|_\mathbb{F}\nn\lesssim \int_0^t\big\|\mathbb{H}\overbrace{|\tilde{\psi}_j^{\mathrm{a}}(t_1)\rangle}^{\mathclap{\substack{{\scriptsize{\underbrace{(0,0,\psi_j^{(2)}(t_1),\psi_j^{(3)}(t_1),\psi_j^{(4)}(t_1),0,\dots)}}}\\\uparrow}}}\big\|_\mathbb{F} \,\mathrm{d}t_1\\
&\label{bc11}\mathop{\lesssim_{\epsilon}}
_{\mathclap{\substack{\text{by}\\\eqref{hest},\\\eqref{53b}-\eqref{53d}}}}
 \int_0^tN^{(-1+3\beta)/2}N^{(-1+\beta)/2+(j-1)(-1+2\beta)}t_1^{(j+1)/2}\log^{6}(1+t_1)\mathrm{d}t_1\end{align}
which implies \eqref{53a}.
\hfill$\Box$\newline

Now let's see what the energy estimate applied to \eqref{52h}
would imply if we were to stop the iteration at the $j$th step:
\bc\label{fee} For $|\tilde{\psi}_j^\mathrm{e}\rangle$ satisfying equation \eqref{52h} which is
\[\Big(\frac{1}{i}\partial_t-\mathcal{L}\Big)|\tilde{\psi}_j^\mathrm{e}\rangle=-N^{-1/2}\mathcal{E}(t)|\tilde{\psi}_j^\mathrm{a}\rangle\text{ with }|\tilde{\psi}_j^\mathrm{e}(0)\rangle=0\]
we have the following estimate
\begin{align}\label{feeest}\big\||\tilde{\psi}_j^\mathrm{e}(t)\rangle\big\|_\mathbb{F}\lesssim_j & \,\big(N^{j(-1+2\beta)}+N^{-1+3\beta+(j-1)(-1+2\beta)+(-1+\beta)/2}\big)t^{(j+3)/2}\log^{6}(1+t).
\end{align}
In particular, $\big\||\tilde{\psi}_j^\mathrm{e}(t)\rangle\big\|_\mathbb{F}=O(N^{(-3+7\beta)/2+(j-1)(-1+2\beta)})$ for $1/3\leq\beta< 1/2$. To ensure a decay we also need to choose 
\[\beta <\frac{1+2j}{3+4j}.\]
Hence, if j is sufficiently large, $\beta$ will be as close as desired to 1/2 in which case we will also have $\big\||\tilde{\psi}_j^\mathrm{e}(t)\rangle\big\|_\mathbb{F}$ decaying as $N\rightarrow\infty$. 
\ec

\noindent\textit{Proof}. Applying the standard energy estimate to \eqref{52h} gives
\begin{align}&\nn\||\tilde{\psi}_j^{\mathrm{e}}(t)\rangle\|_\mathbb{F}\lesssim \int_0^t \|N^{-1/2}\mathcal{E}(t_1)|\tilde{\psi}_j^\mathrm{a}(t_1)\rangle\|_\mathbb{F}\,\mathrm{d}t_1\\
&\lesssim \nn\int_0^t\bigg\{\underbrace{\|\mathbb{H}|\tilde{\psi}_j^{\mathrm{a}}(t_1)\rangle\big\|_\mathbb{F}}_{\substack{\lesssim N^{j(-1+2\beta)}t_1^{(j+1)/2}\log^{6}(1+t_1)\\
\text{ as in line \eqref{bc11}}}}+\hspace{1.5cm}\underbrace{\|\tilde{\mathbb{H}}|\tilde{\psi}_j^{\mathrm{a}}(t_1)\rangle\|_\mathbb{F}}_{\mathclap{\substack{\hspace{2.3cm}\lesssim N^{-1+3\beta}N^{(j-1)(-1+2\beta)}N^{(-1+\beta)/2}t_1^{(j+1)/2}\log^2(1+t_1)\\\hspace{.5 cm}\text{by \eqref{htest} and \eqref{53b}-\eqref{53d}}}}}\hspace{1.3cm}\bigg\} \,\mathrm{d}t_1
\end{align}
which implies estimate \eqref{feeest} .\hfill$\Box$\newline

\noindent\textit{Proof of Theorem \ref{mainresult} (Main result)}. Considering \eqref{52a}, Theorem \ref{rpfs},  Theorem \ref{inductive} and  Corollary \ref{fee}, we have 
\begin{align}&\nn\big\||\tilde{\psi}\rangle\big\|_\mathbb{F}\leq \,\overbrace{\big\||\tilde{\psi}^\mathrm{r}\rangle\big\|_\mathbb{F}}^{\mathclap{\substack{\hspace{.5cm} \text{by Theorem \ref{rpfs}}\\\hspace{2.2cm}\lesssim_\epsilon N^{-1/2+\beta(1+
\epsilon)}t\log^4(1+t)}}}
\hspace{2.6cm}+\hspace{.5cm}\overbrace{\sum_{m=1}^j\big\||\tilde{\psi}_m^\mathrm{a}\rangle\big\|_\mathbb{F}}^{\mathclap{\substack{\hspace{1.3cm}\text{by Theorem \ref{inductive}}\text{ for }\beta<1/2,\\\hspace{1.8cm}\lesssim N^{(-1+\beta)/2} t^{(j+1)/2} \log^2(1+t)}}}\hspace{1cm}\\&\hspace{1.5cm}\nn+\underbrace{\sum_{m=1}^{j-1}\big\||\tilde{\psi}_m^\mathrm{r}\rangle\big\|_\mathbb{F}}_{\mathclap{\substack{\hspace{1.5cm}\lesssim N^{-1+2\beta}t^{(j+3)/2}\log^6(1+t)\\\text{by Theorem \ref{inductive}}}}}\hspace{1.5cm}+\hspace{.5 cm}\underbrace{\big\||\tilde{\psi}_j^\mathrm{e}\rangle\big\|_\mathbb{F}}_{\mathclap{\substack{\hspace{4.3cm}\lesssim N^{(-3+7\beta)/2+(j-1)(-1+2\beta)}t^{(j+3)/2}\log^6(1+t)\\\hspace{2.7cm}\text{for }1/3\leq\beta<1/2\text{ by Corollary \ref{fee}}}}}\hspace{4cm}\\&\label{lastbetaest}
\lesssim_\epsilon \Big(N^{-1/2+\beta(1+\epsilon)}+N^{(-3+7\beta)/2+(j-1)(-1+2\beta)}\Big)t^{(j+3)/2}\log^6(1+t).\end{align}
For $1/3\leq\beta< 2j/(1-2\epsilon+4j)$, \eqref{lastbetaest} will decay as $N^{-1/2+\beta(1+\epsilon)}$ as $N\rightarrow\infty$ and for $2j/(1-2\epsilon+4j)\leq\beta<(1+2j)/(3+4j)$, \eqref{lastbetaest} will decay as $N^{(-3+7\beta)/2+(j-1)(-1+2\beta)}$, which implies estimate  \eqref{mrest}. \hfill$\Box$

\vspace{1cm}

\begin{center}\textbf{Appendix. Proof of Lemma \ref{hopest}} \end{center}
\noindent\textit{Proof of Lemma \ref{hopest}.} Recalling \begin{align}&\nn N^{-1/2}\mathcal{E}(t)=\sum_{j=1}^4\Big(\mathcal{E}_j(t)+\mathcal{E}_j^\ast (t)\Big)+\mathcal{E}_2^{\mathrm{sa}}(t)+\mathcal{E}_4^{\mathrm{sa}}(t) \text{ from } \eqref{ese}\\
&\nn \hspace{1.6cm} = \mathbb{H}+\underbrace{\frac{1}{2N}\int \mathrm{d}y_1\mathrm{d}y_2v_N(y_1-y_2)\mathcal{Q}^\ast_{y_1y_2}\mathcal{Q}_{y_1y_2}}_{\tilde{\mathbb{H}}}\text{ from }\eqref{hdef}-\eqref{htilddef},  \end{align}
it is sufficient to obtain operator norm estimates for the terms listed in \eqref{wad} since from the general theory of bounded linear operators on Hilbert spaces, the adjoint of an operator will have the same operator norm as the operator itself.
 
 A typical contribution to $\mathbb{H}$ coming from the contributions involved in the terms in \eqref{ese} is of the form
\[\int \mathrm{d}y_1\dots\mathrm{d}y_l f(y_1,\dots, y_l)\underbrace{\big(a, a^\ast\big)_l}_{\mathclap{\substack{l\text{th order}\\\text{term in }a,a^\ast}}}\text{ where }l=1,2,3,4.\]

Let's first consider estimating the second and the fourth order terms. 

\eqref{22s} and \eqref{22u} are similar terms. If we consider  \eqref{22s} in which we have $l=4$, $\big(a, a^\ast\big)_4=\mathcal{Q}_{y_1y_2}\mathcal{D}_{y_4y_3}$ and $f$ being equal to 
\[f_{\eqref{22s}}(y_1,y_2,y_3,y_4)=\frac{1}{2N}\int\mathrm{d}x_1\mathrm{d}x_2\big\{\bar u(y_1,x_1)\bar u(x_2,y_2)v_N(x_1-x_2)c(y_3,x_1)u(x_2,y_4)\big\},\]
we can write the contribution to $\mathbb{H}\psi^{(j)}$ coming from  \eqref{22s} as
\be\label{3aia}\int\mathrm{d}y_1\mathrm{d}y_2\mathrm{d}y_3\mathrm{d}y_4\big\{f_{\eqref{22s}}(y_1,y_2,y_3,y_4)\mathcal{Q}_{y_1y_2}\mathcal{D}_{y_4y_3}\big\}\big(\psi^{(j)}\big)\ee
producing a function in sector $j-2$ for $j\geq 2$, $L^2$-norm of which we want to estimate. We have the following typical estimates among others arising from symmetrizations  involved in the definition of the creation operators:
\begin{align}&\nn\text{Type 1: }\Big\|\int \Big(\int f_{\eqref{22s}}(y_1,y_2,y_3,y_2)\mathrm{d}y_2\Big)\psi^{(j)}(y_3,y_1, z_1,\dots,z_{j-2})\mathrm{d}y_1\mathrm{d}y_3\Big\|_{L^2(\mathbb{R}^{3(j-2)})}\\
&\nn\hspace{1.22cm}\leq \underbrace{\Big\|\int f_{\eqref{22s}}(y_1,y_2,y_3,y_2)\mathrm{d}y_2\Big\|_{L^2_{y_1y_3}}}_{\mathclap{\substack{\leq \text{sum of }L^2\text{-norms of}\\ \text{the terms like \eqref{F2b}, \eqref{F2h}}}}
}\|\psi^{(j)}\|_{L^2(\mathbb{R}^{3j})}\\
&\nn\hspace{1.22cm}\mathop{\lesssim_\epsilon}_{\mathclap{\substack{\text{by }\\\text{\eqref{Frest}, }l=2}}} N^{-1+2\beta(1+\epsilon)}\log^4(1+t)\|\psi^{(j)}\|_{L^2(\mathbb{R}^{3j})}\\
&\nn\text{Type 2: } \Big\|\int \Big(\int f_{\eqref{22s}}(y_1,y_2,y_3,y_1)\mathrm{d}y_1\Big)\psi^{(j)}(y_3,y_2, z_1,\dots,z_{j-2})\mathrm{d}y_1\mathrm{d}y_3\Big\|_{L^2(\mathbb{R}^{3(j-2)})}\\
&\nn\hspace{1.22cm}\leq \underbrace{\Big\|\int f_{\eqref{22s}}(y_1,y_2,y_3,y_1)\mathrm{d}y_1\Big\|_{L^2_{y_1y_3}}}_{\mathclap{\substack{\leq \text{sum of }L^2\text{-norms of}\\ \text{the terms like \eqref{F2c}, \eqref{F2j}}}}
}\|\psi^{(j)}\|_{L^2(\mathbb{R}^{3j})}\\
&\nn\hspace{1.22cm}\mathop{\lesssim_\epsilon}_{\mathclap{\substack{\text{by }\\\text{\eqref{Frest}, }l=2}}} N^{-1+2\beta(1+\epsilon)}\log^4(1+t)\|\psi^{(j)}\|_{L^2(\mathbb{R}^{3j})}\\
&\nn\text{Type 3: }\Big\|\int \mathrm{d}y_1\mathrm{d}y_2\mathrm{d}y_3\, f_{\eqref{22s}}(y_1,y_2,y_3,z_1)\,\psi^{(j)}(y_3,y_2,y_1, z_2,\dots,z_{j-2})\Big\|_{L^2(\mathbb{R}^{3(j-2)})}\\
&\nn\hspace{1.22cm}\leq \|\underbrace{f_\eqref{22s}}_{\mathclap{\substack{\text{sum of}\\\text{\eqref{F4c}-\eqref{F4d} like terms}}}}\|_{L^2(\mathbb{R}^{12)}}\|\psi^{(j)}\|_{L^2(\mathbb{R}^{3j})}\mathop{\lesssim_\epsilon}_{\mathclap{\substack{\text{by }\\\text{\eqref{Frest}, }l=4}}} N^{-1+2\beta(1+\epsilon)}\log^4(1+t)\|\psi^{(j)}\|_{L^2(\mathbb{R}^{3j})}.
 \end{align}
With the above estimates we can estimate the contribution  \eqref{3aia} as:
\begin{align}\nn&\Big\|\int\mathrm{d}y_1\mathrm{d}y_2\mathrm{d}y_3\mathrm{d}y_4\big\{f_{\eqref{22s}}(y_1,y_2,y_3,y_4)\mathcal{Q}_{y_1y_1}\mathcal{D}_{y_4y_3}\big\}\big(\psi^{(j)}\big)\Big\|_{L^2(\mathbb{R}^{3(j-2)})}\\&\nn\lesssim_{j,\epsilon}N^{-1+2\beta(1+\epsilon)}\log^4(1+t)\|\psi^{(j)}\|_{L^2(\mathbb{R}^{3j})}.\end{align} 
 
If we consider the contribution to  $\mathbb{H}\psi^{(j)}$ coming from \eqref{22v} and its adjoint, we have 
\begin{align}&\nn \int \mathrm{d}y_1\mathrm{d}y_2\mathrm{d}y_3\mathrm{d}y_4\Big\{f_{\eqref{22v}}(y_1,y_2,y_3,y_4)\overbrace{\mathcal{Q}_{y_1y_2}\mathcal{Q}_{y_3y_4}}^{\mathclap{\substack{j\geq 4\text{ should hold}\\\text{for non-trivial}\\\text{contribution}}}}+\overline{f_{\eqref{22v}}}(y_1,y_2,y_3,y_4)\mathcal{Q}_{y_1y_2}^\ast\mathcal{Q}_{y_3y_4}^\ast\Big\}\big(\psi^{(j)}\big)\\
&\text{with }f_{\eqref{22v}}(y_1,y_2,y_3,y_4)=\frac{1}{2N}\int\mathrm{d}x_1\mathrm{d}x_2\big\{\bar u(y_1,x_1)\bar u(x_2,y_2)v_N(x_1-x_2)c(y_3,x_1)\bar c(x_2,y_4)\big\}
\end{align}
which will produce a contribution to $\mathbb{H}\psi^{(j)}$ of the following type:
\begin{align}\nn\Big(0,\dots,0,\int_{\mathbb{R}^{12}} \mathrm{d}\bi{\mathrm{y}}\big\{f_{\eqref{22v}}(\bi{\mathrm{y}})\psi^{(j)}(\bi{\mathrm{y}}, z_1,\dots,z_{j-4})\big\},0,\dots,0,\big(\overline{f_{\eqref{22v}}}\otimes\psi^{(j)}\big)(z_1,\dots,z_{j+4}),0,\dots\Big)\end{align}
Fock space norm of which is 
\begin{align}&\nn\lesssim_j\underbrace{\|f_{\eqref{22v}}\|_{L^2(\mathbb{R}^{12})}}_{\mathclap{\substack{\leq\text{sum of }L^2\text{-norms of}\\\text{ terms like \eqref{F4a}-\eqref{F4d}}}}}\|\psi^{(j)}\|_{L^2(\mathbb{R}^{3j})}\\
&\nn\lesssim_\epsilon \Big(N^{-1+2\beta(1+\epsilon)}\log^4(1+t)+N^{-1+5\beta/2+\beta\epsilon}\log^2(1+t)\Big)\|\psi^{(j)}\|_{L^2(\mathbb{R}^{3j})}\end{align}
where the last inequality follows by \eqref{Frest}, $l=4$ and also by the following estimate (see \eqref{F4s} for $F^\mathrm{s}_4$):
\begin{align}\|F_4^\mathrm{s}\|_{L^2(\mathbb{R}^{12})}&\nn\lesssim N^{-1}\Big(\int v_N^2(y_1-y_2)\|u(y_3,y_1)\|^2_{L^2_{y_3}}\|u(y_2,y_4)\|^2_{L^2_{y_4}}\mathrm{d}y_1\mathrm{d}y_2\Big)^{1/2}\\
&\nn \lesssim N^{-1}\big\|\|u(y_2,y_4)\|_{L^2_{y_4}}\big\|_{L^\infty_{y_2}}\big\|\big(v_N^2\ast\|u(y_3,\cdot)\|^2_{L^2_{y_3}}\big)(y_2)\big\|_{L^1_{y_2}}^{1/2}\\
&\nn\lesssim_\epsilon N^{-1+\beta(1+\epsilon)}\log(1+t) \|v_N\|_{L^2(\mathbb{R}^3)}\|u\|_{L^2(\mathbb{R}^6)}\scriptsize{\text{ by \eqref{umlp}}}\\
&\label{F4L2}\lesssim N^{-1+5\beta/2+\beta\epsilon}\log^2(1+t).\end{align}

Now let's look at the contribution coming from \eqref{22w} only in the most singular case which corresponds to keeping only the $\delta$-parts of $c$-terms recalling $c(x,y)=\delta(x-y)+p(x,y)$:
\begin{align}&\label{22wms}\int \mathrm{d}y_1\mathrm{d}y_2\mathrm{d}y_4\Big\{\Big(\underbrace{\frac{1}{2N}\int\mathrm{d}x_2\bar u(x_2,y_2)v_N(y_1-x_2)u(x_2,y_4)}_{=:f(y_1,y_2,y_4)}\Big)\mathcal{D}_{y_1y_2}\mathcal{D}_{y_4y_1}\Big\}\big(\psi^{(j)}\big).
\end{align} 
This will not cause any sector shifts. We have the following typical estimates among others arising from symmetrization:
\begin{align}&\nn\text{Type 1: } \Big\|\Big(\int \mathrm{d}y_2 f(z_1,y_2,y_2)\Big)\psi^{(j)}(z_1,\dots,z_j)\Big\|_{L^2(\mathbb{R}^{3j})}\\
&\nn\hspace{1.22cm}=\frac{1}{2N}\Big\|\Big(\int \mathrm{d}x_2(u\circ\bar u)(x_2,x_2)v_N(z_1-x_2)\Big)\psi^{(j)}(z_1,\dots,z_j)\Big\|_{L^2(\mathbb{R}^{3j})}\\
&\nn\hspace{1.22cm}\lesssim N^{-1}\Big\|\Big(v_N\ast(u\circ \bar u)(\cdot,\cdot)\Big)(z_1)\Big\|_{L^\infty_{z_1}}\|\psi^{(j)}\|_{L^2(\mathbb{R}^{3j})}\\
&\nn\hspace{1.22cm}\leq N^{-1} \|v_N\|_{L^1(\mathbb{R}^3)}\big\|\|u(x,z_1)\|_{L^2_x}\big\|_{L^\infty_{z_1}}^2\|\psi^{(j)}\|_{L^2(\mathbb{R}^{3j})}\\
&\nn\hspace{1.22cm}\lesssim_\epsilon N^{-1+2\beta(1+\epsilon)}\log^2(1+t)\|\psi^{(j)}\|_{L^2(\mathbb{R}^{3j})}\quad\scriptstyle{\text{ by \eqref{umlp}}}\\
&\nn\text{Type 2: }\Big\|\int f(z_1,y_2,z_2)\psi^{(j)}(z_1,y_2,z_3,\dots,z_j)\mathrm{d}y_2\Big\|_{L^2(\mathbb{R}^{3j})}\\
&\nn\hspace{1.22cm}=\frac{1}{2N}\Big\|\int\Big(\int\bar u(x_2,y_2)v_N(z_1-x_2)u(x_2,z_2)\mathrm{d}x_2\Big)\psi^{(j)}(z_1,y_2,z_3,\dots,z_j)\mathrm{d}y_2\Big\|_{L^2(\mathbb{R}^{3j})}\\
&\nn\hspace{1.22cm}\lesssim N^{-1}\Big\|\|u(x_2,z_2)\|_{L^2_{z_2}}\Big\|_{L^\infty_{x_2}}\Big\|\int v_N(x_2)\Big(\underbrace{\int \bar u(z_1-x_2,y_2)\psi^{(j)}(z_1,y_2,z_3,\dots,z_j)\mathrm{d}y_2}_{\leq\|u(z_1-x_2,y_2)\|_{L^2_{y_2}}\|\psi^{(j)}(z_1,y_2,z_3,\dots,z_j)\|_{L^2_{y_2}}}\Big)\mathrm{d}x_2\Big\|_{L^2(\mathbb{R}^{3(j-1)})}\\
&\nn\hspace{1.22cm}\leq N^{-1}\|v_N\|_{L^1(\mathbb{R}^3)}\big\|\|u(x_2,z)\|_{L^2_{z}}\big\|_{L^\infty_{x_2}}^2\|\psi^{(j)}\|_{L^2(\mathbb{R}^{3j})}\\
&\nn\hspace{1.22cm}\lesssim_\epsilon N^{-1+2\beta(1+\epsilon)}\log^2(1+t)\|\psi^{(j)}\|_{L^2(\mathbb{R}^{3j})}\quad\scriptstyle{\text{ by \eqref{umlp}}}.
\end{align}
Estimate for the contribution coming from \eqref{22x} is almost the same with the above and \eqref{22q} can be estimated similarly. The other $\mathcal{D}\mathcal{D}$-contribution comes from \eqref{22z} but this term is even less singular due to not having any $c$-factors.  

Contributions to $\mathbb{H}\psi^{(j)}$ coming from \eqref{22t} and \eqref{22r} are similar hence if we look at the contribution from \eqref{22r}, considered only in the most singular case corresponding to keeping only the $\delta$-parts of $c$-terms, it has the form
\begin{align}\frac{1}{2N}\int\mathrm{d}y_1\mathrm{d}y_2\mathrm{d}y_4 \big\{\bar u(y_4,y_2)v_N(y_1-y_4)\mathcal{D}_{y_1y_2}\mathcal{Q}_{y_1y_4}\big\}\big(\psi^{(j)}\big)\end{align}
lowering the sector by two. We can make the following typical estimate for this contribution up to symmetrizations:
\begin{align}&\nn\frac{1}{2N}\Big\|\int v_N(z_1-y_4)\bar u(y_4,y_2)\psi^{(j)}(y_2,y_4,z_1,\dots,z_{j-2})\mathrm{d}y_2\mathrm{d}y_4\Big\|_{L^2(\mathbb{R}^{3(j-2)})}\\
&\nn\lesssim N^{-1}\Big\|\big\|v_N(z_1-y_4)\|u(y_4,y_2)\|_{L^2_{y_2}}\big\|_{L^2_{y_4}}\|\psi^{(j)}(y_2,y_4,z_1,\dots,z_{j-2})\|_{L^2_{y_2y_4}}\Big\|_{L^2_{z_1\dots z_{j-2}}}\\
&\nn\leq N^{-1}\|v_N\|_{L^2(\mathbb{R}^3)}\big\|\|u(y_2,y_4)\|_{L^2_{y_2}}\big\|_{L^\infty_{y_4}}\|\psi^{(j)}\|_{L^2(\mathbb{R}^{3j})}\\
&\nn\lesssim_\epsilon N^{-1+5\beta/2+\beta\epsilon}\log(1+t)\|\psi^{(j)}\|_{L^2(\mathbb{R}^{3j})}\scriptsize{\text{ by \eqref{umlp}}}
\end{align}

Similar estimates can be made for the contributions to $\mathbb{H}\psi^{(j)}$ coming from the term in \eqref{22y} provided we keep the $p$-part of $\bar c(y_1,x_1)$ (or of $c(y_3,x_1)$) and replace the remaining three $c$-factors with their corresponding $\delta$-parts. 

We move on to checking the second order contributions to  $\mathbb{H}\psi^{(j)}$.

\eqref{22c} and \eqref{22d} are similar terms. \eqref{22h} seems to be more singular compared to \eqref{22g}. So let's estimate the contributions to $\mathbb{H}\psi^{(j)}$ coming from  \eqref{22d} and \eqref{22h} which can be considered together in the form:
\begin{align}
&\hspace{2.5cm}\int\mathrm{d}y_1\,\mathrm{d}y_2\big\{f(y_1,y_2)\mathcal{D}_{y_1y_2}\big\}\big(\psi^{(j)}\big)\quad\text{where}
\\ &\nn f(y_1,y_2)=\frac{1}{2N}\int\mathrm{d}x_1\mathrm{d}x_2 \big\{v_N(x_1-x_2) \big[
\nn(\bar u\circ\bar c)(x_1,x_2) u(y_1,x_1)\bar c(x_2,y_2)\\
&\nn\hspace{6.1cm}+2(u\circ\bar u)(x_1,x_1)\bar u(y_2,x_2)u(y_1,x_2)\big]\big\}
\end{align}
and we can estimate it as follows:
\begin{align}\nn&\Big\|\int\mathrm{d}y_1\,\mathrm{d}y_2\big\{f(y_1,y_2)\mathcal{D}_{y_1y_2}\big\}\big(\psi^{(j)}\big)\Big\|_{L^2(\mathbb{R}^{3j})}\\
&\nn\leq \sum_{k=1}^j\Big\|\int f(z_k,y_2)\psi^{(j)}(y_2,\overbrace{z_1,\dots, z_j}^{z_k \text{ missing}})\mathrm{d}y_2 \Big\|_{L^2(\mathbb{R}^{3j})}\\
&\lesssim_{j}\underbrace{\|f\|_{L^2(\mathbb{R}^6)}}_{\mathclap{\substack{\qquad\leq\text{sum of }L^2\text{-norms}\\\hspace{1cm}\text{of terms similar to \eqref{F2b},}\\\hspace{1cm}\text{\eqref{F2d},\eqref{F2f}, \eqref{F2g} }}}}\|\psi^{(j)}\|_{L^2(\mathbb{R}^{3j})}\mathop{\lesssim_{j,\epsilon}}_{\mathclap{\substack{\text{by}\\\text{ \eqref{Frest}, }l=2}}}N^{-1+2\beta(1+\epsilon)}\log^4(1+t)\|\psi^{(j)}\|_{L^2(\mathbb{R}^{3j})}.\end{align}

If we consider the contribution to $\mathbb{H}\psi^{(j)}$ coming from \eqref{22e}-\eqref{22f} and their adjoints, we have 
\begin{align}& \hspace{1.2cm}\int \mathrm{d}y_1\mathrm{d}y_2 \big\{f(y_1,y_2)\mathcal{Q}_{y_1y_2}+\bar {f}(y_1,y_2)\mathcal{Q}^\ast_{y_1y_2}\big\}\big(\psi^{(j)}\big)\quad\text{where}\\
&\hspace{2cm}\nn f(y_1,y_2)=f_{\eqref{22e}}(y_1,y_2)+f_\eqref{22f}(y_1,y_2)\quad\text{and}\\
&\nn f_{\eqref{22e}}(y_1,y_2)=\frac{1}{2N}\int\mathrm{d}x_1\mathrm{d}x_2(\bar u\circ\bar c)(x_1,x_2)v_N(x_1-x_2)c(y_1,x_1)\bar c(x_2,y_2),\\
&\nn f_{\eqref{22f}}(y_1,y_2)=\frac{1}{2N}\int\mathrm{d}x_1\mathrm{d}x_2( u\circ c)(x_1,x_2)v_N(x_1-x_2)\bar u(y_1,x_1)\bar u(x_2,y_2)\end{align}
which will produce a contribution to $\mathbb{H}\psi^{(j)}$ of the following form
\begin{align}\nn \Big(0,\dots,0, \int_{\mathbb{R}^6}\mathrm{d}\textbf{y}\big\{f(\textbf{y})\psi^{(j)}(\textbf{y},z_1,\dots,z_{j-2})\big\},0,0,0, (\bar f\otimes\psi^{(j)})(z_1,\dots,z_{j+2}),0,\dots \Big)\end{align}
Fock space norm of which is 
\begin{align}&\nn\lesssim_j\Big(\underbrace{\|f_\eqref{22e}\|_{L^2(\mathbb{R}^6)}+\|f_\eqref{22f}\|_{L^2(\mathbb{R}^6)}}_{\mathclap{\substack{\leq\text{sum of }L^2\text{-norms of terms like}\\\eqref{F2a}, \eqref{F2e}, \eqref{F2g}, \eqref{F2l}}}}\Big)\|\psi^{(j)}\|_{L^2(\mathbb{R}^{3j})}\\
&\nn\lesssim_\epsilon \big(N^{-1+2\beta(1+\epsilon)}\log^4(1+t)+N^{-1+5\beta/2+\beta\epsilon}\log^2(1+t)\big)\|\psi^{(j)}\|_{L^2(\mathbb{R}^{3j})}\end{align}
where the last inequality follows by \eqref{Frest}, $l=2$ and the following estimate (see \eqref{F2s} for $F^\mathrm{s}_2$):
\begin{align}
 \|F_2^\mathrm{s}\|_{L^2(\mathbb{R}^6)}&\lesssim\nn N^{-1}\bigg(\int v_N^2(y_2)\Big(\int\big\{|u(y_1,y_1-y_2)|^2+|(\bar p\circ u)(y_1,y_1-y_2)|^2\big\}\mathrm{d}y_1\Big)\mathrm{d}y_2\bigg)^{\mathlarger{\frac{1}{2}}}\\
&\nn \lesssim N^{-1} \underbrace{\|v_N\|_{L^2(\mathbb{R}^3)}}_{O(N^{3\beta/2})}\big(\|u\|_{H^{\frac{3}{2}+}(\mathbb{R}^6)}+\big\|\|u(x,y)\|_{L^2_x}\big\|^2_{L^4_y}\big)\\
&\label{F2L2}\lesssim_\epsilon N^{-1+5\beta/2+\beta\epsilon} \log^2(1+t)\scriptsize{\text{ by \eqref{gsn}, \eqref{u24}, \eqref{p24}}}.\end{align}

Next let's deal with the third order terms. \eqref{22i}, \eqref{22j}, \eqref{22k} are providing $\mathcal{D}a$ (or $a^\ast\mathcal{Q}$)-terms which lower the sector by one. The most singular contribution comes from \eqref{22k}. Let's consider its estimate  in the most singular case by keeping the $\delta$-parts of $c$-factors. The corresponding contribution to $\mathbb{H}\psi^{(j)}$ will have the following form:
\begin{align}
\frac{1}{\sqrt{N}}\int \mathrm{d}y_1\mathrm{d}y_3\big\{
v_N(y_1-y_3)\bar \phi(y_3)a_{y_1}^\ast a_{y_1}a_{y_3}\big\}\big(\psi^{(j)}\big)\end{align}
whose $L^2$-norm is
\begin{align}&\nn \simeq_j\frac{1}{\sqrt{N}} \Big\|\int\mathrm{d}y_3 \big\{v_N(z_1-y_3)\bar \phi(y_3)\psi^{(j)}(y_3, z_1,\dots,z_{j-1})\big\}\Big\|_{L^2(\mathbb{R}^{3(j-1)})}\\&\nn \leq \frac{\|\phi(t,\cdot )\|_{L^\infty(\mathbb{R}^3)}}{\sqrt{N}}\|v_N\|_{L^2(\mathbb{R}^3)}\big\|\|\psi^{(j)}(y_3,z_1\dots,z_{j-1})\|_{L^2_{y_3}}\big\|_{L^2_{z_1\dots z_{j-1}}}\\
&\nn\lesssim\frac{N^{(-1+3\beta)/2}}{1+t^{3/2}}\|\psi^{(j)}\|_{L^2(\mathbb{R}^{3j})}.\end{align}

We can write the contributions to $\mathbb{H}\psi^{(j)}$ coming from \eqref{22l} and \eqref{22o} together with their adjoints  in the form:
\begin{align}
&\label{19jncont}\int \mathrm{d}y_1\mathrm{d}y_2\mathrm{d}y_3\big\{f(y_1,y_2,y_3)\mathcal{Q}_{y_1y_2}a_{y_3}+\bar f(y_1,y_2,y_3)a_{y_1}^\ast\mathcal{Q}^\ast_{y_2y_3}\big\}\big(\psi^{(j)}\big)\quad\text{where}\\
&\nn f(y_1,y_2,y_3)=N^{-1/2}\int\mathrm{d}x_1\mathrm{d}x_2v_N(x_1-x_2)\big\{\bar u(y_1,x_1)\phi(x_2)\bar u(x_2,y_2)c(y_3,x_1)\\
&\nn\hspace{6.8cm}+\bar u(y_1,x_1)\bar\phi(x_2)c(y_2,x_1)\bar c(x_2,y_3)\big\}
\end{align}
which will produce a contribution of the following form:
\begin{align}\nn\Big(0,\dots ,0, \int_{\mathbb{R}^9} \mathrm{d}\textbf{y}\big\{f(\textbf{y})\psi^{(j)}(\textbf{y}, z_1,\dots ,z_{j-3})\big\},0,\dots,0,\big(\bar f\otimes\psi^{(j)}\big)(z_1,\dots,z_{j+3}),0,\dots\Big)\end{align}
Fock space norm of which is 
\begin{align}&\nn\lesssim_j\underbrace{\|f\|_{L^2(\mathbb{R}^9)}}_{\mathclap{\substack{\qquad\quad\lesssim\text{sum of }L^2\text{-norms of}\\\qquad\quad\text{tems like \eqref{F3a}-\eqref{F3f}}}}}\|\psi^{(j)}\|_{L^2(\mathbb{R}^{3j})}\\
&\nn\lesssim_\epsilon \Big(N^{-1/2+\beta(1+\epsilon)}\log^3(1+t)/(1+t^{3/2})+N^{(-1+3\beta)/2}\log(1+t)/(1+t^{3/2})\Big)\|\psi^{(j)}\|_{L^2(\mathbb{R}^{3j})}\end{align}
by \eqref{Frest}, $l=3$ and the following estimate (see \eqref{F3s} for $F_3^\mathrm{s}$):
\begin{align}
\nn\|F_3^\mathrm{s}\|_{L^2(\mathbb{R}^9)}&\lesssim N^{-1/2}\|\phi\|_{L^\infty(\mathbb{R}^3)}\Big(\int v_N^2(y_1-y_2)\|u(y_3,y_1)\|_{L^2_{y_3}}^2\mathrm{d}y_1\mathrm{d}y_2\Big)^{1/2}\\
&\nn\lesssim \frac{N^{-1/2}}{1+t^{3/2}}\|v_N\|_{L^2(\mathbb{R}^3)}\|u\|_{L^2(\mathbb{R}^6)}\\
&\label{F3L2}\lesssim N^{(-1+3\beta)/2}\log(1+t)/(1+t^{3/2})\quad\scriptsize{\text{by \eqref{pul2}}}.
\end{align}
Other third order contributions to $\mathbb{H}\psi^{(j)}$ are less singular and can be estimated similarly. The first order contributions in \eqref{22a}-\eqref{22b}  are providing with similar bounds and the estimates for them are similar to  the previous estimates. The estimates so far prove \eqref{hest}.

Finally let's prove the estimate \eqref{htest} on $\tilde{\mathbb{H}}\psi^{(j)}$. This is the contribution coming from \eqref{22y} when all $c$-factors are replaced with their corresponding $\delta$-parts as we can recall from the definition \eqref{htilddef}. We have the following estimate:
\begin{align}&\nn\frac{1}{2N}\Big\|\int \mathrm{d}y_1\mathrm{d}y_2\big\{v_N(y_1-y_2)\mathcal{Q}^\ast_{y_1y_2}\mathcal{Q}_{y_1y_2}\big\}\big(\psi^{(j)}\big)\Big\|_{L^2(\mathbb{R}^{3j})}\\
&\nn\simeq_jN^{-1}\|v_N(z_1-z_2)\psi^{(j)}(z_1,z_2,\dots,z_j)\|_{L^2(\mathbb{R}^{3j})}\\
&\nn\lesssim N^{-1} \|v_N\|_{L^\infty(\mathbb{R}^3)}\|\psi^{(j)}\|_{L^2(\mathbb{R}^{3j})}\lesssim N^{-1+3\beta}\|\psi^{(j)}\|_{L^2(\mathbb{R}^{3j})}\end{align}
with which we completed proving Lemma \ref{hopest}. \hfill$\Box$
\ssection{\large Acknowledgements}

I am grateful to Professor Manoussos Grillakis and Professor Matei Machedon for their valuable advices, encouragement and discussions related to this work.
\vspace{0.09in}
\begin{center}
\textbf{\large References}
\end{center}
{\def\section*#1{} 

}

\end{document}